%% file: pap_rev.tex
\documentclass[final,leqno,onefignum,onetabnum]{siamltex}

\usepackage{fullpage}
\usepackage{comment}
\usepackage{upgreek}
\usepackage{cite}
\usepackage[mathscr]{eucal}
\usepackage{color}
\usepackage[usenames,dvipsnames,svgnames,table]{xcolor}
\usepackage{amsmath,amssymb,amsopn,mathtools}
\usepackage{graphicx}
\usepackage{hyperref}       
\hypersetup{
	colorlinks=true,   
	citecolor=blue,     
	filecolor=blue,     
	linkcolor=red,    
	urlcolor=blue}      

\usepackage{relsize}
\numberwithin{equation}{section}
\usepackage{enumitem}

\newlist{Assumption}{enumerate}{1}
\setlist[Assumption]{label=A\arabic*}

\definecolor{Blue}{rgb}{0,0,1}
\definecolor{Red}{rgb}{1,0,0}
\definecolor{Green}{rgb}{0,1,0}
\definecolor{Cyan}{rgb}{0,0.72,0.92}
\definecolor{Amethyst}{rgb}{0.6,0.4,0.8}
\definecolor{Bronze}{rgb}{0.8,0.5,0.2}
\definecolor{Violet}{rgb}{0.54,0.17,0.89}

\input{commands}

\usepackage{amsmath,amssymb}
\usepackage{epsfig, algorithmic, algorithm}
\usepackage{multirow, booktabs}
\usepackage{siunitx}
\usepackage{xr}
\mathtoolsset{showonlyrefs=true}

\newlist{steps}{enumerate}{1}
\setlist[steps, 1]{label = Step \arabic*:}

\title{Reduced order models for Lagrangian Hydrodynamics}

\author{
  Dylan Matthew Copeland\thanks{Center for Applied Scientific Computing,
  Lawrence Livermore National Laboratory, Livermore, CA 94550
  (copeland11@llnl.gov)   }
  \and
  Siu Wun Cheung\thanks{Center for Applied Scientific Computing, Lawrence
  Livermore National Laboratory, Livermore, CA 94550 (cheung26@llnl.gov)}
  \and
  Kevin Huynh\thanks{Applications, Simulations, and Quality, Lawrence Livermore
  National Laboratory, Livermore, CA 94550 (huynh24@llnl.gov)}
  \and
  Youngsoo Choi\thanks{Center for Applied Scientific Computing, Lawrence
  Livermore National Laboratory, Livermore, CA 94550 (choi15@llnl.gov)}
}

\begin{document}
\setlength{\abovedisplayskip}{3pt}
\setlength{\belowdisplayskip}{3pt} 
\setlength{\abovedisplayshortskip}{3pt} 
\setlength{\belowdisplayshortskip}{3pt}

\maketitle

\begin{abstract}
  As a mathematical model of high-speed flow and shock wave propagation in a
  complex multimaterial setting, Lagrangian hydrodynamics is characterized by
  moving meshes, advection-dominated solutions, and moving shock fronts with
  sharp gradients. These challenges hinder the existing projection-based model
  reduction schemes from being practical. We develop several variations of
  projection-based reduced order model techniques for Lagrangian hydrodynamics
  by introducing three different reduced bases for position, velocity, and
  energy fields. A time-windowing approach is also developed to address the
  challenge imposed by the advection-dominated solutions.  Lagrangian
  hydrodynamics is formulated as a nonlinear problem, which requires a proper
  hyper-reduction technique.  Therefore, we apply the over-sampling DEIM and SNS
  approaches to reduce the complexity due to the nonlinear terms. Finally, we
  also present both a posteriori and a priori error bounds associated with our
  reduced order model.  We compare the performance of the spatial and
  time-windowing reduced order modeling approaches in terms of accuracy and
  speed-up with respect to the corresponding full order model for several
  numerical examples, namely Sedov blast, Gresho vortices, Taylor-Green
  vortices, and triple-point problems. 
\end{abstract}

\begin{keywords} 
reduced order model, hyper-reduction, 
hydrodynamics, compressible flow, 
Lagrangian methods, advection-dominated problems
\end{keywords}


\section{Introduction}\label{sec:intro}

Physical simulations are key to developments in science, engineering, and
technology.  Many physical processes are mathematically modeled by
time-dependent nonlinear partial differential equations.  In many applications,
the analytical solution of such problems can not be obtained, and numerical
methods are developed to approximate the solutions efficiently.  However,
subject to the complexity of the model problem and the size of the problem
domain, the computational cost can be prohibitively high.  It may take a long
time to run one forward simulation even with high performance computing.  In
decision-making applications where multiple forward simulations are needed, such
as parameter study, design optimization \cite{wang2007large, de2020three,
de2018adaptive, white2020dual}, optimal control \cite{choi2015practical,
choi2012simultaneous}, uncertainty quantification \cite{smith2013uncertainty,
biegler2011large}, and inverse problems \cite{galbally2010non,
biegler2011large}, the computationally expensive simulations are not desirable.
To this end, a reduced order model (ROM) can be useful in this context to obtain
sufficiently accurate approximate solutions with considerable speed-up
compared to a corresponding full order model (FOM). 

Many model reduction schemes have been developed to reduce the computational
cost of simulations while minimizing the error introduced in the reduction
process. Most of these approaches seek to extract an intrinsic solution subspace
for condensed solution representation by a linear combination of reduced basis
vectors. The reduced basis vectors are extracted from performing proper
orthogonal decomposition (POD) on the snapshot data of the FOM simulations.  The
number of degrees of freedom is then reduced by substituting the ROM solution
representation into the (semi-)discretized governing equation.  These approaches
take advantage of both the known governing equation and the solution data
generated from the corresponding FOM simulations to form linear subspace reduced
order models (LS-ROM).  Example applications include, but are not limited to, the nonlinear
diffusion equations \cite{hoang2020domain, fritzen2018algorithmic}, the Burgers
equation and the Euler equations in small-scale \cite{choi2019space,
choi2020sns, carlberg2018conservative}, the convection--diffusion equations
\cite{mojgani2017lagrangian, kim2020efficientII}, the Navier--Stokes equations
\cite{xiao2014non, burkardt2006pod}, rocket nozzle shape design
\cite{amsallem2015design}, flutter avoidance wing shape optimization
\cite{choi2020gradient}, topology optimization of wind turbine blades
\cite{choi2019accelerating}, lattice structure design
\cite{mcbane2020component}, porous media flow/reservoir simulations
\cite{ghasemi2015localized, jiang2019implementation, yang2016fast,
wang2020generalized}, computational electro-cardiology \cite{yang2017efficient},
inverse problems \cite{fu2018pod}, shallow water equations \cite{zhao2014pod,
cstefuanescu2013pod}, Boltzmann transport problems \cite{choi2021space},
computing electromyography \cite{mordhorst2017pod}, spatio-temporal dynamics of
a predator--prey system \cite{dimitriu2013application}, acoustic wave-driven
microfluidic biochips \cite{antil2012reduced}, and Schr{\"o}dinger equation
\cite{cheng2016reduced}.  Survey papers for the projection-based LS-ROM
techniques can be found in \cite{gugercin2004survey, benner2015survey}. 

In spite of successes of the classical LS-ROM in many applications, these
approaches are limited to the assumption that the intrinsic solution space falls
into a subspace with a small dimension, i.e., the solution space has a small
Kolmogorov $n$-width.  This assumption is violated in advection-dominated
problems,  such as sharp gradients, moving shock fronts, and turbulence, which
hinders these model reduction schemes from being practical.  Our goal in this
paper is to develop an efficient reduced order model for hydrodynamics
simulation.  Some reduced order model techniques for hydrodynamics or turbulence
models in the literature include \cite{mou2020data, parish2017non,
gadalla2020comparison, bergmann2009enablers, osth2014need, baiges2015reduced,
san2018extreme}, which are mostly built on the Eulerian formulation, i.e., the
computational mesh is stationary with respect to the fluid motion.  In contrast,
numerical methods in the Lagrangian formulation, which are characterized by a
computational mesh that moves along with the fluid velocity, are developed for
better capturing the shocks and preserving the conserved quantities in
advection-dominated problems.  It therefore becomes natural to develop
Lagrangian-based reduced order models to overcome the challenges posed by
advection-dominated problems.  Some existing work in this research direction
include \cite{mojgani2017lagrangian,lu2020lagrangian}, where a Lagrangian POD
and dynamic mode decomposition (DMD) reduced order model are introduced
respectively for the one-dimensional nonlinear advection-diffusion equation.  We
remark that there are some similarities and differences between our work and
\cite{mojgani2017lagrangian} in using POD for developing Lagrangian-based
reduced order models.  It is important to note that our work is based on the
more complicated and challenging two-dimensional or three-dimensional
compressible Euler equations.  Additionally, the reduced bases are built
independently for each state variable in our work, while a single basis is built
for the whole state in \cite{mojgani2017lagrangian}.  Furthermore, we introduce
the time-windowing concept so as to ensure adequate ROM speed-up by decomposing
the time frame into small time windows and building a temporally-local ROM for
each window. 

Recently, there have been many attempts to develop efficient ROMs for the
advection-dominated or sharp gradient problems.  The attempts can be divided
mainly into two categories.  The first category enhances the solution
representability of the linear subspace by introducing some special treatments
and adaptive schemes.  A dictionary-based model reduction method for the
approximation of nonlinear hyperbolic equations is developed in
\cite{abgrall2016robust}, where the reduced approximation is obtained from the
minimization of the residual in the $L_1$ norm for the reduced linear subspace.
In \cite{carlberg2015adaptive}, a fail-safe $h$-adaptive algorithm is developed.
The algorithm enables ROMs to be incrementally refined to capture the shock
phenomena which are unobserved in the original reduced basis through
a-posteriori online enrichment of the reduced-basis space by decomposing a given
basis vector into several vectors with disjoint support.  The windowed
least-squares Petrov--Galerkin model reduction for dynamical systems with
implicit time integrators is introduced in \cite{parish2019windowed,
shimizu2020windowed}, which can overcome the challenges arising from the
advection-dominated problems by representing only a small time window with a
local ROM.  Another active research direction is to exploit the sharp gradients
and represent spatially local features in ROM, such as the online adaptivity
bases and adaptive sampling approach \cite{peherstorfer2018model} and the shock
reconstruction surrogate approach \cite{constantine2012reduced}.  In
\cite{taddei2020space}, an adaptive space-time registration-based model
reduction is used to align local features of parameterized hyperbolic PDEs in a
fixed one-dimensional reference domain.  Some new approaches have been developed
for aligning the sharp gradients by using a superposition of snapshots with
shifts or transforms.  In \cite{reiss2018shifted}, the shifted proper orthogonal
decomposition (sPOD) introduces time-dependent shifts of the snapshot matrix in
POD in an attempt to separate different transport velocities in advection-dominated
problems.  The practicality of this approach relies heavily on accurate
determination of shifted velocities.  In \cite{rim2018transport}, an iterative
transport reversal algorithm is proposed to decompose the snapshot matrix into
multiple shifting profiles.  In \cite{welper2020transformed},  inspired by the
template fitting \cite{kirby1992reconstructing}, a high resolution transformed
snapshot interpolation with an appropriate behavior near singularities is
considered. 

The second category replaces the linear subspace solution representation with
the nonlinear manifold, which is a very active research direction.  Recently, a
neural network-based reduced order model is developed in \cite{lee2020model} and
extended to preserve the conserved quantities in the physical conservation laws
\cite{lee2019deep}.  In these approaches, the weights and biases in the neural
network are determined in the training phase, and existing numerical methods,
such as finite difference and finite element methods, are utilized.  However,
since the nonlinear terms need to be updated every time step or Newton step, and
the computation of the nonlinear terms still scale with the FOM size, these
approaches do not achieve any speed-up with respect to the corresponding FOM.
Recently, Kim, et al., have achieved a considerable speed-up with the nonlinear
manifold reduced order model \cite{kim2020fast, kim2020efficient}, but it was
only applied to small problems.  Manifold approximations via transported
subspaces in \cite{rim2019manifold} introduced a low-rank approximation to the
transport dynamics by approximating the solution manifold with a transported
subspace generated by low-rank transport modes.  However, their approach is
limited to one-dimensional problem setting.  In \cite{rim2020depth}, a depth
separation approach for reduced deep networks in nonlinear model reduction is
presented, in which the reduced order model is composed with hidden layers with
low-rank representation.

In this paper, we present an alternative reduced order model technique for
advection-dominated problems.  We consider the advection-dominated problems
arising in compressible gas dynamics. The Euler equation is used to model the
high-speed flow and shock wave propagation in a complex multimaterial setting,
and numerically solved in a moving Lagrangian frame, where the computational
mesh is moved along with the fluid velocity. In computational fluid dynamics,
the Lagrangian method is widely used to model different fluid phenomena, for
instance, the immersed boundary method
\cite{peskin2002immersed,boffi2016discrete,cheung2018mass} for fluid-structure
interaction, the Lagrangian particle tracking method
\cite{yeung2002Lagrangian,luthi2005Lagrangian} for turbulence, and the Arbitrary
Lagrangian-Eulerian method \cite{hirt1974arbitrary,benson1989efficient,
guermond2017invariant} for general deformation problems.  In
\cite{dobrev2012high}, a general framework of high-order curvilinear finite
elements and adaptive time stepping of explicit time integrators is proposed for
numerical discretization of the Lagrangian hydrodynamics problem over general
unstructured two-dimensional and three-dimensional computational domains.  Using
general high-order polynomial basis functions for approximating the state
variables, and curvilinear meshes for capturing the geometry of the flow and
maintaining robustness with respect to mesh motion, the method achieves
high-order accuracy.  On the other hand, a modification is made to the
second-order Runge-Kutta method to compensate for the lack of total energy
conservation in standard high-order time integration techniques.  The
introduction of an artificial viscosity tensor further generates the appropriate
entropy and ensures the Rankine-Hugoniot jump conditions at a shock boundary.
Although the method shows great capability and lots of advantages, forward
simulations can be computationally very expensive, especially in
three-dimensional applications with high-order finite elements over fine meshes.
It is therefore desirable to develop efficient ROM techniques for the Lagrangian
hydrodynamics simulation.

We adopt a time-windowing approach for reduced order
modeling of Lagrangian hydrodynamics. The time-windowing approach is introduced to
handle the difficulties arising from advection-dominated problems.  Two
different time window division mechanisms will be considered, namely by the physical time
or the number of snapshots. Several techniques of
construction of offset variables which serve as reference points in the time
windows will be introduced.  For a given time window, proper orthogonal
decomposition is used to extract the dominant modes in solution
representability, and an oversampling hyper-reduction technique is employed to
reduce the complexity due to the nonlinear terms in the governing equations.  We
will present error estimates to theoretically justify our method and numerical
examples to exhibit the capabilities of our method.  For the purpose of
reproducible research, open source codes are available as a branch of the Laghos
GitHub repository\footnote{GitHub page, {\it
https://github.com/CEED/Laghos/tree/rom}.} and libROM GitHub
repository\footnote{GitHub page, {\it https://github.com/LLNL/libROM}.}.

The main contributions of this paper are summarized as follows:
\begin{itemize}
  \item We present a parametric time-windowing reduced order model for
    Lagrangian hydrodynamics.
  \item We follow the SNS procedure in \cite{choi2020sns} to derive an efficient 
  hyper-reduction technique for the compressible Euler equations.
  \item We introduce a new mechanism of decomposing the time domain in adaptive
    time stepping schemes to ensure uniform reduced order model size.
  \item We introduce different procedures of offset variables 
    which serve as reference points in the time windows. 
  \item We derive several error bounds for the reduced order model of
    Lagrangian hydrodynamics.
  \item We present numerical results of the reduced order model on 
  two-dimensional or three-dimensional compressible Euler equations.
\end{itemize}

\subsection{Organization of the paper}\label{sec:organization} 
In Section~\ref{sec:FOM}, we introduce the semidiscrete Lagrangian conservation
laws. A projection-based ROM is described in Section~\ref{sec:ROM}, and the
time-windowing approach is introduced in Section~\ref{sec:timewindowing}. A
posteriori and a priori error bounds are derived for our Lagrangian
hydrodynamics ROM in Section~\ref{sec:errorbound}.  Numerical results are
presented in Section~\ref{sec:numericalresults}, and the conclusion is summarized
in Section~\ref{sec:conclusion}. In Appendix~\ref{sec:appendix}, Laghos command
line options are provided for each of the numerical experiments presented.

\section{Lagrangian hydrodynamics}\label{sec:FOM}
We consider the system of Euler equations of gas dynamics in a Lagrangian
reference frame \cite{harlow1971fluid}, assuming no external body force is
exerted:
\begin{equation}\label{eq:euler}
  \begin{aligned}
    \text{momentum conservation}:& & \densitySymbol\frac{d\velocitySymbol}{d\timeSymbol} &=
    \gradientSymbol \cdot \stressSymbol \\
    \text{mass conservation}:& & \dfrac{1}{\densitySymbol}\frac{d\densitySymbol}{d\timeSymbol} &=
    -\gradientSymbol \cdot \velocitySymbol \\
    \text{energy conservation}:& & \densitySymbol\frac{d\energySymbol}{d\timeSymbol} &=
    \stressSymbol :  \gradientSymbol \velocitySymbol \\
    \text{equation of motion}:& & \frac{d\positionSymbol}{d\timeSymbol} &= \velocitySymbol.
  \end{aligned}
\end{equation}
Here, $\frac{d}{d\timeSymbol} = \frac{\partial}{\partial \timeSymbol} + 
\velocitySymbol \cdot \gradientSymbol$ is the material derivative, 
$\densitySymbol$ denotes the density of the fluid, $\positionSymbol$ and
$\velocitySymbol$ denote the position and the velocity of the particles in a
deformable medium $\solDomainSymbol(\timeSymbol)$ in the Eulerian coordinates,
$\stressSymbol$ denotes the deformation stress tensor, and $\energySymbol$
denotes the internal energy per unit mass.  These physical quantities can be
treated as functions of the time $\timeSymbol$ and the particle
$\initialPosition \in \initialDomain = \solDomainSymbol(0)$.  In gas dynamics,
the stress tensor is isotropic, and we write $\stressSymbol = -\pressureSymbol
\identitySymbol + \artificialStressSymbol$, where $\pressureSymbol$ denotes the
thermodynamic pressure, and $\artificialStressSymbol$ denotes the artificial
viscosity stress.  The thermodynamic pressure is described by the equation of
state, and can be expressed as a function of the density and the internal
energy.  In our work, we focus on the case of polytropic ideal gas with an
adiabatic index $\adiabaticIndexSymbol > 1$, which yields the equation of state 
\begin{equation}\label{eq:EOS}
  \pressureSymbol = (\adiabaticIndexSymbol - 1) \densitySymbol \energySymbol.
\end{equation}
The system is prescribed with an initial condition and a boundary condition
$\velocitySymbol \cdot \normalSymbol = \neumannSymbol$, where $\normalSymbol$ is
the outward normal unit vector on the domain boundary.  Moreover, a set of
problem parameters $\param \in \paramDomain$ determines certain physical data in
the system of Euler equations \eqref{eq:euler}, and therefore the physical
quantities are parametrized with the data $\param$.  

\subsection{Spatial discretization}
Following \cite{dobrev2012high}, we adopt a spatial discretization for
\eqref{eq:euler} using a kinematic space $\kinematicFE \subset
[H^1(\initialDomain)]^\dimensionSymbol$ for approximating the position and the
velocity, and a thermodynamic space $\thermodynamicFE \subset
L_2(\initialDomain)$ for approximating the energy. The density can be
eliminated, and the equation of mass conservation can be decoupled from
\eqref{eq:euler}. We assume high-order finite element (FEM) discretization in
space, so that the finite dimensions $\sizeKinematicFE$ and
$\sizeThermodynamicFE$ are the global numbers of FEM degrees of freedom in the
corresponding discete FEM spaces.  For more details, see \cite{dobrev2012high}.
The FEM coefficient vector functions for velocity and position are denoted as
$\velocity, \position:[0,\finalTime]\times \paramDomain \mapto
\RR{\sizeKinematicFE}$,
and the coefficient vector function for energy is denoted as
$\energy:[0,\finalTime] \times \paramDomain \mapto \RR{\sizeThermodynamicFE}$.
The semidiscrete Lagrangian conservation laws can be expressed as a nonlinear
system of differential equations in the coefficients with respect to the bases
for the kinematic and thermodynamic spaces:
\begin{equation}\label{eq:fom}
  \begin{aligned}
    \text{momentum conservation}:& & \kinematicMassMat\frac{d\velocity}{d\timeSymbol} &=
    -\forceMat(\velocity, \energy, \position; \param) \cdot\oneVec \\
    \text{energy conservation}:& & \thermodynamicMassMat\frac{d\energy}{d\timeSymbol} &=
    \forceMat(\velocity, \energy, \position; \param)^T\cdot\velocity \\
    \text{equation of motion}:& & \frac{d\position}{d\timeSymbol} &= \velocity.
  \end{aligned}
\end{equation}

Let $\state \equiv (\velocity; \energy; \position)^T\in\RR{\sizeWholeFE}$,
$\sizeWholeFE = 2\sizeKinematicFE + \sizeThermodynamicFE$, be the hydrodynamic
state vector. Then the semidiscrete conservation equation of \eqref{eq:fom} can
be written in a compact form as
\begin{equation}\label{eq:combinedFOM}
  \frac{d\state}{d\timeSymbol} = \forceSys(\state,\timeSymbol;\param),  
\end{equation}
where the nonlinear force term, $\forceSys:\RR{\sizeWholeFE} \times
\paramDomain \mapto \RR{\sizeWholeFE} $, is defined as 
\begin{align}\label{eq:combinedForce}
  \forceSys(\state;\param) &\equiv
  \pmat{\forceSysVelocity(\velocity, \energy, \position) \\
        \forceSysEnergy(\velocity, \energy, \position) \\
        \forceSysPosition(\velocity, \energy, \position)} 
        \equiv
  \pmat{-\kinematicMassMat^{-1}\forceOne \\
        \thermodynamicMassMat^{-1}\forceTv \\
        \velocity},
\end{align}
where $\forceOne: \RR{\sizeWholeFE} \times 
\paramDomain \mapto \RR{\sizeKinematicFE}$ and $\forceTv: \RR{\sizeWholeFE}
\times \paramDomain \mapto \RR{\sizeThermodynamicFE}$ are
nonlinear vector functions that are defined respectively as 
\begin{align}\label{eq:forceOneforceTv}
  \forceOne &\equiv \forceMat\cdot\oneVec, & 
  \forceTv &\equiv \forceMat^T\cdot\velocity.
\end{align}

\subsection{Time integrators}\label{sec:time_integrator}
In order to obtain a fully discretized system of equations, one needs to apply a
time integrator. We consider two different explicit Runge-Kutta schemes: the
RK2-average and RK4 schemes. The temporal domain is discretized as $\{
  \timek{\timeIndex} \}_{\timeIndex=0}^{\ntimestep}$, where $\timek{n}$ denotes
a discrete moment in time with $\timek{0} = 0$, $\timek{\ntimestep} =
\finalTime$, and $\timek{n-1} < \timek{n}$ for
$\timeIndex\in\nat{\ntimestep}$, where $\nat{N}\equiv \{1,\ldots,N\}$.  The
computational domain at time $\timek{\timeIndex}$ is denoted as
$\Omega^{\timeIndex} \equiv \Omega(\timek{\timeIndex})$. We denote the
quantities of interest defined on $\Omega^{\timeIndex}$ with a subscript
$\timeIndex$.  

\subsubsection{The RK2-average scheme}\label{sec:RK2-avg}
The midpoint Runge--Kutta second-order scheme is written as
\begin{align}\label{eq:midpointRK2}
  \statet{\timeIndex+\frac{1}{2}} &= \statet{\timeIndex} +
    \frac{\timestepk{\timeIndex}}{2}\forceSys(\statet{\timeIndex}), &
  \statet{\timeIndex+1} &= \statet{\timeIndex} +
    \timestepk{\timeIndex}\forceSys(\statet{\timeIndex+\frac{1}{2}}),
\end{align}
where $\timestepk{\timeIndex} \equiv t_{\timeIndex+1} - t_{\timeIndex}$.
In practice, the midpoint RK2 scheme can be unstable even for simple test
problems. Therefore, the following RK2-average scheme is used:
\begin{align}\label{eq:RK2-avg}
  \velocityt{\timeIndex+\frac{1}{2}} &= \velocityt{\timeIndex} - (\timestepk{\timeIndex}/2)
    \kinematicMassMat^{-1} \forceOnek{\timeIndex}, &
    \velocityt{\timeIndex+1} &= \velocityt{\timeIndex} - \timestepk{\timeIndex}
    \kinematicMassMat^{-1} \forceOnek{\timeIndex+\frac{1}{2}}, \\
  \energyt{\timeIndex+\frac{1}{2}} &= \energyt{\timeIndex} + (\timestepk{\timeIndex}/2)
    \thermodynamicMassMat^{-1} \forceTvk{\timeIndex}, &
    \energyt{\timeIndex+1} &= \energyt{\timeIndex} + \timestepk{\timeIndex}
    \thermodynamicMassMat^{-1} \avgforceTvk{\timeIndex+\frac{1}{2}}, \\
  \positiont{\timeIndex+\frac{1}{2}} &= \positiont{\timeIndex} + (\timestepk{\timeIndex}/2)
    \velocityt{\timeIndex+\frac{1}{2}}, & \positiont{\timeIndex+1} &=
    \positiont{\timeIndex} + \timestepk{\timeIndex} \avgvelocityt{\timeIndex+\frac{1}{2}},
\end{align}
where the state $\statet{\timeIndex} 
= (\velocityt{\timeIndex}; \energyt{\timeIndex}; 
\positiont{\timeIndex})^T \in\RR{\sizeWholeFE}$ 
is used to compute the updates 
\begin{align}
\forceOnek{\timeIndex} & = \left(\forceMat (\statet{\timeIndex}) \right ) \cdot\oneVec, &
\forceTvk{\timeIndex} & = \left(\forceMat (\statet{\timeIndex}) \right )^T \cdot 
\velocityt{\timeIndex+\frac{1}{2}},
\end{align}
in the first stage. Similarly, $\statet{\timeIndex+\frac{1}{2}} = 
(\velocityt{\timeIndex+\frac{1}{2}}; \energyt{\timeIndex+\frac{1}{2}}; 
\positiont{\timeIndex+\frac{1}{2}})^T \in\RR{\sizeWholeFE}$ 
is used to compute the updates 
\begin{align}
\forceOnek{\timeIndex+\frac{1}{2}} 
& = \left(\forceMat (\statet{\timeIndex+\frac{1}{2}}) \right ) \cdot \oneVec, &
\avgforceTvk{\timeIndex+\frac{1}{2}} & = \left
(\forceMat (\statet{\timeIndex + \frac{1}{2}}) \right )^T \cdot
\avgvelocityt{\timeIndex+\frac{1}{2}},
\end{align}
with $\avgvelocityt{\timeIndex+\frac{1}{2}} = (\velocityt{\timeIndex} +
\velocityt{\timeIndex+1})/2$ in the second stage.  Note that the RK2-average
scheme is different from the midpoint RK2 scheme in the updates for energy and
position.  The RK2-average scheme uses the midpoint velocity
$\velocityt{\timeIndex+\frac{1}{2}}$ and the average velocity
$\avgvelocityt{\timeIndex+\frac{1}{2}}$ to update energy and position in the
first stage and the second stage respectively, while the midpoint RK2 uses the
initial velocity $\velocityt{\timeIndex}$ and the midpoint velocity
$\velocityt{\timeIndex+\frac{1}{2}}$ respectively. The RK2-average scheme is
proved to conserve the discrete total energy exactly (see Proposition 7.1 of
\cite{dobrev2012high}). 

\subsubsection{The RK4 scheme}\label{sec:RK4}
The traditional RK4 scheme can be applied to the Lagrangian hydrodynamics in
\eqref{eq:combinedFOM} as
\begin{align}\label{eq:RK4}
  \statet{\timeIndex+1} = \statet{\timeIndex} + \frac{1}{6} \timestepk{\timeIndex}
  (\rkcoeffI + 2\rkcoeffII + 2\rkcoeffIII + \rkcoeffIV),
\end{align}
where 
\begin{align}\label{eq:RK4coeff}
  \rkcoeffI   &= \forceSys(\statet{\timeIndex}), &
  \rkcoeffII  &=
    \forceSys(\statet{\timeIndex}+\frac{\timestepk{\timeIndex}}{2}\rkcoeffI),&
  \rkcoeffIII &= \forceSys(\statet{\timeIndex}+\frac{\timestepk{\timeIndex}}{2}\rkcoeffII), &
  \rkcoeffIV  &= \forceSys(\statet{\timeIndex}+\timestepk{\timeIndex}\rkcoeffIII).
\end{align}

\subsubsection{Adaptive time stepping}\label{sec:adaptiveDelta_t}
Since explicit Runge-Kutta methods are used, we need to control the time step
size in order to maintain the stability of the fully discrete schemes.  We
follow the automatic time step control algorithm described in Section 7.3 of
\cite{dobrev2012high}, which we briefly describe here.  At the time step
$\timeIndex$, the algorithm starts with a time step estimate
$\timestepestimatek{\timeIndex}$ defined as 
\begin{equation}\label{eq:adaptive-timestep}
\timestepestimatek{\timeIndex} = \min_{\positionSymbol} 
\CFLconst \left(\dfrac{c_s(\positionSymbol)}{h_{\text{min}}(\positionSymbol)} 
+ \CFLconst_\mu \dfrac{\mu_s(\positionSymbol)}{\densitySymbol(\positionSymbol) 
h^2_{\text{min}}(\positionSymbol)}\right)^{-1}, 
\end{equation}
where the minimum is taken over all quadrature points used in the evaluation of
the local force matrices and over all Runge-Kutta stages. Here, $\CFLconst$ and
$\CFLconst_\mu$ are certain Courant--Friedrichs--Lewy (CFL) constants.  The
default values we use are $\CFLconst = 0.5$ and $\CFLconst_\mu = 2.5$.
Moreover, $c_s$ is the speed of sound, $\mu_s$ is the viscosity coefficient, and
$h_{\text{min}}$ is the minimal singular value of the zone Jacobian.  These
quantities are used for the evaluation of the artificial stress tensor for
modeling shock wave propagation. The discussion of artificial viscosity is beyond
the scope of this paper, and the reader is referred to Section 6 of
\cite{dobrev2012high} for more details.  With the estimate
$\timestepestimatek{\timeIndex}$, we use the following algorithm to control the
time-step:
\begin{steps}
  \item Given a time-step $\timestep$ and state $\statet{\timeIndex}$, evaluate
    the state $\statet{\timeIndex+1}$ and the corresponding time step estimate
    $\timestepestimatek{\timeIndex}$.
  \item If $\timestep \geq \timestepestimatek{\timeIndex}$, set $\timestep \leftarrow
    \beta_1 \timestep$ and go to Step 1.
  \item If $\timestep \leq \gamma \timestepestimatek{\timeIndex}$, set
    $\timestep \leftarrow \beta_2 \timestep$.
  \item Set $\timestepk{\timeIndex} = \timestep$, 
  $\timeIndex \leftarrow \timeIndex + 1$ and continue with the next time step.
\end{steps}
Here, $\beta_1$, $\beta_2$, and $\gamma$ denote given constants. The default
values we use are $\beta_1 = 0.85$, $\beta_2 = 1.02$, and $\gamma = 0.8$.  We
remark that, due to the automatic time step control, the temporal discretization
is not determined a-priori and depends on various model inputs, such as the
underlying state space, the CFL constant and the time integrator. 

\section{Reduced order model}\label{sec:ROM}
In this section, we present the details of the projection-based reduced oder
model for the semi-discrete Lagrangian conservation laws \eqref{eq:fom}. We
start with the usage of the reduced order model, which we refer to as the online phase.
We sequentially discuss the solution representation by a subspace in
Section~\ref{sec:solrepresentation}, Galerkin projection in Section~\ref{sec:Galerkin} and
hyper-reduction in Section~\ref{sec:hyperreduction}, resulting in a continuous-in-time reduced order
model. Using a time integrator in Section~\ref{sec:rom_time_integrator}, 
we obtain a fully discrete reduced order model. 
Then we move on to discuss the construction of the reduced order model, 
which is precomputed once in the offline phase. 
We sequentially present the details of construction 
of the solution subspace by proper orthogonal decomposition on snapshot matrices in 
Section~\ref{sec:POD}, a technique of construction of the nonlinear term bases by 
a conforming subspace relation in Section~\ref{sec:SNS}, and the construction of the
sampling indices by the discrete empirical interpolation method (DEIM) in Section~\ref{sec:DEIM}.

\subsection{Solution representation}\label{sec:solrepresentation}
We restrict our solution space to a subspace spanned by a reduced basis for each
field. That is, the subspace for velocity, energy, and position fields are
defined as
\begin{align}\label{eq:subspces}
  \velocitySubspace &\equiv
  \Span{\velocityBasisVeck{\basisIndex}}_{\basisIndex=1}^{\sizeROMvelocity} \subseteq
  \RR{\sizeKinematicFE},&
  \energySubspace &\equiv
  \Span{\energyBasisVeck{\basisIndex}}_{\basisIndex=1}^{\sizeROMenergy} \subseteq
  \RR{\sizeThermodynamicFE},&
  \positionSubspace &\equiv
  \Span{\positionBasisVeck{\basisIndex}}_{\basisIndex=1}^{\sizeROMposition} \subseteq
  \RR{\sizeKinematicFE}, 
\end{align}
with $\dim(\velocitySubspace)=\sizeROMvelocity\ll\sizeKinematicFE$,
$\dim(\energySubspace)=\sizeROMenergy\ll\sizeThermodynamicFE$,  and
$\dim(\positionSubspace)=\sizeROMposition\ll\sizeKinematicFE$. Using these
subspaces, each discrete field is approximated in trial subspaces, $\velocity
\approx \velocityApprox \in \velocityOS + \velocitySubspace$, $\energy \approx
\energyApprox \in \energyOS + \energySubspace$, and $\position \approx
\positionApprox \in \positionOS + \positionSubspace$, or equivalently
\begin{align}\label{eq:solrepresentation}
  \velocityApprox(\timeSymbol; \param) &= \velocityOS(\param) + \velocityBasis\ROMvelocity(\timeSymbol; \param), \\
  \energyApprox(\timeSymbol; \param) &= \energyOS(\param) + \energyBasis\ROMenergy(\timeSymbol; \param), \\
  \positionApprox(\timeSymbol; \param) &= \positionOS(\param) + \positionBasis\ROMposition(\timeSymbol; \param), 
\end{align}
where 
$\velocityOS(\param)\in\RR{\sizeKinematicFE}$,
$\energyOS(\param)\in\RR{\sizeThermodynamicFE}$, and
$\positionOS(\param)\in\RR{\sizeKinematicFE}$ denote the prescribed offset
vectors for velocity, energy, and position fields respectively; the orthonormal basis
matrices
$\velocityBasis\in\RR{\sizeKinematicFE\times\sizeROMvelocity}$,
$\energyBasis\in\RR{\sizeThermodynamicFE\times\sizeROMenergy}$, and
$\positionBasis\in\RR{\sizeKinematicFE\times\sizeROMposition}$
are defined as
\begin{align}\label{eq:basisMats}
  \velocityBasis &\equiv \bmat{\velocityBasisVeck{1} & \cdots & \velocityBasisVeck{\sizeROMvelocity}},&
  \energyBasis &\equiv \bmat{\energyBasisVeck{1} & \cdots & \energyBasisVeck{\sizeROMenergy}},&
  \positionBasis &\equiv \bmat{\positionBasisVeck{1} & \cdots & \positionBasisVeck{\sizeROMposition}};
\end{align}
and $\ROMvelocity:[0,\finalTime]\times \paramDomain \mapto
\RR{\sizeROMvelocity}$, $\ROMenergy:[0,\finalTime]\times \paramDomain \mapto
\RR{\sizeROMenergy}$, and $\ROMposition:[0,\finalTime]\times \paramDomain \mapto
\RR{\sizeROMposition}$ denote the time-dependent generalized coordinates for
velocity, energy, and position fields, respectively.  One natural choice of the
offset vectors is to use the initial values, i.e.  $\velocityOS(\param) = \velocity(0;
\param)$, $\energyOS(\param) = \energy(0; \param)$, and $\positionOS(\param) =
\position(0; \param)$.  Replacing $\velocity$, $\energy$, and $\position$ with
$\velocityApprox$, $\energyApprox$, and $\positionApprox$ in
Eqs.~\eqref{eq:fom}, the semi-discretized Lagrangian hydrodynamics system
becomes the following over-determined semi-discrete system:
\begin{equation}\label{eq:rom-sol}
  \begin{aligned}
    \kinematicMassMat \velocityBasis \frac{d\ROMvelocity}{d\timeSymbol} &=
    -\forceOne(\velocityOS+\velocityBasis\ROMvelocity,
    \energyOS+\energyBasis\ROMenergy, \positionOS+\positionBasis\ROMposition,
    \timeSymbol; \param)\\
    \thermodynamicMassMat \energyBasis \frac{d\ROMenergy}{d\timeSymbol} &=
    \forceTv(\velocityOS+\velocityBasis\ROMvelocity, \energyOS+\energyBasis\ROMenergy, \positionOS+\positionBasis\ROMposition, \timeSymbol; \param) \\
    \positionBasis \frac{d\ROMposition}{d\timeSymbol} &= \velocityOS + \velocityBasis \ROMvelocity.
  \end{aligned}
\end{equation}
The system of ordinary differential equations is closed by defining the initial
condition at $\timeSymbol = 0$.  With the initial values as the offset vectors,
using the solution representation \eqref{eq:solrepresentation}, one can derive
the initial values of the ROM coefficient vectors given by  zero vectors in the
corresponding ROM spaces, i.e.  $\ROMvelocity(0; \param) =
\zeroVec{\sizeROMvelocity}$, $\ROMenergy(0; \param) = \zeroVec{\sizeROMenergy}$,
and $\ROMposition(0; \param) = \zeroVec{\sizeROMposition}$, where $\zeroVec{m}
\in \RR{m}$ is the zero vector in $\RR{m}$. 

\subsection{Galerkin projection}\label{sec:Galerkin}
Note that Eqs.~\eqref{eq:rom-sol} has more equations than unknowns (i.e., an
over-determined system).  It is likely that there is no solution satisfying
Eq.~\eqref{eq:rom-sol}. The system must be closed in order to obtain a solution.
One can invert the mass matrices in the momentum conservation and energy
conservation of Eqs.~\eqref{eq:rom-sol} and apply Galerkin projection to obtain
the following reduced system of semi-discrete Lagrangian hydrodynamics
equations:
\begin{equation}\label{eq:rom-denseRHS}
  \begin{aligned}
    \frac{d\ROMvelocity}{d\timeSymbol} &= -\velocityBasis^T
      \kinematicMassMat^{-1}\forceOne(\velocityOS+\velocityBasis\ROMvelocity,
      \energyOS+\energyBasis\ROMenergy, \positionOS+\positionBasis\ROMposition,
      \timeSymbol; \param) \\ 
    \frac{d\ROMenergy}{d\timeSymbol} &=
      \energyBasis^T \thermodynamicMassMat^{-1}
      \forceTv(\velocityOS+\velocityBasis\ROMvelocity,
      \energyOS+\energyBasis\ROMenergy, \positionOS+\positionBasis\ROMposition,
      \timeSymbol; \param) \\
    \frac{d\ROMposition}{d\timeSymbol} &= \positionBasis^T \velocityOS +
      \positionBasis^T \velocityBasis \ROMvelocity.
  \end{aligned}
\end{equation}
Here, we use the assumption of basis orthonormality.
There is one major issue with this Galerkin formulation. The nonlinear term
involves the inverse of the mass matrices, which are dense matrices. This makes
each component of the nonlinear term vector dependent on potentially every state
variable component. Even though a gappy POD is used with a small number of sample points,
the evaluation of each point requires possibly all the state variables,
resulting in inefficient computational complexity. Therefore, we choose an
alternative Galerkin projection, in which we do not invert mass matrices.
Instead, we form reduced mass matrices by applying Galerkin projection directly
to Eqs.~\eqref{eq:rom-sol}, resulting in
\begin{equation}\label{eq:rom}
  \begin{aligned}
   \ROMKinematicMassMat \frac{d\ROMvelocity}{d\timeSymbol} &= -\velocityBasis^T
    \forceOne(\velocityOS+\velocityBasis\ROMvelocity, \energyOS +
    \energyBasis\ROMenergy, \positionOS + \positionBasis\ROMposition,
    \timeSymbol; \param) \\
   \ROMThermodynamicMassMat\frac{d\ROMenergy}{d\timeSymbol} &=
    \energyBasis^T\forceTv(\velocityOS + \velocityBasis\ROMvelocity, \energyOS
    + \energyBasis\ROMenergy, \positionOS +
    \positionBasis\ROMposition, \timeSymbol; \param) \\
    \frac{d\ROMposition}{d\timeSymbol} &= \positionBasis^T \velocityOS +
    \positionBasis^T \velocityBasis \ROMvelocity,
  \end{aligned}
\end{equation}
where the reduced kinematic matrix, $\ROMKinematicMassMat \in
\RR{\sizeROMvelocity \times \sizeROMvelocity}$, and thermodynamic mass matrix,
$\ROMThermodynamicMassMat \in \RR{\sizeROMenergy \times \sizeROMenergy}$, are
defined as
\begin{align}\label{eq:rmass}
  \ROMKinematicMassMat &\equiv \velocityBasis^T \kinematicMassMat \velocityBasis,&
  \ROMThermodynamicMassMat &\equiv \energyBasis^T \thermodynamicMassMat \energyBasis.
\end{align}
Note that the reduced mass matrices can be precomputed once, since the mass
matrices are constant in time (cf. \cite{dobrev2012high}). Furthermore, the
dimensions of the mass matrices are small, so computing their factorizations,
e.g., LU or Cholesky factorizations, can be done efficiently and once for all
time. Alternatively, one can orthogonalize $\velocityBasis$ and $\energyBasis$
with respect to $\kinematicMassMat$ and $\thermodynamicMassMat$, respectively,
i.e., $\ROMKinematicMassMat = \identityMat{\sizeROMvelocity}$ and
$\ROMThermodynamicMassMat = \identityMat{\sizeROMenergy}$, where
$\identityMat{\sizeROMsymbol}\in\RR{\sizeROMsymbol \times \sizeROMsymbol}$ is an
identity matrix. In conclusion, the handling of the reduced mass matrices is
straightforward and efficient.

\subsection{Hyper-reduction}\label{sec:hyperreduction}
The nonlinear matrix function, $\forceMat$, changes every time the state
variables evolve.  Additionally, it needs to be multiplied by the
basis matrices whenever the update in the nonlinear term occurs, which scales
with the full-order model (FOM) size. Therefore, we cannot expect any speed-up
without special treatment of the nonlinear terms.  To overcome this issue,
a hyper-reduction technique needs to be applied, where $\forceOne$ and
$\forceTv$ are approximated as 
 \begin{align}\label{eq:DEIM_approx}
   \forceOne &\approx \forceOneBasis \ROMforceOne,& 
   \forceTv  &\approx \forceTvBasis \ROMforceTv.
 \end{align}
That is, $\forceOne$ and $\forceTv$ are projected onto subspaces
$\forceOneSubspace \equiv
\Span{\forceOneBasisVeck{\basisIndex}}_{\basisIndex=1}^{\sizeROMforceOne}$ and
$\forceTvSubspace \equiv
\Span{\forceTvBasisVeck{\basisIndex}}_{\basisIndex=1}^{\sizeROMforceTv}$, where
$\forceOneBasis \equiv \bmat{\forceOneBasisVeck{1} & \ldots &
\forceOneBasisVeck{\sizeROMforceOne}} \in
\RR{\sizeKinematicFE\times\sizeROMforceOne}$, $\sizeROMforceOne \ll
\sizeKinematicFE$ and $\forceTvBasis \equiv \bmat{\forceTvBasisVeck{1} & \ldots
& \forceTvBasisVeck{\sizeROMforceTv}} \in
\RR{\sizeThermodynamicFE\times\sizeROMforceTv}$, $\sizeROMforceTv \ll
\sizeThermodynamicFE$, denote the nonlinear term basis matrices, and $\ROMforceOne
\in \RR{\sizeROMforceOne}$ and $\ROMforceTv \in \RR{\sizeROMforceTv}$ denote the
generalized coordinates of the nonlinear terms. Following
\cite{chaturantabut2010nonlinear}, the nonlinear term bases $\forceOneBasis$ and
$\forceTvBasis$ can be constructed by applying another POD on the nonlinear term
snapshots obtained from the FOM simulation at every time step.  This implies
that additional snapshot storage and basis construction are required for the
nonlinear terms.  An alternative approach of avoiding the extra cost without
losing the quality of the hyper-reduction is discussed in \cite{choi2020sns},
which will be discussed in more detail in Section~\ref{sec:POD} and Section~\ref{sec:SNS}. 

Now we will show how the generalized coordinates, $\ROMforceOne$, can be
determined by the following interpolation:
 \begin{align}\label{eq:DEIM_interpolation}
   \forceOneSamplingMat^T \forceOne = \forceOneSamplingMat^T
   \forceOneBasis\ROMforceOne,
 \end{align}
where $\forceOneSamplingMat \equiv \bmat{\unitveck{p_1},
\ldots,\unitveck{p_{\sizeROMforceOneSample}}} \in \RR{\sizeKinematicFE \times
\sizeROMforceOneSample}$, $\sizeROMforceOne \leq \sizeROMforceOneSample \ll
\sizeKinematicFE$, is the sampling matrix and $\unitveck{p_i}$ is the $p_i$-th
column vector of the identity matrix
$\identityMat{\sizeKinematicFE}\in\RR{\sizeKinematicFE\times\sizeKinematicFE}$.
Note that Eq.~\eqref{eq:DEIM_interpolation} is an over-determined system.  Thus,
we solve the least-squares problem, i.e.,
 \begin{align}\label{eq:ls-hyper}
   \ROMforceOne &= \argmin_{\dummyVec\in\RR{\sizeROMforceOne}} \|
   \forceOneSamplingMat^T ( \forceOne - \forceOneBasis\dummyVec ) \|_2^2.    
 \end{align}
The solution to the least-squares problem~\eqref{eq:ls-hyper} is 
 \begin{align}\label{eq:sol-ls-hyper}
   \ROMforceOne = (\forceOneSamplingMat^T\forceOneBasis)^{\dagger}\forceOneSamplingMat^T  \forceOne,
 \end{align}
where the Moore--Penrose inverse of a matrix $\dummyMat \in \RR{I \times J}$,
$I\geq J$, with full column rank is defined as $\dummyMat^{\dagger} :=
(\dummyMat^T\dummyMat)^{-1}\dummyMat^T$.  Therefore, Eq.~\eqref{eq:DEIM_approx}
becomes $\forceOne \approx \forceOneObliqProjMat \forceOne$, where 
 \begin{align}\label{eq:DEIM_f}
   \forceOneObliqProjMat = \forceOneBasis
 \left(\forceOneSamplingMat^T\forceOneBasis\right)^{\dagger}\forceOneSamplingMat^T 
 \in \RR{\sizeKinematicFE\times\sizeKinematicFE}
 \end{align} 
is the oblique projection matrix.  Instead of constructing the sampling matrix
$\forceOneSamplingMat$, for efficiency we simply store the sampling indices
$\{p_1,\ldots,p_{\sizeROMforceOneSample}\} \subset \nat{\sizeKinematicFE}$.
More precisely, the reduced matrix $\left(\forceOneSamplingMat^T\forceOneBasis\right)^{\dagger}
  \in \RR{\sizeROMforceOne\times\sizeROMforceOneSample}$
 can be precomputed and stored in the offline phase,
 and is multiplied to the sampled entries
 $\forceOneSamplingMat^T \forceOne \in \RR{\sizeROMforceOneSample}$
 to obtain $\ROMforceOne$ by \eqref{eq:sol-ls-hyper} in the online phase.
The sampling indices are generated by the discrete
empirical interpolation method (DEIM), which is explained in
Section~\ref{sec:DEIM}.  Similarly, for the nonlinear term of the energy
conservation equation, we have
 \begin{align}\label{eq:sol-ls-hyper-2}
   \ROMforceTv = (\forceTvSamplingMat^T\forceTvBasis)^{\dagger}\forceTvSamplingMat^T  \forceTv. 
 \end{align}
We denote the sampling matrix by $\forceTvSamplingMat \in
\RR{\sizeThermodynamicFE \times \sizeROMforceTvSample}$ and the oblique
projection matrix by
\begin{align}\label{eq:rhsproj}
    \forceTvObliqProjMat & = \forceTvBasis \left(\forceTvSamplingMat^T
      \forceTvBasis\right)^{\dagger} \forceTvSamplingMat^T
 \in \RR{\sizeThermodynamicFE\times\sizeThermodynamicFE}.
\end{align}
Then, the hyper-reduced system is obtained by replacing $\forceOne$ and
$\forceTv$ in Eq.~\eqref{eq:rom} with $\forceOneObliqProjMat\forceOne$ and
$\forceTvObliqProjMat\forceTv$, respectively:
\begin{align}\label{eq:rom-hr}
   \ROMKinematicMassMat \frac{d\ROMvelocity}{d\timeSymbol} &=
    -\velocityBasis^T \forceOneObliqProjMat \forceOne(\velocityOS + \velocityBasis\ROMvelocity, \energyOS + \energyBasis\ROMenergy, \positionOS + \positionBasis\ROMposition, \timeSymbol; \param) \\
    \ROMThermodynamicMassMat\frac{d\ROMenergy}{d\timeSymbol} &=
    \energyBasis^T \forceTvObliqProjMat \forceTv(\velocityOS + \velocityBasis\ROMvelocity, \energyOS + \energyBasis\ROMenergy, \positionOS + \positionBasis\ROMposition, \timeSymbol; \param) \\
    \frac{d\ROMposition}{d\timeSymbol} &= \positionBasis^T \velocityOS +
    \positionBasis^T \velocityBasis \ROMvelocity,
\end{align}
Let $\ROMstate \equiv (\ROMvelocity; \ROMenergy;
\ROMposition)^T\in\RR{\sizeWholeROM}$, $\sizeWholeROM = \sizeROMvelocity +
\sizeROMenergy + \sizeROMposition$, be the reduced order hydrodynamic state
vector. Then the semidiscrete hyper-reduced system  \eqref{eq:rom-hr} can
be written in a compact form as
\begin{equation}\label{eq:combinedROM}
  \frac{d\ROMstate}{d\timeSymbol} = \ROMforceSys(\ROMstate,\timeSymbol;\param),  
\end{equation}
where the nonlinear force term, $\forceSys:\RR{\sizeWholeFE} \times
\paramDomain \mapto \RR{\sizeWholeFE} $, is defined as 
\begin{align}\label{eq:combinedROMForce}
  \ROMforceSys(\ROMstate;\param) &\equiv
  \pmat{\ROMforceSysVelocity(\ROMvelocity, \ROMenergy, \ROMposition) \\
        \ROMforceSysEnergy(\ROMvelocity, \ROMenergy, \ROMposition) \\
        \ROMforceSysPosition(\ROMvelocity, \ROMenergy, \ROMposition)} 
        \equiv
  \pmat{-\ROMKinematicMassMat^{-1}\velocityBasis^T \forceOneObliqProjMat \forceOne\\
        \ROMThermodynamicMassMat^{-1}\energyBasis^T \forceTvObliqProjMat \forceTv \\
        \positionBasis^T \velocityApprox}.
\end{align}

\subsection{Time integrators}\label{sec:rom_time_integrator}
Applying the RK2-average scheme in Section~\ref{sec:RK2-avg} to the
hyper-reduced system \eqref{eq:rom-hr}, the RK2-average fully discrete
hyper-reduced system reads:
\begin{align}\label{eq:RK2-avg-rom-hr}
  \ROMvelocityt{\timeIndex+\frac{1}{2}} &= \ROMvelocityt{\timeIndex} - (\timestepk{\timeIndex}/2)
    \ROMKinematicMassMat^{-1} \velocityBasis^T \forceOneObliqProjMat \forceOneApproxk{\timeIndex}, &
    \ROMvelocityt{\timeIndex+1} &= \ROMvelocityt{\timeIndex} - \timestepk{\timeIndex}
    \ROMKinematicMassMat^{-1} \velocityBasis^T \forceOneObliqProjMat \forceOneApproxk{\timeIndex+\frac{1}{2}}, \\
  \ROMenergyt{\timeIndex+\frac{1}{2}} &= \ROMenergyt{\timeIndex} + (\timestepk{\timeIndex}/2)
    \ROMThermodynamicMassMat^{-1}  \energyBasis^T \forceTvObliqProjMat \forceTvApproxk{\timeIndex}, &
    \ROMenergyt{\timeIndex+1} &= \ROMenergyt{\timeIndex} + \timestepk{\timeIndex}
    \ROMThermodynamicMassMat^{-1}  \energyBasis^T \forceTvObliqProjMat \avgforceTvApproxk{\timeIndex+\frac{1}{2}}, \\
  \ROMpositiont{\timeIndex+\frac{1}{2}} &= \ROMpositiont{\timeIndex} + (\timestepk{\timeIndex}/2)
    \positionBasis^T \velocityApproxt{\timeIndex+\frac{1}{2}}, & \ROMpositiont{\timeIndex+1} &=
    \ROMpositiont{\timeIndex} + \timestepk{\timeIndex} \positionBasis^T \avgvelocityApproxt{\timeIndex+\frac{1}{2}},
\end{align}
where the lifted ROM approximation $\stateApproxt{\timeIndex} 
= (\velocityApproxt{\timeIndex}; \energyApproxt{\timeIndex}; 
\positionApproxt{\timeIndex})^T \in\RR{\sizeWholeFE}$ given by 
\begin{align}\label{eq:discretesolrepresentation}
  \velocityApproxt{\timeIndex} &= \velocityOS + \velocityBasis\ROMvelocityt{\timeIndex},&
  \energyApproxt{\timeIndex} &= \energyOS + \energyBasis\ROMenergyt{\timeIndex},&
  \positionApproxt{\timeIndex} &= \positionOS + \positionBasis\ROMpositiont{\timeIndex}, 
\end{align}
is used to compute the updates 
\begin{align}
\forceOneApproxk{\timeIndex} & = \left
(\forceMat (\stateApproxt{\timeIndex}) \right ) \cdot \oneVec, &  
\forceTvApproxk{\timeIndex} & = \left
(\forceMat (\stateApproxt{\timeIndex}) \right )^T \cdot
\velocityApproxt{\timeIndex+\frac{1}{2}}, 
\end{align}
in the first stage. Similarly, $\stateApproxt{\timeIndex+\frac{1}{2}}
= (\velocityApproxt{\timeIndex+\frac{1}{2}}; \energyApproxt{\timeIndex+\frac{1}{2}}; 
\positionApproxt{\timeIndex+\frac{1}{2}})^T \in\RR{\sizeWholeFE}$ 
is used to computed the updates  
\begin{align}
\forceOneApproxk{\timeIndex+\frac{1}{2}} & = \left
(\forceMat (\stateApproxt{\timeIndex+\frac{1}{2}}) \right ) \cdot \oneVec, &
\avgforceTvApproxk{\timeIndex+\frac{1}{2}} & = \left
(\forceMat (\stateApproxt{\timeIndex + \frac{1}{2}}) \right )^T \cdot
\avgvelocityApproxt{\timeIndex+\frac{1}{2}}, 
\end{align}
with $\avgvelocityApproxt{\timeIndex+\frac{1}{2}} =
(\velocityApproxt{\timeIndex} + \velocityApproxt{\timeIndex+1})/2$ in the second
stage.  The lifting is computed only for the sampled degrees of freedom, avoiding
full order computation. Here, the time step size $\timestepk{\timeIndex}$ is determined
adaptively using the automatic time step control algorithm described in
Section~\ref{sec:adaptiveDelta_t}, with the state $\statet{\timeIndex}$ replaced
by the lifted ROM approximation $\stateApproxt{\timeIndex}$.  We remark that the
RK4 scheme in Section~\ref{sec:RK4} can also be applied in a similar manner to
derive the RK4 fully discrete hyper-reduced system.  As noted in
Section~\ref{sec:adaptiveDelta_t}, the temporal discretization depends on the
state space, which implies that it is very likely that the temporal
discretization used in the hyper-reduced system is different from the full
order model even with the same problem setting.  To this end, we denote by
$\ntimestepROM$ the number of time steps in the fully discrete hyper-reduced
system, to differentiate it from the notation $\ntimestep$ for the full order
model.

\subsection{Proper orthogonal decomposition}\label{sec:POD}
This section describes how we obtain the reduced basis matrices, i.e.,
$\velocityBasis$,  $\energyBasis$, $\positionBasis$, $\forceOneBasis$, and
$\forceTvBasis$.  Proper orthogonal decomposition (POD) is commonly used to
construct a reduced basis. It suffices to describe how to construct
the reduced basis for the energy field only, i.e., $\energyBasis$, because other
bases will be constructed in the same way using POD.  POD
\cite{berkooz1993proper} obtains $\energyBasis$ from a truncated singular
value decomposition (SVD) approximation to a FOM solution snapshot matrix. It is
related to principal component analysis in statistical analysis
\cite{hotelling1933analysis} and Karhunen--Lo\`{e}ve expansion \cite{loeve1955}
in stochastic analysis.  In order to collect solution data for performing POD,
we run FOM simulations on a set of problem parameters, namely
$\{\param_{\paramIndex}\}_{\paramIndex=1}^{\nparam}$.  For
$\paramIndex\in\nat{\nparam}$, let $\ntimestep(\param_{\paramIndex})$ be the
number of time steps in the FOM simulation with the problem parameter
$\param_\paramIndex$.  By choosing $\energyOS(\param_\paramIndex) =
\energy(0;\param_\paramIndex)$, a solution snapshot matrix is formed by
assembling all the FOM solution data including the intermediate
Runge-Kutta stages, i.e.
 \begin{equation}
 \snapshots\equiv\bmat{\energyt{1}(\param_1)-\energyOS(\param_1) & \cdots &
 \energyt{\ntimestep^{\param_{\nparam}}}(\param_{\nparam})-\energyOS(\param_{\nparam})} \in
 \RR{\sizeThermodynamicFE\times\sizeSnapshot},
 \end{equation}
where $\energyt{\timeIndex}(\param_{\paramIndex})$ is the energy state at
$\timeIndex$th time step with problem parameter $\param_\paramIndex$ for
$\timeIndex\in\nat{\ntimestep(\param_{\paramIndex})}$ computed from the FOM
simulation, e.g. the fully discrete RK2-average scheme \eqref{eq:RK2-avg}, and
$\sizeSnapshot = \rkStage \sum_{\paramIndex=1}^{\nparam}
\ntimestep(\param_{\paramIndex})$ with $\rkStage$ being the number of 
Runge-Kutta stages in a time step. 
Then, POD computes its thin SVD: 
 \begin{align}\label{eq:SVD} 
   \snapshots &= \leftSingularMat\singularValueMat\rightSingularMat^T,
 \end{align} 
where $\leftSingularMat\in\RR{\sizeThermodynamicFE\times\sizeSnapshot}$ and
$\rightSingularMat\in\RR{\sizeSnapshot\times\sizeSnapshot}$ are orthogonal matrices,
and $\singularValueMat\in\RR{\sizeSnapshot\times\sizeSnapshot}$ is the diagonal
singular value matrix.  Then POD chooses the leading
$\sizeROMenergy$ columns of $\leftSingularMat$ to set
$\energyBasis = \bmat{\leftSingularVeck{1} & \ldots &
\leftSingularVeck{\sizeROMenergy}}$, where $\leftSingularVeck{\basisIndex}$ is
$\basisIndex$-th column vector of $\leftSingularMat$. The basis size,
$\sizeROMenergy$, is determined by the energy criteria, i.e., we find the
minimum $\sizeROMenergy\in\nat{\sizeSnapshot}$ such that the following
condition is satisfied:
 \begin{align}\label{eq:energy_criteria}
   \frac{\sum_{\basisIndex=1}^{\sizeROMenergy}
   \singularValue_{\basisIndex}}{\sum_{\basisIndex=1}^{\sizeSnapshot}
   \singularValue_{\basisIndex}} &\geq \singularValueThreshold,
 \end{align}
where $\singularValue_{\basisIndex}$ is a $\basisIndex$-th largest singular
value in the singular matrix, $\singularValueMat$, and $\singularValueThreshold
\in \RR{}_{+}$ denotes a threshold.\footnote{We use the default value $\singularValueThreshold =
0.9999$ unless stated otherwise.}

The POD basis minimizes $\|\snapshots - \energyBasis\energyBasis^T\snapshots
\|_F^2$ over all $\energyBasis\in\RR{\sizeThermodynamicFE \times
\sizeROMenergy}$ with orthonormal columns, where $\|\dummyMat\|_F$ denotes the
Frobenius norm of a matrix $\dummyMat\in\RR{I\times J}$, defined as
$\|\dummyMat\|_F = \sqrt{\sum_{i=1}^{I}\sum_{j=1}^{J} \dummyElemSymbol_{ij}^2}$
with $\dummyElemSymbol_{ij}$ being the $(i,j)$ element of $\dummyMat$.  Since
the objective function does not change if $\energyBasis$ is post-multiplied by
an arbitrary $\sizeROMenergy\times\sizeROMenergy$ orthogonal matrix, the POD
procedure seeks the optimal $\sizeROMenergy$-dimensional subspace that
captures the snapshots in the least-squares sense.  For more details on POD, we
refer to \cite{hinze2005proper,kunisch2002galerkin}.

The same procedure can be used to construct the other bases $\velocityBasis$,
$\positionBasis$, $\forceOneBasis$, and $\forceTvBasis$.
Then the reduced mass matrices in \eqref{eq:rmass} can be precomputed once and stored.

\subsection{Solution nonlinear subspace}\label{sec:SNS}
While the snapshot SVD can be applied to construct the nonlinear term bases
$\forceOneBasis$ and $\forceTvBasis$, an alternative way to obtain these
basis matrices is to use the solution nonlinear subspace (SNS)
method in \cite{choi2020sns}. The subspace relations $\forceOneBasis =
\kinematicMassMat \velocityBasis$ and $\forceTvBasis = \thermodynamicMassMat
\energyBasis$ are established by Eq.~\eqref{eq:rom-sol}.  In this case, we have 
$\ROMKinematicMassMat = \velocityBasis^T \forceOneBasis$ and 
$\ROMThermodynamicMassMat = \energyBasis^T \forceTvBasis$, and 
therefore \eqref{eq:rom-hr} becomes 
\begin{align}\label{eq:rom-hr-sns}
   \frac{d\ROMvelocity}{d\timeSymbol} &=
    -\ROMforceOne(\velocityOS + \velocityBasis\ROMvelocity, \energyOS + \energyBasis\ROMenergy, \positionOS + \positionBasis\ROMposition, \timeSymbol; \param) \\
    \frac{d\ROMenergy}{d\timeSymbol} &=
    \ROMforceTv(\velocityOS + \velocityBasis\ROMvelocity, \energyOS + \energyBasis\ROMenergy, \positionOS + \positionBasis\ROMposition, \timeSymbol; \param) \\
    \frac{d\ROMposition}{d\timeSymbol} &= \positionBasis^T \velocityOS +
    \positionBasis^T \velocityBasis \ROMvelocity,
\end{align}
where $\ROMforceOne$ and $\ROMforceTv$ are defined in \eqref{eq:sol-ls-hyper}
and \eqref{eq:sol-ls-hyper-2} respectively. 

\subsection{Discrete Empirical Interpolation Method}\label{sec:DEIM}
This section describes how we obtain the sampling matrices, i.e.
$\forceOneSamplingMat$ and $\forceTvSamplingMat$.  The discrete empirical
interpolation method (DEIM) is a popular choice for nonlinear model reduction.
It suffices to describe how to construct the sampling matrix for the
momentum nonlinear term only, i.e., $\forceOneSamplingMat$, as the other
matrix will be constructed in the same way.  The sampling matrix
$\forceOneSamplingMat$ is characterized by the sampling indices
$\{p_1,\ldots,p_{\sizeROMforceOneSample}\}$, which can be found either by a row
pivoted LU decomposition \cite{chaturantabut2010nonlinear} or the strong column
pivoted rank-revealing QR (sRRQR) decomposition \cite{drmac2016new,
drmac2018discrete}.  Algorithm 1 of \cite{chaturantabut2010nonlinear} uses the
greedy algorithm to sequentially seek additional interpolating indices
corresponding to the entry with the largest magnitude of the residual of
projecting an active POD basis vector onto the preceding basis vectors at the
preceding interpolating indices.  The number of interpolating indices
returned is the same as the number of basis vectors, i.e. $\sizeROMforceOneSample =
\sizeROMforceOne$.  Algorithm 3 of \cite{carlberg2013gnat} and Algorithm 5 of
\cite{carlberg2011efficient} use the greedy procedure to minimize the error in
the gappy reconstruction of the POD basis vectors $\forceOneBasis$.  These
algorithms allow over-sampling, i.e. $\sizeROMforceOneSample \geq
\sizeROMforceOne$, and can be regarded as extensions of Algorithm 1 of
\cite{chaturantabut2010nonlinear}.  Instead of using the greedy algorithm,
Q-DEIM is introduced in \cite{drmac2016new} as a new framework for constructing
the DEIM projection operator via the QR factorization with column pivoting.
Depending on the algorithm for selecting the sampling indices, the DEIM
projection error bound is determined. For example, the row pivoted LU
decomposition in \cite{chaturantabut2010nonlinear} results in the following
error bound:
 \begin{align}\label{eq:errorbound_DEIM}
   \|\forceOne - \forceOneObliqProjMat \forceOne \|_2 \leq \conditionnumber
   \|(\identityMat{\sizeKinematicFE}-\forceOneBasis\forceOneBasis^T)\forceOne \|_2,
 \end{align}
where $\|\cdot\|_2$ denotes $\ell_2$ norm of a vector, and $\conditionnumber$ is
the condition number of $(\forceOneSamplingMat^T\forceTvBasis)^{-1}$, bounded by
 \begin{align}\label{eq:crudeBound}
   \conditionnumber \leq
   (1+\sqrt{2\sizeKinematicFE})^{\sizeROMforceOne-1}\|\forceOneBasisVeck{1}\|_\infty^{-1}.
 \end{align}
On the other hand, the sRRQR factorization in \cite{drmac2018discrete} reveals a
tighter bound than \eqref{eq:crudeBound}:
 \begin{align}\label{eq:tighterBound}
   \conditionnumber \leq
   \sqrt{1+\tuningparam^2\sizeROMforceOne(\sizeKinematicFE-\sizeROMforceOne)}
 \end{align}
where $\tuningparam$ is a tuning parameter in the sRRQR factorization (i.e.,
$f$ in Algorithm 4 of \cite{gu1996efficient}). 
Once the sampling indices are determined, the reduced matrix $\left(\forceOneSamplingMat^T\forceOneBasis\right)^{\dagger}
  \in \RR{\sizeROMforceOne\times\sizeROMforceOneSample}$ 
 can be precomputed and stored.

\section{Time-windowing approach}\label{sec:timewindowing}
Section~\ref{sec:ROM} presented a spatial ROM, in which the spatial unknowns are
reduced to small subspaces. For time-dependent advection-dominated problems, the
dimensions of the ROM spaces required to maintain accuracy grow with the
simulation time. This is intuitively clear, considering that the solution
changes over time, across many elements in the domain, resulting in many linearly
independent snapshots. Consequently, the computational cost of the ROM
simulation grows with the simulation time, so that straightforward use of ROM
may not be faster than the FOM simulation. In order to ensure
adequate speed-up for the ROM simulation, in this section we introduce a
time-windowing approach that decomposes the time domain $[\timek{0},\finalTime]$
into small windows with temporally-local ROM spaces defined on each window, such
that each ROM space is small but accurate for its window.  The time windows can
be prescribed for some given times $\windowk{\windowIndex}$, i.e. window
$\windowIndex$ is $[\windowk{\windowIndex-1},\windowk{\windowIndex}]$ for $1
\leq \windowIndex< \nwindow$, with $\windowk{0} = 0$ and $\windowk{\nwindow} =
\finalTime$.  The time integration should solve for the windows
$1 \leq \windowIndex \leq \nwindow$ sequentially.  First, we
discuss the solution representation and the hyper-reduced system in a time window
in Sections~\ref{sec:solrepresentation-tw} and \ref{sec:hyperreduction-tw}, respectively.
We complete the discussion of the online phase with the initial conditions in a time window in
Section~\ref{sec:init-tw}. Then we move on to discuss the offline phase in the time windowing approach.
We present the mechanisms of decomposing the temporal domain into time windows
in Section~\ref{sec:decompose_tw} and the details
of construction of the temporally local solution subspaces by proper orthogonal
decomposition on snapshot matrices in Section~\ref{sec:POD_tw}.
The sampling indices in a time window can be constructed by DEIM,
described in Section~\ref{sec:DEIM}.
Finally, we discuss the choices of offset vectors in each
time window and the corresponding online computation in Section~\ref{sec:offset_tw}.

\subsection{Solution representation}\label{sec:solrepresentation-tw}
In the time window $\windowIndex$, suppose a basis is 
constructed for each variable and nonlinear term by using FOM snapshots from
that window. Then we restrict our solution space 
to a subspace spanned by a reduced basis for
each field. That is, the subspace for velocity, energy, and position fields are
defined as 
\begin{align}\label{eq:subspces-tw}
  \velocitySubspaceWindow{\windowIndex} &\equiv
  \Span{\velocityBasisVeck{\windowIndex,\basisIndex}}_{\basisIndex=1}^{\sizeROMvelocityWindow{\windowIndex}} 
  \subseteq \RR{\sizeKinematicFE},&
  \energySubspaceWindow{\windowIndex} &\equiv
  \Span{\energyBasisVeck{\windowIndex,\basisIndex}}_{\basisIndex=1}^{\sizeROMenergyWindow{\windowIndex}} 
  \subseteq \RR{\sizeThermodynamicFE},&
  \positionSubspaceWindow{\windowIndex} &\equiv
  \Span{\positionBasisVeck{\windowIndex,\basisIndex}}_{\basisIndex=1}^{\sizeROMpositionWindow{\windowIndex}} 
  \subseteq \RR{\sizeKinematicFE}, 
\end{align}
with $\dim(\velocitySubspaceWindow{\windowIndex})=\sizeROMvelocityWindow{\windowIndex} \ll\sizeKinematicFE$,
$\dim(\energySubspaceWindow{\windowIndex})=\sizeROMenergyWindow{\windowIndex} \ll\sizeThermodynamicFE$,  and
$\dim(\positionSubspaceWindow{\windowIndex})=\sizeROMpositionWindow{\windowIndex} \ll\sizeKinematicFE$. Using these subspaces, each discrete field is approximated in trial subspaces, 
$\velocity\approx \velocityApprox \in \velocityOSWindow{\windowIndex} + 
\velocitySubspaceWindow{\windowIndex}$, 
$\energy \approx\energyApprox \in \energyOSWindow{\windowIndex} + 
\energySubspaceWindow{\windowIndex}$, and 
$\position \approx\positionApprox \in \positionOSWindow{\windowIndex} + 
\positionSubspaceWindow{\windowIndex}$, or equivalently
\begin{align}\label{eq:solrepresentation-tw}
  \velocityApprox(\timeSymbol; \param) &= \velocityOSWindow{\windowIndex}(\param) + 
  \velocityBasisWindow{\windowIndex}\ROMvelocityWindow{\windowIndex}(\timeSymbol; \param), \\
  \energyApprox(\timeSymbol; \param) &= \energyOSWindow{\windowIndex}(\param) + 
  \energyBasisWindow{\windowIndex}\ROMenergyWindow{\windowIndex}(\timeSymbol; \param), \\
  \positionApprox(\timeSymbol; \param) &= \positionOSWindow{\windowIndex}(\param) + 
  \positionBasisWindow{\windowIndex}\ROMpositionWindow{\windowIndex}(\timeSymbol; \param), 
\end{align}
where 
$\velocityOSWindow{\windowIndex}(\param)\in\RR{\sizeKinematicFE}$, 
$\energyOSWindow{\windowIndex}(\param)\in\RR{\sizeThermodynamicFE}$, and
$\positionOSWindow{\windowIndex}(\param)\in\RR{\sizeKinematicFE}$ 
denote the prescribed offset vectors for velocity, energy, and position fields respectively, 
which will be discussed in detail in Section~\ref{sec:offset_tw}. 
The basis matrices 
$\velocityBasisWindow{\windowIndex}\in\RR{\sizeKinematicFE\times\sizeROMvelocityWindow{\windowIndex}}$,
$\energyBasisWindow{\windowIndex}\in\RR{\sizeThermodynamicFE\times\sizeROMenergyWindow{\windowIndex}}$, and
$\positionBasisWindow{\windowIndex}\in\RR{\sizeKinematicFE\times\sizeROMpositionWindow{\windowIndex}}$
are defined as
\begin{align}\label{eq:basisMats-tw}
  \velocityBasisWindow{\windowIndex} &\equiv \bmat{\velocityBasisVeck{\windowIndex,1} & \cdots & \velocityBasisVeck{\windowIndex,\sizeROMvelocityWindow{\windowIndex}}},&
  \energyBasisWindow{\windowIndex} &\equiv \bmat{\energyBasisVeck{\windowIndex,1} & \cdots & \energyBasisVeck{\windowIndex,\sizeROMenergyWindow{\windowIndex}}},&
  \positionBasisWindow{\windowIndex} &\equiv \bmat{\positionBasisVeck{\windowIndex,1} & \cdots & \positionBasisVeck{\windowIndex,\sizeROMpositionWindow{\windowIndex}}},
\end{align}
and 
$\ROMvelocityWindow{\windowIndex}:[\windowk{\windowIndex-1},\windowk{\windowIndex}]\times \paramDomain \mapto \RR{\sizeROMvelocityWindow{\windowIndex}}$,
$\ROMenergyWindow{\windowIndex}:[\windowk{\windowIndex-1},\windowk{\windowIndex}]\times \paramDomain \mapto \RR{\sizeROMenergyWindow{\windowIndex}}$, and
$\ROMpositionWindow{\windowIndex}:[\windowk{\windowIndex-1},\windowk{\windowIndex}]\times \paramDomain \mapto \RR{\sizeROMpositionWindow{\windowIndex}}$
denote the time-dependent generalized coordinates for velocity, energy, and
position fields, respectively. Replacing $\velocity$, $\energy$, and $\position$ with
$\velocityApprox$, $\energyApprox$, and $\positionApprox$ in
Eqs.~\eqref{eq:fom}, the semi-discretized Lagrangian hydrodynamics system 
in the time window $[\windowk{\windowIndex-1},\windowk{\windowIndex}]$
becomes the following over-determined semi-discrete system:
\begin{align}\label{eq:rom-sol-tw}
   \kinematicMassMat \velocityBasisWindow{\windowIndex} 
   \frac{d\ROMvelocityWindow{\windowIndex}}{d\timeSymbol} &=
    - \forceOne(\velocityOSWindow{\windowIndex} + \velocityBasisWindow{\windowIndex}\ROMvelocityWindow{\windowIndex}, \energyOSWindow{\windowIndex} + \energyBasisWindow{\windowIndex}\ROMenergyWindow{\windowIndex}, \positionOSWindow{\windowIndex} + \positionBasisWindow{\windowIndex}\ROMpositionWindow{\windowIndex}, \timeSymbol; \param) \\
    \thermodynamicMassMat \energyBasisWindow{\windowIndex} \frac{d\ROMenergyWindow{\windowIndex}}{d\timeSymbol} &=
    \forceTv(\velocityOSWindow{\windowIndex} + \velocityBasisWindow{\windowIndex}\ROMvelocityWindow{\windowIndex}, \energyOSWindow{\windowIndex} + \energyBasisWindow{\windowIndex}\ROMenergyWindow{\windowIndex}, \positionOSWindow{\windowIndex} + \positionBasisWindow{\windowIndex}\ROMpositionWindow{\windowIndex}, \timeSymbol; \param) \\
    \positionBasisWindow{\windowIndex} \frac{d\ROMpositionWindow{\windowIndex}}{d\timeSymbol} &= (\positionBasisWindow{\windowIndex})^T \velocityOSWindow{\windowIndex} +
    \velocityBasisWindow{\windowIndex} \ROMvelocityWindow{\windowIndex}.
\end{align}
The system \eqref{eq:rom-sol-tw} of ordinary differential equations is closed by defining the initial
condition at $\timeSymbol = \windowk{\windowIndex-1}$, 
which will be discussed in Section~\ref{sec:init-tw}. 

\subsection{Hyper-reduction}\label{sec:hyperreduction-tw}
Following Section~\ref{sec:Galerkin} and Section~\ref{sec:hyperreduction}, 
we derive the hyper-reduced system in a time window $\windowIndex$.
We approximate $\forceOne$ and $\forceTv$ as
 \begin{align}\label{eq:DEIM_approx-tw}
   \forceOne &\approx \forceOneBasisWindow{\windowIndex} \ROMforceOneWindow{\windowIndex},& 
   \forceTv  &\approx \forceTvBasisWindow{\windowIndex} \ROMforceTvWindow{\windowIndex}.
 \end{align}
That is, $\forceOne$ and $\forceTv$ are projected onto subspaces
$\forceOneSubspaceWindow{\windowIndex} \equiv
\Span{\forceOneBasisVeck{\windowIndex,\basisIndex}}_{\basisIndex=1}^{\sizeROMforceOneWindow{\windowIndex}}$ and
$\forceTvSubspaceWindow{\windowIndex} \equiv
\Span{\forceTvBasisVeck{\windowIndex,\basisIndex}}_{\basisIndex=1}^{\sizeROMforceTvWindow{\windowIndex}}$,
 where $\forceOneBasisWindow{\windowIndex} \equiv \bmat{\forceOneBasisVeck{\windowIndex,1} & \ldots &
 \forceOneBasisVeck{\windowIndex,\sizeROMforceOneWindow{\windowIndex}}} \in
 \RR{\sizeKinematicFE\times\sizeROMforceOneWindow{\windowIndex}}$, 
 $\sizeROMforceOneWindow{\windowIndex} \ll \sizeKinematicFE$ and 
 $\forceTvBasisWindow{\windowIndex} \equiv \bmat{\forceTvBasisVeck{\windowIndex,1} &
 \ldots & \forceTvBasisVeck{\windowIndex,\sizeROMforceTvWindow{\windowIndex}}} \in
 \RR{\sizeThermodynamicFE\times\sizeROMforceTvWindow{\windowIndex}}$, 
 $\sizeROMforceTvWindow{\windowIndex} \ll \sizeThermodynamicFE$ denote the nonlinear term basis matrices and
 $\ROMforceOneWindow{\windowIndex} \in \RR{\sizeROMforceOneWindow{\windowIndex}}$ and $\ROMforceTvWindow{\windowIndex} \in
 \RR{\sizeROMforceTvWindow{\windowIndex}}$ denote the generalized coordinates of the nonlinear
 terms. Then, the hyper-reduced system for $[\windowk{\windowIndex-1},\windowk{\windowIndex}]$ is given by
\begin{align}\label{eq:rom-hr-tw}
   \ROMKinematicMassMatWindow{\windowIndex} \frac{d\ROMvelocityWindow{\windowIndex}}{d\timeSymbol} &=
    -(\velocityBasisWindow{\windowIndex})^T \forceOneObliqProjMatWindow{\windowIndex} \forceOne(\velocityOSWindow{\windowIndex} + \velocityBasisWindow{\windowIndex}\ROMvelocityWindow{\windowIndex}, \energyOSWindow{\windowIndex} + \energyBasisWindow{\windowIndex}\ROMenergyWindow{\windowIndex}, \positionOSWindow{\windowIndex} + \positionBasisWindow{\windowIndex}\ROMpositionWindow{\windowIndex}, \timeSymbol; \param) \\
    \ROMThermodynamicMassMatWindow{\windowIndex}\frac{d\ROMenergyWindow{\windowIndex}}{d\timeSymbol} &=
    (\energyBasisWindow{\windowIndex})^T \forceTvObliqProjMatWindow{\windowIndex} \forceTv(\velocityOSWindow{\windowIndex} + \velocityBasisWindow{\windowIndex}\ROMvelocityWindow{\windowIndex}, \energyOSWindow{\windowIndex} + \energyBasisWindow{\windowIndex}\ROMenergyWindow{\windowIndex}, \positionOSWindow{\windowIndex} + \positionBasisWindow{\windowIndex}\ROMpositionWindow{\windowIndex}, \timeSymbol; \param) \\
    \frac{d\ROMpositionWindow{\windowIndex}}{d\timeSymbol} &= (\positionBasisWindow{\windowIndex})^T \velocityOSWindow{\windowIndex} +
    (\positionBasisWindow{\windowIndex})^T \velocityBasisWindow{\windowIndex} \ROMvelocityWindow{\windowIndex},
\end{align}
where the reduced kinematic matrix, $\ROMKinematicMassMatWindow{\windowIndex}
\in \RR{\sizeROMvelocityWindow{\windowIndex} \times
\sizeROMvelocityWindow{\windowIndex}}$, and thermodynamic mass matrix,
$\ROMThermodynamicMassMatWindow{\windowIndex} \in
\RR{\sizeROMenergyWindow{\windowIndex} \times
\sizeROMenergyWindow{\windowIndex}}$, are defined as
\begin{align}\label{eq:rmass-tw}
  \ROMKinematicMassMatWindow{\windowIndex} & = (\velocityBasisWindow{\windowIndex})^T \kinematicMassMat \velocityBasisWindow{\windowIndex} ,&
  \ROMThermodynamicMassMatWindow{\windowIndex} & = (\energyBasisWindow{\windowIndex})^T \thermodynamicMassMat \energyBasisWindow{\windowIndex},
\end{align}
 the oblique projection matrices $\forceOneObliqProjMatWindow{\windowIndex} \in \RR{\sizeKinematicFE\times\sizeKinematicFE}$ and 
$\forceTvObliqProjMatWindow{\windowIndex} \in 
\RR{\sizeThermodynamicFE\times\sizeThermodynamicFE}$ are defined as 
\begin{align}\label{eq:rhsproj-tw}
\forceOneObliqProjMatWindow{\windowIndex} & = \forceOneBasisWindow{\windowIndex} 
 ((\forceOneSamplingMatWindow{\windowIndex} )^T\forceOneBasisWindow{\windowIndex} )^{\dagger}(\forceOneSamplingMatWindow{\windowIndex})^T, &
    \forceTvObliqProjMat & = \forceTvBasisWindow{\windowIndex}  \left((\forceTvSamplingMatWindow{\windowIndex})^T
      \forceTvBasisWindow{\windowIndex} \right)^{\dagger} (\forceTvSamplingMatWindow{\windowIndex})^T, 
\end{align}
with the sampling matrices 
$\forceOneSamplingMatWindow{\windowIndex} \in \RR{\sizeKinematicFE \times
 \sizeROMforceOneSampleWindow{\windowIndex}}$ and 
$\forceTvSamplingMatWindow{\windowIndex} \in \RR{\sizeThermodynamicFE \times
 \sizeROMforceTvSampleWindow{\windowIndex}}$. 
As in the case of spatial ROM, the reduced mass matrices
$\ROMKinematicMassMatWindow{\windowIndex}$ and 
$\ROMThermodynamicMassMatWindow{\windowIndex}$ and 
the reduced matrices $\forceOneBasisWindow{\windowIndex} 
 ((\forceOneSamplingMatWindow{\windowIndex} )^T\forceOneBasisWindow{\windowIndex} )^{\dagger}$ 
and $\forceTvBasisWindow{\windowIndex}  \left((\forceTvSamplingMatWindow{\windowIndex})^T
      \forceTvBasisWindow{\windowIndex} \right)^{\dagger}$ can be precomputed and stored.

\subsection{Initial condition of window}\label{sec:init-tw}
In this subsection, we discuss the initial conditions of 
the variables in a time window. 
The corresponding fully discretized system can be obtained by the time
integrators introduced in Section~\ref{sec:rom_time_integrator}. 
Using the relation \eqref{eq:solrepresentation-tw}, one can derive the initial condition by
lifting the ROM solution to the FOM spaces, using the ROM bases in the time
window $\windowIndex-1$, and then projecting onto the ROM spaces in the time window $\windowIndex$, i.e.
\begin{align}\label{eq:ic-tw}
\ROMvelocityWindow{\windowIndex}(\windowk{\windowIndex-1}) & = 
(\velocityBasisWindow{\windowIndex})^T ( 
\velocityOSWindow{\windowIndex-1} + 
  \velocityBasisWindow{\windowIndex-1}\ROMvelocityWindow{\windowIndex-1}(\windowk{\windowIndex-1}) - \velocityOSWindow{\windowIndex} ), \\
\ROMenergyWindow{\windowIndex}(\windowk{\windowIndex-1}) & = 
(\energyBasisWindow{\windowIndex})^T ( 
\energyOSWindow{\windowIndex-1} + 
  \energyBasisWindow{\windowIndex-1}\ROMenergyWindow{\windowIndex-1}(\windowk{\windowIndex-1}) - \energyOSWindow{\windowIndex} ), \\
\ROMpositionWindow{\windowIndex}(\windowk{\windowIndex-1}) & = 
(\positionBasisWindow{\windowIndex})^T ( 
\positionOSWindow{\windowIndex-1} + 
  \positionBasisWindow{\windowIndex-1}\ROMpositionWindow{\windowIndex-1}(\windowk{\windowIndex-1}) - \positionOSWindow{\windowIndex} ).
\end{align}
However, in view of Theorem~\ref{thm:apriori-semi}, in order to minimize the
error bound in the induced norm, we can project the lifted solution obliquely
onto the ROM spaces by 
\begin{align}\label{eq:ic2-tw}
\ROMvelocityWindow{\windowIndex}(\windowk{\windowIndex-1}) & = 
(\ROMKinematicMassMatWindow{\windowIndex})^{-1} (\velocityBasisWindow{\windowIndex})^T \kinematicMassMat ( 
\velocityOSWindow{\windowIndex-1} + 
  \velocityBasisWindow{\windowIndex-1}\ROMvelocityWindow{\windowIndex-1}(\windowk{\windowIndex-1}) - \velocityOSWindow{\windowIndex} ), \\
\ROMenergyWindow{\windowIndex}(\windowk{\windowIndex-1}) & = 
(\ROMThermodynamicMassMatWindow{\windowIndex})^{-1} (\energyBasisWindow{\windowIndex})^T \thermodynamicMassMat ( 
\energyOSWindow{\windowIndex-1} + 
  \energyBasisWindow{\windowIndex-1}\ROMenergyWindow{\windowIndex-1}(\windowk{\windowIndex-1}) - \energyOSWindow{\windowIndex} ), \\
\ROMpositionWindow{\windowIndex}(\windowk{\windowIndex-1}) & = 
(\positionBasisWindow{\windowIndex})^T ( 
\positionOSWindow{\windowIndex-1} + 
  \positionBasisWindow{\windowIndex-1}\ROMpositionWindow{\windowIndex-1}(\windowk{\windowIndex-1}) - \positionOSWindow{\windowIndex} ). 
\end{align}
We remark that for either case, the determination of the initial condition 
only involves inexpensive operations in reduced dimensions in the online case, 
since the contribution of offset vectors can be precomputed and stored in the offline phase. 

\subsection{Mechanism of decomposing time domain}\label{sec:decompose_tw}
The remainder of this section is devoted to discussing the
offline computations of the reduced order model.
It should be noted that the time windows are set
before constructing the reduced bases in the offline phase,
which will be discussed in Section~\ref{sec:POD_tw}.
We consider two approaches of decomposing the time domain into windows.

\subsubsection{Physical time windowing}\label{sec:physical_tw}
A natural approach of decomposing the domain into windows is to use a fixed
partition, that is, the end points $\{ \windowk{\windowIndex}
\}_{\windowIndex=1}^{\nwindow}$ are fixed and user-defined.  Using this naive
approach, the number of snapshots for a variable in a window is proportional to
the window size if a uniform time step is used in the FOM simulation.  In
addition, if a uniform partition is used to decompose the domain into windows,
then the basis sizes will be heuristically balanced among the windows.  However,
in the setting of adaptive time stepping, it becomes unclear how to
prescribe window sizes to balance the ROM basis size among time windows. 

\subsubsection{Time windowing by number of samples}\label{sec:sample_tw}
We consider an alternative way to divide the time domain into windows using a
prescribed number of snapshots per window, which can be heuristically related to
the basis size.  Let $\ntimestepWindow$ be the number of FOM solution 
snapshots per window.  Recall that
for $\paramIndex\in\nat{\nparam}$, $\ntimestep(\param_{\paramIndex})$ is the
number of time steps in the FOM simulation with the problem parameter
$\param_\paramIndex$, and the temporal domain is discretized as
$\{\timek{\timeIndex}(\param_{\paramIndex}) \}
_{\timeIndex=0}^{\ntimestep(\param_{\paramIndex})}$.  We determine the end
points of the time windows sequentially.  Given the end point
$\windowk{\windowIndex-1}$ of the time window $\windowIndex-1$, we denote by
$\timeIndexWindow{\windowIndex-1}(\param_{\paramIndex})$ the latest time instance
in the FOM simulation with the problem parameter $\param_\paramIndex$ which lies
in the time window $\windowIndex-1$, i.e.
$\timek{\timeIndexWindow{\windowIndex-1}(\param_{\paramIndex})}(\param_{\paramIndex})
< \windowk{\windowIndex-1} \leq \timek{\timeIndexWindow{\windowIndex-1}
(\param_{\paramIndex})+1}(\param_{\paramIndex})$.  Then we determine the end
point of the time window $\windowIndex$ as 
\begin{equation}\label{eq:sample_window}
\windowk{\windowIndex} = \min\left\{ \timek{\timeIndexWindow{\windowIndex-1}
  (\param_{\paramIndex})+\ntimestepWindow+1}(\param_{\paramIndex}):
  \paramIndex\in\nat{\nparam} \right\}, 
\end{equation}
i.e. the shortest time among all the problem parameters $\{ \param_\paramIndex
\}_{\paramIndex=1}^{\nparam}$ such that $\ntimestepWindow+1$ new time steps are
taken.  Here, we extend the notation $\timek{\timeIndex}(\param_\paramIndex)$ to
denote $\timek{\timeIndex}(\param_\paramIndex) = \finalTime$ for all $\timeIndex
\geq \ntimestep(\param_\paramIndex)$.  The iterative process is terminated when
$\windowk{\windowIndex} = \finalTime$, at which we set $\nwindow =
\windowIndex$. 

\subsection{Temporally local solution subspaces}\label{sec:POD_tw}
In this subsection, we discuss the construction of the solution subspaces in a
time window.  Again, it will be sufficient to discuss how to construct the
reduced basis for the energy field only, i.e., $\energyBasisWindow{\windowIndex}
\in \RR{\sizeThermodynamicFE\times\sizeROMenergyWindow{\windowIndex}}$, because
other bases will be constructed in the same way.  In order to collect solution
data for performing POD, we run FOM simulations on a set of problem parameters,
namely $\{\param_{\paramIndex}\}_{\paramIndex=1}^{\nparam}$.
Using the mechanisms described in
Section~\ref{sec:physical_tw} or Section~\ref{sec:sample_tw}, the temporal
domain is decomposed into time windows.  The collected FOM solution vectors are
then clustered into time windows according to the time instance when the sample
was taken.  By choosing $\energyOSWindow{\windowIndex}(\param_\paramIndex)$
appropriately,  a solution snapshot matrix is formed by assembling all the FOM
solution data, i.e.
\begin{equation}\label{eq:window_sample}
\snapshots^\windowIndex\equiv\bmat{\energyt{\timeIndexWindow{\windowIndex-1}(\param_1)}(\param_{\paramIndex})-\energyOSWindow{\windowIndex}(\param_1) & \cdots &
 \energyt{\timeIndexWindow{\windowIndex}(\param_{\nparam})+1}(\param_{\nparam})-\energyOSWindow{\windowIndex}(\param_{\nparam})}  \in
 \RR{\sizeThermodynamicFE\times\sizeSnapshotS^{\windowIndex}},
\end{equation}
where $\sizeSnapshotS^{\windowIndex} = \sum_{\paramIndex=1}^{\nparam}
\left[ \rkStage \left(\timeIndexWindow{\windowIndex}(\param_{\paramIndex})-
\timeIndexWindow{\windowIndex-1}(\param_{\paramIndex})\right)+1\right]$, and then
its thin SVD is computed to obtain $\energyBasisWindow{\windowIndex}$ as in
Section~\ref{sec:POD}. Similar to the discussion in Section~\ref{sec:POD} and Section~\ref{sec:SNS} 
for the spatial ROM, the nonlinear term bases $\forceOneBasisWindow{\windowIndex}$ and 
$\forceTvBasisWindow{\windowIndex}$ in the time window $\windowIndex$ can
also be obtained either by POD on nonlinear term snapshots or the SNS method in
\cite{choi2020sns} with the subspace relations
$\forceOneBasisWindow{\windowIndex} = \kinematicMassMat
\velocityBasisWindow{\windowIndex}$ and $\forceTvBasisWindow{\windowIndex} =
\thermodynamicMassMat \energyBasisWindow{\windowIndex}$. 

\subsection{Window offset vectors}\label{sec:offset_tw}
To end this section, we discuss the choices of the offset vectors in each of the
time windows.  Again, we present how to construct the offset vectors for the
energy field only, i.e. $\energyOSWindow{\windowIndex}$, as the velocity and position
fields are treated similarly.
We remark that, in obtaining the reduced basis matrices from the problem
parameters $\{ \param_{\paramIndex} \}_{\paramIndex=1}^{\nparam}$, the offset
vector $\energyOSWindow{\windowIndex}(\param_\paramIndex)$ has to be determined
and subtracted from the FOM solution vectors in \eqref{eq:window_sample}.  On
the other hand, in employing the reduced order model \eqref{eq:rom-hr-tw} for a
generic problem parameter $\param \in \paramDomain$, the offset vector
$\energyOSWindow{\windowIndex}(\param)$ with the same physical meaning has to be
calculated or approximated. 

\subsubsection{Initial states}\label{sec:offset_initial}
The first choice is to use the initial states as in Section~\ref{sec:ROM}, i.e.
we set $\energyOSWindow{\windowIndex}(\param_\paramIndex) = \energy(0;
\param_\paramIndex)$ for all $\windowIndex \in \nat{\nwindow}$ and  $\paramIndex
\in \nat{\nparam}$.  Then, for a generic problem parameter $\param \in
\paramDomain$, we simply take $\energyOSWindow{\windowIndex}(\param) =
\energy(0; \param)$ for all $\windowIndex \in \nat{\nwindow}$. 

\subsubsection{Previous window} \label{sec:offset_previous}
Another approach is to choose the offset vectors 
depending on the window.  A natural choice is to choose
the offset as an approximate solution around the time
$\windowk{\windowIndex-1}$, i.e.
$\energyOSWindow{\windowIndex}(\param_\paramIndex) =
\energyt{\timeIndexWindow{\windowIndex-1}(\param_{\paramIndex})}(\param_{\paramIndex})$
for all $\windowIndex \in \nat{\nwindow}$ and $\paramIndex \in \nat{\nparam}$.
In that case, in employing the reduced order model in the time window $\windowIndex$ for a
generic problem parameter $\param \in \paramDomain$, an approximation of the
energy field $\energy(\windowk{\windowIndex-1}; \param)$ at the time
$\timeSymbol = \windowk{\windowIndex-1}$ can be used as the offset vector
$\energyOSWindow{\windowIndex}(\param)$.  One such approximation is given by
using the final state in the previous window.  More precisely, we lift the ROM
solution in the time window $\windowIndex-1$ to the FOM spaces by using the ROM
bases, i.e. 
\begin{equation}
\energyOSWindow{\windowIndex}(\param) = \energyOSWindow{\windowIndex-1}(\param) + 
  \energyBasisWindow{\windowIndex-1}\ROMenergyWindow{\windowIndex-1}(\windowk{\windowIndex-1}; \param).
\label{eq:previous_os}
\end{equation}
The advantage of this approach is that it provides the best approximation of the
final state in the previous window in the solution subspace of the current
window.  However, the drawback of this approach is that it cannot be precomputed,
and it involves the lifting operation which scales with the dimension of the FOM
state space. In the case of advection-dominated problems, time windows are
typically small and change frequently, so the lifting operation
may limit the overall speed-up. 

\subsubsection{Parametric interpolation}\label{sec:offset_interpolate}
The last choice we present in this section is to make use of the available data
from the problem parameters $\{ \param_\paramIndex
\}_{\paramIndex=1}^{\nparam}$.  We consider the choice of the offset
vector $\energyOSWindow{\windowIndex}(\param_\paramIndex) =
\energyt{\timeIndexWindow{\windowIndex-1}(\param_{\paramIndex})}(\param_{\paramIndex})$
for all $\windowIndex \in \nat{\nwindow}$ and $\paramIndex \in \nat{\nparam}$.
In order to approximate the offset vector
$\energyOSWindow{\windowIndex}(\param)$ for a generic parameter $\param \in
\paramDomain$, we can interpolate the known offset vectors
$\{\energyOSWindow{\windowIndex}(\param_{\paramIndex})\}_{\paramIndex=1}^{\nparam}$,
which are computed and stored in the process of snapshot collection.  Here, we
present an interpolation procedure using the inverse distance weighting (IDW).
Suppose $\paramMetric: \paramDomain \times \paramDomain \to \RR{+}$ is a
distance metric on the domain $\paramDomain$ and $\interpolationPowerSymbol \geq
1$ is a given number. The interpolated value is given by the convex combination 
\begin{equation}
\energyOSWindow{\windowIndex}(\param) = \begin{cases}
\energyOSWindow{\windowIndex}(\param_{\paramIndex}) & 
\text{ if } \paramMetric(\param, \param_{\paramIndex}) = 0 \text{ for some } \paramIndex \in \nat{\nparam}, \\ 
\dfrac{\sum_{\paramIndex=1}^{\nparam} \interpolationWeightSymbol_{\paramIndex}(\param) 
\energyOSWindow{\windowIndex}(\param_{\paramIndex})}
{\sum_{\paramIndex=1}^{\nparam} \interpolationWeightSymbol_{\paramIndex}(\param)}
& \text{ otherwise}, 
\end{cases}
\label{eq:interpolate_os}
\end{equation}
where, for $\paramIndex \in \nat{\nparam}$, the coefficients 
$\interpolationWeightSymbol_{\paramIndex}(\param) \in \RR{+}$ are given by
\begin{equation}
\interpolationWeightSymbol_{\paramIndex}(\param) = 
\dfrac{1}{\paramMetric(\param, \param_{\paramIndex})^\interpolationPowerSymbol}.
\end{equation}
In our work, we use the Euclidean distance as the metric $\paramMetric$ and $r=2$.  
One advantage of the IDW interpolation
scheme is that the offset vectors are exactly the stored interpolating values in
the reproductive cases, which provides the most accurate approximations.
Another advantage is that the interpolation can be precomputed and is
inexpensive as the interpolating coefficients depend solely on the problem
parameters, which typically lie in low-dimensional structures.  However, the IDW
interpolation scheme lacks the ability to extrapolate.  In practice, other
schemes allowing inexpensive and accurate extrapolation can be considered.

\section{Error bounds}\label{sec:errorbound}
In this section, we present error estimates for our proposed reduced order model
in Section~\ref{sec:ROM}.  Due to the use of adaptive time-step control, we
cannot directly compare the temporally discrete full order model solution and
the reduced order approximation. Instead, we will use the continuous-in-time
full order model solution in \eqref{eq:fom} as the reference solution.
Throughout the section, when there is no ambiguity, we drop the time symbol
$\timeSymbol$ and the problem parameter symbol $\param$ for simplifying the
notations.  The error analysis in decomposed into two parts.  The first part
accounts for the approximation error of the reduced order model.  More
specifically, we analyze the error between the continuous-in-time full order
model solution in \eqref{eq:fom} and the continuous-in-time reduced order
approximation in \eqref{eq:rom-hr}, defined by 
\begin{equation}\label{eq:rom-cerror}
  \begin{aligned}
    \continuousVelocityError & = \velocity - \velocityApprox = \velocity -
    (\velocityOS+\velocityBasis\ROMvelocity) \\
    \continuousEnergyError & = \energy - \energyApprox = \energy -
    (\energyOS+\energyBasis\ROMenergy) \\ 
    \continuousPositionError & = \position - \positionApprox = \position -
    (\positionOS+\positionBasis\ROMposition).
  \end{aligned}
\end{equation}
The second part accounts for the truncation error of the temporal
discretization.  More specifically, we analyze the error between the
continuous-in-time reduced order approximation in \eqref{eq:rom-hr} and the
RK2-average fully discrete reduced order approximation in
\eqref{eq:RK2-avg-rom-hr}, defined by 
\begin{equation}\label{eq:rom-derror}
  \begin{aligned}
    \discreteVelocityErrort{\timeIndex} & = \velocityApprox(\timek{\timeIndex})
    - \velocityApproxt{\timeIndex} =
    \velocityBasis(\ROMvelocity(\timek{\timeIndex}) - \ROMvelocityt{\timeIndex})
    \\
    \discreteEnergyErrort{\timeIndex} & = \energyApprox(\timek{\timeIndex}) -
    \energyApproxt{\timeIndex} = \energyBasis(\ROMenergy(\timek{\timeIndex}) -
    \ROMenergyt{\timeIndex}) \\ 
    \discretePositionErrort{\timeIndex} & = \positionApprox(\timek{\timeIndex})
    - \positionApproxt{\timeIndex} =
    \positionBasis(\ROMposition(\timek{\timeIndex}) -
    \ROMpositiont{\timeIndex}).
  \end{aligned}
\end{equation}
We remark that the error analysis can be extended to the time windowing approach
in Section~\ref{sec:timewindowing} by viewing the final solution in the previous
window as the initial condition and stacking up the error in a sequence of time
windows. 

Before we begin with the analysis, we introduce a few tools and notations that
will facilitate our discussion.  Since the mass matrices $\kinematicMassMat$ and
$\thermodynamicMassMat$ are symmetric and positive definite, they possess the
Cholesky factorizations
\begin{equation}\label{eq:cholesky}
\begin{aligned}
\kinematicMassMat & = \kinematicCholeskyMat \kinematicCholeskyMat^T \\
\thermodynamicMassMat & = \thermodynamicCholeskyMat \thermodynamicCholeskyMat^T.
\end{aligned}
\end{equation}
The mass matrices also induce the weighted functional $L^2$ norm 
on the kinematic and thermodynamic finite element space respectively:
\begin{equation}
\begin{aligned}
\velocityInducedNorm{\velocity}^2  & = \velocity \cdot \kinematicMassMat \cdot
  \velocity \\
\energyInducedNorm{\energy}^2 & = \energy \cdot \thermodynamicMassMat \cdot
  \energy.
\end{aligned}
\end{equation}
On the Cartesian product space $\kinematicFE \times \thermodynamicFE \times
\kinematicFE$, we define the product norm 
\begin{equation}\label{eq:fullnorm}
\fullNorm{(\velocity,\energy,\position)}^2 = 
\velocityInducedNorm{\velocity}^2 + \energyInducedNorm{\energy}^2 + \euclideanNorm{\position}^2,
\end{equation}
where $\euclideanNorm{\cdot}$ denotes the standard Euclidean norm. 

\subsection{A-priori estimate for approximation error}
We start with providing an a-priori error estimate between the
continuous-in-time full order model solution in \eqref{eq:fom} and the
continuous-in-time reduced order approximation in \eqref{eq:rom-hr}.  The error
bound depends on the full order model solution and is controlled by a
combination of several quantities, namely the mismatch of initial generalized
coordinates in the reduced subspaces, the oblique projection error of the
solution onto the reduced subspaces over time, and the projection error of the
nonlinear terms in the Euclidean norm. 
\begin{theorem}\label{thm:apriori-semi}
Assume there holds the following Lipchitz continuity conditions for the force
  matrix $\forceMat$: there exists $K_1, K_2 > 0$ such that for any
  $(\velocity, \energy, \position), (\velocity', \energy', \position') \in
  \kinematicFE \times \thermodynamicFE \times \kinematicFE$, 
\begin{equation}\label{eq:LipRHS-1}
  \begin{aligned}
\euclideanNorm{\forceOne(\velocity, \energy, \position)-\forceOne(\velocity', \energy', \position')}
& \leq K_1 \fullNorm{(\velocity-\velocity',\energy-\energy',\position-\position')} \\
\euclideanNorm{ \forceTv(\velocity, \energy, \position) - \forceTv(\velocity', \energy', \position')} 
& \leq K_2 \fullNorm{(\velocity-\velocity',\energy-\energy',\position-\position')}.
  \end{aligned}
\end{equation}
Then there exists a generic constant $C>0$ such that for $\timeSymbol > 0$, we
  have 
\begin{equation}\label{eq:apriori-semi} 
  \begin{aligned}
    \fullNorm{(\continuousVelocityError(\timeSymbol), \continuousEnergyError(\timeSymbol), \continuousPositionError(\timeSymbol))} & \leq 
     C e^{C\timeSymbol}  \Big[ \fullNorm{(\velocityInitialError, \energyInitialError, \positionInitialError)} + 
    \max_{0 \leq \dummyTimeSymbol \leq \timeSymbol} \fullNorm{(\velocityDifference(\dummyTimeSymbol), \energyDifference(\dummyTimeSymbol), \positionDifference(\dummyTimeSymbol)} + \\
   & \qquad \int_0^\timeSymbol 
\euclideanNorm{(\velocityIdentity - \forceOneObliqProjMat)\forceOne(\velocity, \energy, \position)} +
\euclideanNorm{((\energyIdentity - \forceTvObliqProjMat) \forceTv(\velocity, \energy, \position)} d\dummyTimeSymbol  \Big],
  \end{aligned}
\end{equation}
where $\continuousVelocityError(\timeSymbol)$,
  $\continuousEnergyError(\timeSymbol)$, and
  $\continuousPositionError(\timeSymbol)$ are defined in
  Eq.~\eqref{eq:rom-cerror} and 
\begin{align}\label{eq:def-diff}
    \velocityInitialError & = \velocityBasis \left[\ROMvelocity(0) - \ROMKinematicMassMat^{-1}\velocityBasis^T \kinematicMassMat (\velocity(0)-\velocityOS)\right] & 
    \velocityDifference(\timeSymbol) & = (\velocityIdentity - \velocityBasis \ROMKinematicMassMat^{-1}\velocityBasis^T \kinematicMassMat) (\velocity(\timeSymbol)-\velocityOS) \\
    \energyInitialError & = \energyBasis \left[\ROMenergy(0) - \ROMThermodynamicMassMat^{-1}\energyBasis^T \thermodynamicMassMat (\energy(0)-\energyOS)\right] &
    \energyDifference(\timeSymbol) & = (\energyIdentity - \energyBasis \ROMThermodynamicMassMat^{-1}\energyBasis^T  \thermodynamicMassMat) (\energy(\timeSymbol)-\energyOS) \\
   \positionInitialError & = \positionBasis\left[\ROMposition(0) -  \positionBasis^T (\position(0)-\positionOS)\right] & 
   \positionDifference(\timeSymbol) & = (\positionIdentity - \positionBasis \positionBasis^T) (\position(\timeSymbol)-\positionOS).
\end{align}
\begin{proof}
Since the reduced order mass matrices $\ROMKinematicMassMat$ and
  $\ROMThermodynamicMassMat$ are nonsingular, we rewrite \eqref{eq:rom-hr} as
\begin{equation}\label{eq:rom-hr-re}
  \begin{aligned}  
   \kinematicMassMat \velocityBasis \frac{d\ROMvelocity}{d\timeSymbol} &=
    - \kinematicMassMat \velocityBasis \ROMKinematicMassMat^{-1}\velocityBasis^T 
    \forceOneObliqProjMat \forceOne(\velocityApprox, \energyApprox, \positionApprox) \\
    \thermodynamicMassMat \energyBasis \frac{d\ROMenergy}{d\timeSymbol} &=
    \thermodynamicMassMat \energyBasis \ROMThermodynamicMassMat^{-1} \energyBasis^T 
    \forceTvObliqProjMat \forceTv(\velocityApprox, \energyApprox, \positionApprox) \\
    \positionBasis \frac{d\ROMposition}{d\timeSymbol} &= \positionBasis \positionBasis^T 
    \velocityApprox.
  \end{aligned}
\end{equation}
Noting that $(\velocityOS, \energyOS, \positionOS)$ is constant with respect to
  time and subtracting \eqref{eq:rom-hr-re} from \eqref{eq:fom}, we obtain 
\begin{equation}\label{eq:error-eq}
  \begin{aligned}
   \kinematicMassMat \frac{d\continuousVelocityError}{d\timeSymbol} &
    = \velocityResidual{1}  + \velocityResidual{2} \\
   \thermodynamicMassMat \frac{d\continuousEnergyError}{d\timeSymbol} &
    = \energyResidual{1}  + \energyResidual{2} \\
    \frac{d\continuousPositionError}{d\timeSymbol} &
     = \positionResidual{1}  + \positionResidual{2} \\
  \end{aligned}
\end{equation}
where the residuals are defined by 
\begin{equation}\label{eq:res-def-a}
\begin{aligned}
\velocityResidual{1} & = 
-(\velocityIdentity - \kinematicMassMat \velocityBasis \ROMKinematicMassMat^{-1}\velocityBasis^T)\forceOne(\velocity, \energy, \position) \\
\velocityResidual{2} & = 
-\kinematicMassMat \velocityBasis \ROMKinematicMassMat^{-1}\velocityBasis^T (\forceOne(\velocity, \energy, \position)-
\forceOneObliqProjMat\forceOne(\velocityApprox, \energyApprox, \positionApprox)) \\
\energyResidual{1} & = 
(\energyIdentity -  \thermodynamicMassMat \energyBasis \ROMThermodynamicMassMat^{-1} \energyBasis^T) \forceTv(\velocity, \energy, \position) \\
\energyResidual{2} & = 
\thermodynamicMassMat \energyBasis \ROMThermodynamicMassMat^{-1} \energyBasis^T  (\forceTv(\velocity, \energy, \position) - 
\forceTvObliqProjMat\forceTv(\velocityApprox, \energyApprox, \positionApprox)) \\
\positionResidual{1} & = (\positionIdentity - \positionBasis \positionBasis^T) \velocity \\
\positionResidual{2} & = \positionBasis \positionBasis^T 
(\velocity - \velocityApprox).
\end{aligned}
\end{equation}
Integrating \eqref{eq:error-eq} over $(0,\timeSymbol)$, we have 
\begin{equation}\label{eq:error-eq-int}
  \begin{aligned}
   \kinematicMassMat (\continuousVelocityError(\timeSymbol) - \continuousVelocityError(0)) &=  
   \int_0^\timeSymbol \velocityResidual{1}(\dummyTimeSymbol) d\dummyTimeSymbol + \int_0^\timeSymbol \velocityResidual{2}(\dummyTimeSymbol) d\dummyTimeSymbol \\
   \thermodynamicMassMat (\continuousEnergyError(\timeSymbol) - \continuousEnergyError(0)) &= 
   \int_0^\timeSymbol \energyResidual{1}(\dummyTimeSymbol) d\dummyTimeSymbol + \int_0^\timeSymbol \energyResidual{2}(\dummyTimeSymbol) d\dummyTimeSymbol \\
   \continuousPositionError(\timeSymbol) - \continuousPositionError(0) &= 
   \int_0^\timeSymbol \positionResidual{1}(\dummyTimeSymbol) d\dummyTimeSymbol + \int_0^\timeSymbol \positionResidual{2}(\dummyTimeSymbol) d\dummyTimeSymbol.
  \end{aligned}
\end{equation}
On the other hand, integrating \eqref{eq:fom} over $(0,\timeSymbol)$, we have 
\begin{equation}\label{eq:fom-int}
  \begin{aligned}
    \kinematicMassMat (\velocity(\timeSymbol) - \velocity(0)) &=\int_0^\timeSymbol 
    -\forceOne(\velocity, \energy, \position) d\dummyTimeSymbol \\
    \thermodynamicMassMat (\energy(\timeSymbol) - \energy(0)) &=\int_0^\timeSymbol
    \forceTv(\velocity, \energy, \position) d\dummyTimeSymbol\\
    \position(\timeSymbol) - \position(0) &= \int_0^\timeSymbol \velocity(\dummyTimeSymbol) d\dummyTimeSymbol.
  \end{aligned}
\end{equation}
This implies 
\begin{equation}\label{eq:res1a}
  \begin{aligned}
    \int_0^\timeSymbol \velocityResidual{1}(\dummyTimeSymbol) d\dummyTimeSymbol & = 
    \kinematicMassMat (\velocityIdentity - \velocityBasis \ROMKinematicMassMat^{-1}\velocityBasis^T \kinematicMassMat ) 
    \left[(\velocity(\timeSymbol)-\velocityOS) - (\velocity(0)-\velocityOS)\right] \\
    \int_0^\timeSymbol \energyResidual{1}(\dummyTimeSymbol) d\dummyTimeSymbol & = 
    \thermodynamicMassMat (\energyIdentity - \energyBasis \ROMThermodynamicMassMat^{-1} \energyBasis^T \thermodynamicMassMat) 
    \left[(\energy(\timeSymbol)-\energyOS) - (\energy(0)-\energyOS)\right] \\
    \int_0^\timeSymbol \positionResidual{1}(\dummyTimeSymbol) d\dummyTimeSymbol & = 
    (\positionIdentity - \positionBasis \positionBasis^T) 
    \left[(\position(\timeSymbol)-\positionOS) - (\position(0)-\positionOS)\right].
  \end{aligned}
\end{equation}
Substituting \eqref{eq:res1a} into \eqref{eq:error-eq-int}, we obtain
\begin{equation}\label{eq:error-eq-int2a}
  \begin{aligned}
   \kinematicMassMat \continuousVelocityError(\timeSymbol) &=  
    \kinematicMassMat (\velocityIdentity - \velocityBasis \ROMKinematicMassMat^{-1}\velocityBasis^T \kinematicMassMat) (\velocity(\timeSymbol)-\velocityOS) -
   \kinematicMassMat \velocityBasis \left[\ROMvelocity(0) - \ROMKinematicMassMat^{-1}\velocityBasis^T \kinematicMassMat (\velocity(0)-\velocityOS)\right] + 
   \int_0^\timeSymbol \velocityResidual{2}(\dummyTimeSymbol) d\dummyTimeSymbol \\
   \thermodynamicMassMat \continuousEnergyError(\timeSymbol) &= 
    \thermodynamicMassMat (\energyIdentity - \energyBasis \ROMThermodynamicMassMat^{-1}\energyBasis^T  \thermodynamicMassMat) (\energy(\timeSymbol)-\energyOS) -
   \thermodynamicMassMat \energyBasis \left[\ROMenergy(0) - \ROMThermodynamicMassMat^{-1}\energyBasis^T \thermodynamicMassMat (\energy(0)-\energyOS)\right] + \int_0^\timeSymbol \energyResidual{2}(\dummyTimeSymbol) d\dummyTimeSymbol \\
   \continuousPositionError(\timeSymbol) &= 
    (\positionIdentity - \positionBasis \positionBasis^T) (\position(\timeSymbol)-\positionOS) - \positionBasis\left[\ROMposition(0) -  \positionBasis^T (\position(0)-\positionOS)\right] + \int_0^\timeSymbol \positionResidual{2}(\dummyTimeSymbol) d\dummyTimeSymbol.
  \end{aligned}
\end{equation}
By the triangle inequality, we have 
\begin{equation}\label{eq:error-eq-int3a}
  \begin{aligned}
    \velocityInducedNorm{\continuousVelocityError(\timeSymbol)} & \leq
    \velocityInducedNorm{\velocityInitialError} + 
    \velocityInducedNorm{\velocityDifference(\timeSymbol)} +
    \int_0^\timeSymbol \euclideanNorm{\kinematicCholeskyMat^{-1} \velocityResidual{2}(\dummyTimeSymbol)} d\dummyTimeSymbol \\
    \energyInducedNorm{\continuousEnergyError(\timeSymbol)} & \leq
    \energyInducedNorm{\energyInitialError} + 
    \energyInducedNorm{\energyDifference(\timeSymbol)} +
    \int_0^\timeSymbol \euclideanNorm{\thermodynamicCholeskyMat^{-1}\energyResidual{2}(\dummyTimeSymbol)} d\dummyTimeSymbol \\
    \euclideanNorm{\continuousPositionError(\timeSymbol)} & \leq
    \euclideanNorm{\positionInitialError} + 
    \euclideanNorm{\positionDifference(\timeSymbol)} + 
    \int_0^\timeSymbol \euclideanNorm{\positionResidual{2}(\dummyTimeSymbol)} d\dummyTimeSymbol.
  \end{aligned}
\end{equation}
Next, we estimate the last terms on the right hand side of
\eqref{eq:error-eq-int3a}.  We rewrite the residuals as 
\begin{equation}\label{eq:res2a}
\begin{aligned}
\euclideanNorm{\kinematicCholeskyMat^{-1} \velocityResidual{2}} & = 
\euclideanNorm{(\kinematicCholeskyMat^T \velocityBasis \ROMKinematicMassMat^{-1}\velocityBasis^T \kinematicCholeskyMat) \kinematicCholeskyMat^{-1} (\forceOne(\velocity, \energy, \position)-
\forceOneObliqProjMat\forceOne(\velocityApprox, \energyApprox, \positionApprox))} \\
\euclideanNorm{\thermodynamicCholeskyMat^{-1}\energyResidual{2}} & = 
\euclideanNorm{(\thermodynamicCholeskyMat^T \energyBasis \ROMThermodynamicMassMat^{-1} \energyBasis^T \thermodynamicCholeskyMat) \thermodynamicCholeskyMat^{-1} (\forceTv(\velocity, \energy, \position) - 
\forceTvObliqProjMat\forceTv(\velocityApprox, \energyApprox, \positionApprox))} \\
\euclideanNorm{\positionResidual{2}} & = 
\euclideanNorm{(\positionBasis \positionBasis^T) \kinematicCholeskyMat^{-T} \kinematicCholeskyMat^{T}  
(\velocity - \velocityApprox)}.
\end{aligned}
\end{equation}
Note that the matrices 
\begin{equation}
\begin{split}
\ProjectionSymbol_{\velocitySymbol} & = 
\kinematicCholeskyMat^T \velocityBasis \ROMKinematicMassMat^{-1}\velocityBasis^T \kinematicCholeskyMat \\
\ProjectionSymbol_{\energySymbol} & = 
\thermodynamicCholeskyMat^T \energyBasis \ROMThermodynamicMassMat^{-1} \energyBasis^T \thermodynamicCholeskyMat \\
\ProjectionSymbol_{\positionSymbol} & = 
\positionBasis \positionBasis^T
\end{split}
\end{equation}
are orthogonal projection matrices and hence have norm $1$.  Therefore, we have
the following estimates for the residuals 
\begin{equation}\label{eq:res3a}
\begin{aligned}
\euclideanNorm{\kinematicCholeskyMat^{-1} \velocityResidual{2}} & \leq
\euclideanNorm{\kinematicCholeskyMat^{-1} (\forceOne(\velocity, \energy, \position)-
\forceOneObliqProjMat\forceOne(\velocityApprox, \energyApprox, \positionApprox))} \\
\euclideanNorm{\thermodynamicCholeskyMat^{-1}\energyResidual{2}} & \leq
\euclideanNorm{\thermodynamicCholeskyMat^{-1} (\forceTv(\velocity, \energy, \position) - 
\forceTvObliqProjMat\forceTv(\velocityApprox, \energyApprox, \positionApprox))} \\
\euclideanNorm{\positionResidual{2}} & \leq 
\euclideanNorm{\kinematicCholeskyMat^{-T} \kinematicCholeskyMat^{T}  
(\velocity - \velocityApprox)}.
\end{aligned}
\end{equation}
We further invoke the fact that the norm of a symmetric and positive definite
matrix is the square of that of its Cholesky factor and its square root, to
obtain 
\begin{equation}\label{eq:res4a}
\begin{aligned}
\euclideanNorm{\kinematicCholeskyMat^{-1} \velocityResidual{2}} & \leq
\euclideanNorm{ \kinematicMassMat^{-1/2}}  \euclideanNorm{\forceOne(\velocity, \energy, \position)-
\forceOneObliqProjMat\forceOne(\velocityApprox, \energyApprox, \positionApprox)} \\
\euclideanNorm{\thermodynamicCholeskyMat^{-1}\energyResidual{2}} & \leq
\euclideanNorm{ \thermodynamicMassMat^{-1/2}} \euclideanNorm{\forceTv(\velocity, \energy, \position) - 
\forceTvObliqProjMat\forceTv(\velocityApprox, \energyApprox, \positionApprox)} \\
\euclideanNorm{\positionResidual{2}} & \leq 
\euclideanNorm{ \kinematicMassMat^{-1/2}} \euclideanNorm{\kinematicCholeskyMat^{T}  
(\velocity - \velocityApprox)}.
\end{aligned}
\end{equation}
Using triangular inequalities, we get
\begin{equation}\label{eq:res5a}
  \begin{aligned}
\euclideanNorm{\kinematicCholeskyMat^{-1} \velocityResidual{2}} & \leq 
\euclideanNorm{ \kinematicMassMat^{-1/2}} \left(
\euclideanNorm{(\velocityIdentity - \forceOneObliqProjMat)\forceOne(\velocity, \energy, \position)} + \euclideanNorm{ \forceOneObliqProjMat (\forceOne(\velocity, \energy, \position)-
\forceOne(\velocityApprox, \energyApprox, \positionApprox))} \right) \\
\euclideanNorm{\thermodynamicCholeskyMat^{-1}\energyResidual{2}} & \leq
\euclideanNorm{ \thermodynamicMassMat^{-1/2}} \left(
\euclideanNorm{ (\energyIdentity - \forceTvObliqProjMat) \forceTv(\velocity, \energy, \position)} + \euclideanNorm{ \forceTvObliqProjMat ( \forceTv(\velocity, \energy, \position) - 
\forceTv(\velocityApprox, \energyApprox, \positionApprox))} \right) \\
\euclideanNorm{\positionResidual{2}} & \leq \euclideanNorm{ \kinematicMassMat^{-1/2}} \velocityInducedNorm
{\velocity - \velocityApprox}.
  \end{aligned}
\end{equation}
With the assumptions \eqref{eq:LipRHS-1}, we substitute \eqref{eq:res5a} into
\eqref{eq:error-eq-int3a} and obtain 
\begin{equation}\label{eq:error-eq-int4a}
  \begin{aligned}
    \velocityInducedNorm{\continuousVelocityError(\timeSymbol)} & \leq
    \velocityInducedNorm{\velocityInitialError} + 
    \velocityInducedNorm{\velocityDifference(\timeSymbol)} +
    \euclideanNorm{ \kinematicMassMat^{-1/2}} \int_0^\timeSymbol 
\euclideanNorm{(\velocityIdentity - \forceOneObliqProjMat)\forceOne(\velocity, \energy, \position)} + 
\euclideanNorm{\forceOneObliqProjMat } K_1 \fullNorm{(\continuousVelocityError(\dummyTimeSymbol), \continuousEnergyError(\dummyTimeSymbol), \continuousPositionError(\dummyTimeSymbol))} d\dummyTimeSymbol \\
    \energyInducedNorm{\continuousEnergyError(\timeSymbol)} & \leq
    \energyInducedNorm{\energyInitialError} + 
    \energyInducedNorm{\energyDifference(\timeSymbol)} +
    \euclideanNorm{ \thermodynamicMassMat^{-1/2}} \int_0^\timeSymbol 
\euclideanNorm{((\energyIdentity - \forceTvObliqProjMat) \forceTv(\velocity, \energy, \position)} +
\euclideanNorm{\forceTvObliqProjMat } K_2 \fullNorm{(\continuousVelocityError(\dummyTimeSymbol), \continuousEnergyError(\dummyTimeSymbol), \continuousPositionError(\dummyTimeSymbol))} d\dummyTimeSymbol \\
    \euclideanNorm{\continuousPositionError(\timeSymbol)} & \leq
    \euclideanNorm{\positionInitialError} + 
    \euclideanNorm{\positionDifference(\timeSymbol)} + 
    \euclideanNorm{ \kinematicMassMat^{-1/2}} \int_0^\timeSymbol \fullNorm{(\continuousVelocityError(\dummyTimeSymbol), \continuousEnergyError(\dummyTimeSymbol), \continuousPositionError(\dummyTimeSymbol))} d\dummyTimeSymbol.
  \end{aligned}
\end{equation}
which further implies 
\begin{equation}\label{eq:error-eq-int5a}
  \begin{aligned}
    \fullNorm{(\continuousVelocityError(\timeSymbol), \continuousEnergyError(\timeSymbol), \continuousPositionError(\timeSymbol))} & \leq
    \sqrt{3} \Big[ \fullNorm{(\velocityInitialError, \energyInitialError, \positionInitialError)} + 
    \max_{0 \leq \dummyTimeSymbol \leq \timeSymbol} \fullNorm{(\velocityDifference(\dummyTimeSymbol), \energyDifference(\dummyTimeSymbol), \positionDifference(\dummyTimeSymbol)} + \\
   & \qquad \euclideanNorm{ \kinematicMassMat^{-1/2}} \int_0^\timeSymbol 
\euclideanNorm{(\velocityIdentity - \forceOneObliqProjMat)\forceOne(\velocity, \energy, \position)} d\dummyTimeSymbol + \\
   &\qquad \euclideanNorm{ \thermodynamicMassMat^{-1/2}} \int_0^\timeSymbol 
\euclideanNorm{((\energyIdentity - \forceTvObliqProjMat) \forceTv(\velocity, \energy, \position)} d\dummyTimeSymbol + \\ 
&\qquad C_1
\int_0^\timeSymbol \fullNorm{(\continuousVelocityError(\dummyTimeSymbol), \continuousEnergyError(\dummyTimeSymbol), \continuousPositionError(\dummyTimeSymbol))} d\dummyTimeSymbol \Big], 
  \end{aligned}
\end{equation}
where $C_1 = \max\left\{ 
\euclideanNorm{ \kinematicMassMat^{-1/2}} \euclideanNorm{\forceOneObliqProjMat } K_1, 
\euclideanNorm{ \thermodynamicMassMat^{-1/2}} \euclideanNorm{\forceTvObliqProjMat } K_2, 
\euclideanNorm{ \kinematicMassMat^{-1/2}}
\right\}$. 
By Gronwall's inequality, we conclude 
\begin{equation}\label{eq:error-eq-int6a}
  \begin{aligned}
    \fullNorm{(\continuousVelocityError(\timeSymbol), \continuousEnergyError(\timeSymbol), \continuousPositionError(\timeSymbol))} & \leq
    \sqrt{3} \exp\left(C_1  \timeSymbol \right) \Big[ \fullNorm{(\velocityInitialError, \energyInitialError, \positionInitialError)} + 
    \max_{0 \leq \dummyTimeSymbol \leq \timeSymbol} \fullNorm{(\velocityDifference(\dummyTimeSymbol), \energyDifference(\dummyTimeSymbol), \positionDifference(\dummyTimeSymbol)} + \\
   & \qquad \euclideanNorm{ \kinematicMassMat^{-1/2}} \int_0^\timeSymbol 
\euclideanNorm{(\velocityIdentity - \forceOneObliqProjMat)\forceOne(\velocity, \energy, \position)} d\dummyTimeSymbol + \\
   &\qquad \euclideanNorm{ \thermodynamicMassMat^{-1/2}} \int_0^\timeSymbol 
\euclideanNorm{((\energyIdentity - \forceTvObliqProjMat) \forceTv(\velocity, \energy, \position)} d\dummyTimeSymbol  \Big], 
  \end{aligned}
\end{equation}
which provides the desired result. 
\end{proof}
\end{theorem}

\subsection{A-posteriori estimate for approximation error}
Next, we present an a-posteriori error estimate between the continuous-in-time
full order model solution in \eqref{eq:fom} and the continuous-in-time reduced
order approximation in \eqref{eq:rom-hr}.  The error bound depends on the
reduced order approximation and is controlled by a combination of the mismatch
of initial condition  and the oblique projection error of the nonlinear terms in
the Euclidean norm.  
\begin{theorem}\label{thm:aposteriori-semi} 
  Assume the Lipchitz continuity conditions \eqref{eq:LipRHS-1} hold.  Then
  there exists a generic constant $C>0$ such that for $\timeSymbol > 0$, we have 
\begin{equation}\label{eq:aposteriori-semi}
  \begin{aligned}
    \fullNorm{(\continuousVelocityError(\timeSymbol), \continuousEnergyError(\timeSymbol), \continuousPositionError(\timeSymbol))} & \leq 
    C \exp\left(C  \timeSymbol \right) \Big[ \fullNorm{(\continuousVelocityError(0), \continuousEnergyError(0), \continuousPositionError(0))} + \\
   &\qquad \int_0^\timeSymbol 
\euclideanNorm{(\velocityIdentity - \kinematicMassMat \velocityBasis
    \ROMKinematicMassMat^{-1}\velocityBasis^T
    \forceOneObliqProjMat)\forceOne(\velocityApprox, \energyApprox,
    \positionApprox)} d\dummyTimeSymbol + \\
   &\qquad \int_0^\timeSymbol \euclideanNorm{(\energyIdentity -
    \thermodynamicMassMat \energyBasis \ROMThermodynamicMassMat^{-1}
    \energyBasis^T \forceTvObliqProjMat) \forceTv(\velocityApprox,
    \energyApprox, \positionApprox)} d\dummyTimeSymbol
    \Big]
  \end{aligned}
\end{equation}
where $\continuousVelocityError(\timeSymbol)$,
  $\continuousEnergyError(\timeSymbol)$, and
  $\continuousPositionError(\timeSymbol)$ are defined in
  Eq.~\eqref{eq:rom-cerror}.
\begin{proof}
Instead of decomposing the residuals as in \eqref{eq:res-def-a}, we define 
\begin{equation}\label{eq:res-def-b}
\begin{aligned}
\velocityResidual{1} & = 
-(\velocityIdentity - \kinematicMassMat \velocityBasis \ROMKinematicMassMat^{-1}\velocityBasis^T  \forceOneObliqProjMat)\forceOne(\velocityApprox, \energyApprox, \positionApprox) \\
\velocityResidual{2} & = 
-(\forceOne(\velocity, \energy, \position)-
\forceOne(\velocityApprox, \energyApprox, \positionApprox)) \\
\energyResidual{1} & = 
(\energyIdentity -  \thermodynamicMassMat \energyBasis \ROMThermodynamicMassMat^{-1} \energyBasis^T \forceTvObliqProjMat) \forceTv(\velocityApprox, \energyApprox, \positionApprox) \\
\energyResidual{2} & = 
\forceTv(\velocity, \energy, \position) - 
\forceTv(\velocityApprox, \energyApprox, \positionApprox) \\
\positionResidual{1} & = (\positionIdentity - \positionBasis \positionBasis^T) \velocityApprox \\
\positionResidual{2} & = \velocity - \velocityApprox.
\end{aligned}
\end{equation}
Then \eqref{eq:error-eq-int} still holds true.  By triangle inequality, we have 
\begin{equation}\label{eq:error-eq-int2b}
  \begin{aligned}
    \velocityInducedNorm{\continuousVelocityError(\timeSymbol)} & \leq
    \velocityInducedNorm{\continuousVelocityError(0)} + 
    \int_0^\timeSymbol \euclideanNorm{\kinematicCholeskyMat^{-1} \velocityResidual{1}(\dummyTimeSymbol)} d\dummyTimeSymbol +
    \int_0^\timeSymbol \euclideanNorm{\kinematicCholeskyMat^{-1} \velocityResidual{2}(\dummyTimeSymbol)} d\dummyTimeSymbol \\
    \energyInducedNorm{\continuousEnergyError(\timeSymbol)} & \leq
    \energyInducedNorm{\continuousEnergyError(0)} + 
    \int_0^\timeSymbol \euclideanNorm{\thermodynamicCholeskyMat^{-1}\energyResidual{1}(\dummyTimeSymbol)} d\dummyTimeSymbol +
    \int_0^\timeSymbol \euclideanNorm{\thermodynamicCholeskyMat^{-1}\energyResidual{2}(\dummyTimeSymbol)} d\dummyTimeSymbol \\
    \euclideanNorm{\continuousPositionError(\timeSymbol)} & \leq
     \euclideanNorm{\continuousPositionError(0)} + 
    \int_0^\timeSymbol \euclideanNorm{\positionResidual{1}(\dummyTimeSymbol)} d\dummyTimeSymbol + 
    \int_0^\timeSymbol \euclideanNorm{\positionResidual{2}(\dummyTimeSymbol)} d\dummyTimeSymbol.
  \end{aligned}
\end{equation}
Using the fact that the norm of a symmetric and positive definite matrix is the
  square of that of its Cholesky factor and its square root, and invoking the
  assumptions \eqref{eq:LipRHS-1}, we have 
\begin{equation}\label{eq:error-eq-int3b}
\begin{aligned}
    \velocityInducedNorm{\continuousVelocityError(\timeSymbol)} & \leq
    \velocityInducedNorm{\continuousVelocityError(0)} + \euclideanNorm{ \kinematicMassMat^{-1/2}}
    \left [ \int_0^\timeSymbol \euclideanNorm{(\velocityIdentity - \kinematicMassMat
    \velocityBasis \ROMKinematicMassMat^{-1}\velocityBasis^T
    \forceOneObliqProjMat)\forceOne(\velocityApprox, \energyApprox,
    \positionApprox)} + K_1
    \fullNorm{(\continuousVelocityError(\dummyTimeSymbol),
    \continuousEnergyError(\dummyTimeSymbol), \continuousPositionError(\dummyTimeSymbol))}
    d\dummyTimeSymbol \right ] \\
    \energyInducedNorm{\continuousEnergyError(\timeSymbol)} & \leq
    \energyInducedNorm{\continuousEnergyError(0)} + \euclideanNorm{ \thermodynamicMassMat^{-1/2}}
    \left [ \int_0^\timeSymbol \euclideanNorm{(\energyIdentity -  \thermodynamicMassMat
    \energyBasis \ROMThermodynamicMassMat^{-1} \energyBasis^T
    \forceTvObliqProjMat) \forceTv(\velocityApprox, \energyApprox,
    \positionApprox)} + K_2
    \fullNorm{(\continuousVelocityError(\dummyTimeSymbol),
    \continuousEnergyError(\dummyTimeSymbol), \continuousPositionError(\dummyTimeSymbol))}
    d\dummyTimeSymbol \right ] \\
    \euclideanNorm{\continuousPositionError(\timeSymbol)} & \leq
     \euclideanNorm{\continuousPositionError(0)} + 
    \int_0^\timeSymbol \euclideanNorm{(\positionIdentity - \positionBasis \positionBasis^T) \velocityApprox} + \euclideanNorm{ \kinematicMassMat^{-1/2}} \fullNorm{(\continuousVelocityError(\dummyTimeSymbol), \continuousEnergyError(\dummyTimeSymbol), \continuousPositionError(\dummyTimeSymbol))} d\dummyTimeSymbol.
\end{aligned}
\end{equation}
This implies 
\begin{equation}\label{eq:error-eq-int4b}
  \begin{aligned}
    \fullNorm{(\continuousVelocityError(\timeSymbol), \continuousEnergyError(\timeSymbol), \continuousPositionError(\timeSymbol))} & \leq
    \sqrt{3} \Big[ \fullNorm{(\continuousVelocityError(0), \continuousEnergyError(0), \continuousPositionError(0))} + \\
   &\qquad \euclideanNorm{ \kinematicMassMat^{-1/2}} \int_0^\timeSymbol 
\euclideanNorm{(\velocityIdentity - \kinematicMassMat \velocityBasis
    \ROMKinematicMassMat^{-1}\velocityBasis^T
    \forceOneObliqProjMat)\forceOne(\velocityApprox, \energyApprox,
    \positionApprox)} d\dummyTimeSymbol + \\
   &\qquad \euclideanNorm{ \thermodynamicMassMat^{-1/2}} \int_0^\timeSymbol 
\euclideanNorm{(\energyIdentity -  \thermodynamicMassMat \energyBasis
    \ROMThermodynamicMassMat^{-1} \energyBasis^T \forceTvObliqProjMat)
    \forceTv(\velocityApprox, \energyApprox, \positionApprox)}
    d\dummyTimeSymbol + \\ 
&\qquad C_2
\int_0^\timeSymbol \fullNorm{(\continuousVelocityError(\dummyTimeSymbol), \continuousEnergyError(\dummyTimeSymbol), \continuousPositionError(\dummyTimeSymbol)} d\dummyTimeSymbol \Big], 
  \end{aligned}
\end{equation}
where $C_2 = \max\left\{ 
\euclideanNorm{ \kinematicMassMat^{-1/2}} K_1, 
\euclideanNorm{ \thermodynamicMassMat^{-1/2}} K_2, 
\euclideanNorm{ \kinematicMassMat^{-1/2}}
\right\}$. 
By Gronwall's inequality, we conclude 
\begin{equation}\label{eq:error-eq-int5b}
  \begin{aligned}
    \fullNorm{(\continuousVelocityError(\timeSymbol), \continuousEnergyError(\timeSymbol), \continuousPositionError(\timeSymbol))} & \leq
    \sqrt{3} \exp\left(C_2  \timeSymbol \right) \Big[ \fullNorm{(\continuousVelocityError(0), \continuousEnergyError(0), \continuousPositionError(0))} + \\
   &\qquad \euclideanNorm{ \kinematicMassMat^{-1/2}} \int_0^\timeSymbol
    \euclideanNorm{(\velocityIdentity - \kinematicMassMat \velocityBasis
    \ROMKinematicMassMat^{-1}\velocityBasis^T
    \forceOneObliqProjMat)\forceOne(\velocityApprox, \energyApprox,
    \positionApprox)} d\dummyTimeSymbol + \\
   &\qquad \euclideanNorm{ \thermodynamicMassMat^{-1/2}} \int_0^\timeSymbol
    \euclideanNorm{(\energyIdentity -  \thermodynamicMassMat \energyBasis
    \ROMThermodynamicMassMat^{-1} \energyBasis^T \forceTvObliqProjMat)
    \forceTv(\velocityApprox, \energyApprox, \positionApprox)}
    d\dummyTimeSymbol  \Big], 
  \end{aligned}
\end{equation}
which provides the desired result. 
\end{proof}
\end{theorem}

\subsection{Estimate for truncation error}
Finally, we analyze the error between the continuous-in-time reduced order
approximation in \eqref{eq:rom-hr} and the RK2-average fully discrete reduced
order approximation in \eqref{eq:RK2-avg-rom-hr}.  The error bound is controlled
by a combination of the mismatch of initial condition  and the maximum time step
size.

\begin{theorem}\label{thm:fully-discrete-errbound}
Assume $\forceMat$ is of class $C^2$ on $\kinematicFE \times \thermodynamicFE
  \times \kinematicFE$ and the Lipchitz continuity conditions
  \eqref{eq:LipRHS-1} holds.  In addition, assume there holds the following
  Lipchitz continuity conditions for the Jacobian of the force matrix
  $\forceMat$: there exists $K_3, K_4 > 0$ such that for any
  $\state = (\velocity, \energy, \position), \state' = (\velocity', \energy', \position') \in
  \kinematicFE \times \thermodynamicFE \times \kinematicFE$, 
\begin{equation}\label{eq:LipRHS-2}
  \begin{aligned}
\euclideanNorm{\jacobianState \forceOne(\state) \forceSysApprox(\state) -
    \jacobianState \forceOne(\state') \forceSysApprox(\state')} & \leq K_3
    \fullNorm{(\velocity-\velocity',\energy-\energy',\position-\position')} \\
\euclideanNorm{\jacobianState \forceTv(\state) \forceSysApprox(\state) -
    \jacobianState \forceTv(\state') \forceSysApprox(\state')} & \leq K_4
    \fullNorm{(\velocity-\velocity',\energy-\energy',\position-\position')}, 
  \end{aligned}
\end{equation}
where $\jacobianState \forceOne$ and $\jacobianState \forceTv$ are the Jacobian
  matrix of the vector-valued functions $\forceOne(\state)$ and
  $\forceTv(\state)$ respectively, and $\forceSysApprox$ denotes the lifted
  Euler update 
\begin{align}\label{eq:combinedApproxForce}
  \forceSysApprox(\stateApprox) &
        \equiv
  \pmat{-\velocityBasis\ROMKinematicMassMat^{-1}\velocityBasis^T
  \forceOneObliqProjMat \forceOne(\stateApprox) \\
        \energyBasis\ROMThermodynamicMassMat^{-1}\energyBasis^T
        \forceTvObliqProjMat \forceTv(\stateApprox) \\
        \positionBasis\positionBasis^T \velocityApprox}.
\end{align}
Then there exists a constant $C>0$ such that for $\timeIndex \geq 0$, we have 
\begin{equation}\label{eq:derror}
\begin{aligned}
    \fullNorm{(\discreteVelocityErrort{\timeIndex}, \discreteEnergyErrort{\timeIndex}, \discretePositionErrort{\timeIndex})} \leq e^{C \timek{\timeIndex}} (\fullNorm{(\discreteVelocityErrort{0}, \discreteEnergyErrort{0}, \discretePositionErrort{0})} + C(\timestep)^2),
\end{aligned}
\end{equation}
where $\timestep = \max_{0 \leq \dummyTimeIndex \leq \ntimestepROM-1}
  \timestepk{\dummyTimeIndex}$ and $\discreteVelocityErrort{\timeIndex}$,
  $\discreteEnergyErrort{\timeIndex}$, and $\discretePositionErrort{\timeIndex}$
  are defined in Eq.~\eqref{eq:rom-derror}.
\begin{proof}
First, we note that 
\begin{equation}\label{eq:derror-1}
\begin{aligned}
    \discreteVelocityErrort{\timeIndex+1} & = \discreteVelocityErrort{\timeIndex} + 
    \velocityBasis(\ROMvelocity(\timek{\timeIndex+1}) - \ROMvelocity(\timek{\timeIndex})) - \velocityBasis(\ROMvelocityt{\timeIndex+1} - \ROMvelocityt{\timeIndex}) \\
    \discreteEnergyErrort{\timeIndex+1} & = \discreteEnergyErrort{\timeIndex} + 
    \energyBasis(\ROMenergy(\timek{\timeIndex+1}) - \ROMenergy(\timek{\timeIndex})) - \energyBasis(\ROMenergyt{\timeIndex+1} - \ROMenergyt{\timeIndex}) \\ 
    \discretePositionErrort{\timeIndex+1} & = \discretePositionErrort{\timeIndex} + 
    \positionBasis(\ROMposition(\timek{\timeIndex+1}) - \ROMposition(\timek{\timeIndex})) - \positionBasis(\ROMpositiont{\timeIndex+1} - \ROMpositiont{\timeIndex}). 
\end{aligned}
\end{equation}
We will estimate the second term and the last term 
on the right hand side of each equation in \eqref{eq:derror-1}. 
By Taylor's remainder theorem, we have 
\begin{equation}\label{eq:taylor-sol1}
\begin{aligned}
\ROMvelocity(\timek{\timeIndex+1}) & = \ROMvelocity(\timek{\timeIndex}) + 
\timestepk{\timeIndex} \frac{d\ROMvelocity}{d\timeSymbol}(\timek{\timeIndex}) + 
\dfrac{(\timestepk{\timeIndex})^2}{2} \frac{d^2\ROMvelocity}{d\timeSymbol^2}(\timek{\timeIndex}) + (\timestepk{\timeIndex})^3 \velocityResidual{1} \\
\ROMenergy(\timek{\timeIndex+1}) & = \ROMenergy(\timek{\timeIndex}) + 
\timestepk{\timeIndex} \frac{d\ROMenergy}{d\timeSymbol}(\timek{\timeIndex}) + 
\dfrac{(\timestepk{\timeIndex})^2}{2} \frac{d^2\ROMenergy}{d\timeSymbol^2}(\timek{\timeIndex}) + (\timestepk{\timeIndex})^3 \energyResidual{1} \\
\ROMposition(\timek{\timeIndex+1}) & = \ROMposition(\timek{\timeIndex}) + 
\timestepk{\timeIndex} \frac{d\ROMposition}{d\timeSymbol}(\timek{\timeIndex}) + 
\dfrac{(\timestepk{\timeIndex})^2}{2} \frac{d^2\ROMposition}{d\timeSymbol^2}(\timek{\timeIndex}) + (\timestepk{\timeIndex})^3 \positionResidual{1},
\end{aligned}
\end{equation}
where the remainder is given by 
\begin{equation}\label{eq:res-def-c}
\begin{aligned}
\velocityResidual{1} & = 
\frac{1}{2(\timestepk{\timeIndex})^3} \int_{\timek{\timeIndex}}^{\timek{\timeIndex+1}} 
\frac{d^3\ROMvelocity}{d\timeSymbol^3}(\dummyTimeSymbol) (\dummyTimeSymbol - \timek{\timeIndex})^2 d\dummyTimeSymbol \\
\energyResidual{1} & = 
\frac{1}{2(\timestepk{\timeIndex})^3} \int_{\timek{\timeIndex}}^{\timek{\timeIndex+1}} 
\frac{d^3\ROMenergy}{d\timeSymbol^3}(\dummyTimeSymbol) (\dummyTimeSymbol - \timek{\timeIndex})^2 d\dummyTimeSymbol \\
\positionResidual{1} & = 
\frac{1}{2(\timestepk{\timeIndex})^3} \int_{\timek{\timeIndex}}^{\timek{\timeIndex+1}} 
\frac{d^3\ROMposition}{d\timeSymbol^3}(\dummyTimeSymbol) (\dummyTimeSymbol - \timek{\timeIndex})^2 d\dummyTimeSymbol. 
\end{aligned}
\end{equation}
We rewrite \eqref{eq:rom-hr} as 
\begin{align}\label{eq:rom-hr-rewrite}
   \frac{d\ROMvelocity}{d\timeSymbol} &=
    -\ROMKinematicMassMat^{-1} \velocityBasis^T \forceOneObliqProjMat \forceOne(\stateApprox) \\
    \frac{d\ROMenergy}{d\timeSymbol} &=
    \ROMThermodynamicMassMat^{-1} \energyBasis^T \forceTvObliqProjMat \forceTv(\stateApprox) \\
    \frac{d\ROMposition}{d\timeSymbol} &= \positionBasis^T \velocityApprox,
\end{align}
and differentiate with respect to time to obtain 
\begin{align}\label{eq:rom-hr-2nd-diff}
   \frac{d^2\ROMvelocity}{d\timeSymbol^2} &=
    -\ROMKinematicMassMat^{-1} \velocityBasis^T \forceOneObliqProjMat 
    \jacobianState \forceOne(\stateApprox) \frac{d\stateApprox}{d\timeSymbol} \\
    \frac{d^2\ROMenergy}{d\timeSymbol^2} &=
    \ROMThermodynamicMassMat^{-1} \energyBasis^T \forceTvObliqProjMat 
    \jacobianState \forceTv(\stateApprox) \frac{d\stateApprox}{d\timeSymbol} \\
    \frac{d^2\ROMposition}{d\timeSymbol^2} &= \positionBasis^T \frac{d\velocityApprox}{d\timeSymbol}. 
\end{align}
Substituting back to \eqref{eq:taylor-sol1}, we observe that 
\begin{equation}\label{eq:taylor-sol2}
\begin{aligned}
\ROMvelocity(\timek{\timeIndex+1}) & = \ROMvelocity(\timek{\timeIndex}) - 
\timestepk{\timeIndex} \ROMKinematicMassMat^{-1} \velocityBasis^T \forceOneObliqProjMat 
\left(\forceOne(\stateApprox(\timek{\timeIndex})) + \dfrac{\timestepk{\timeIndex}}{2}   
    \jacobianState \forceOne(\stateApprox(\timek{\timeIndex})) \forceSysApprox(\stateApprox(\timek{\timeIndex})) \right)  + (\timestepk{\timeIndex})^3 \velocityResidual{1} \\
\ROMenergy(\timek{\timeIndex+1}) & = \ROMenergy(\timek{\timeIndex}) + 
\timestepk{\timeIndex}     \ROMThermodynamicMassMat^{-1} \energyBasis^T \forceTvObliqProjMat 
\left(\forceTv(\stateApprox(\timek{\timeIndex})) + \dfrac{\timestepk{\timeIndex}}{2}   
    \jacobianState \forceTv(\stateApprox(\timek{\timeIndex})) \forceSysApprox(\stateApprox(\timek{\timeIndex})) \right) 
 + (\timestepk{\timeIndex})^3 \energyResidual{1} \\
\ROMposition(\timek{\timeIndex+1}) & = \ROMposition(\timek{\timeIndex}) + 
\timestepk{\timeIndex} \positionBasis^{T} 
\left(\velocityApprox(\timek{\timeIndex}) - \dfrac{\timestepk{\timeIndex}}{2}   \velocityBasis\ROMKinematicMassMat^{-1}\velocityBasis^T \forceOneObliqProjMat \forceOne(\stateApprox(\timek{\timeIndex}))\right) 
  + (\timestepk{\timeIndex})^3 \positionResidual{1},
\end{aligned}
\end{equation}
where $\frac{d\stateApprox}{d\timeSymbol}(\timek{\timeIndex}) =
\forceSysApprox(\stateApprox(\timek{\timeIndex}))$ is used through
Eq.~\eqref{eq:combinedApproxForce}.
Next, we are going to find a similar relation for the time-discrete
reduced-order coefficients.  First, we rewrite the first Runge-Kutta stage of
\eqref{eq:RK2-avg-rom-hr} to obtain 
\begin{align}\label{eq:RK2-avg-rom-hr-rewrite-1}
  \ROMvelocityt{\timeIndex+\frac{1}{2}} &= \ROMvelocityt{\timeIndex} -
  \frac{\timestepk{\timeIndex}}{2} \ROMKinematicMassMat^{-1} \velocityBasis^T
  \forceOneObliqProjMat \forceOne (\stateApproxt{\timeIndex}) \\
  \ROMenergyt{\timeIndex+\frac{1}{2}} &= \ROMenergyt{\timeIndex} +
  \frac{\timestepk{\timeIndex}}{2} \ROMThermodynamicMassMat^{-1}  \energyBasis^T
  \forceTvObliqProjMat \forceTv(\stateApproxt{\timeIndex}) +
  (\timestepk{\timeIndex})^2 \energyResidual{2} \\
  \ROMpositiont{\timeIndex+\frac{1}{2}} &= \ROMpositiont{\timeIndex} +
  \frac{\timestepk{\timeIndex}}{2} \positionBasis^T \velocityApproxt{\timeIndex}
  + (\timestepk{\timeIndex})^2 \positionResidual{2},
\end{align}
where $\forceOne$ and $\forceTv$ are defined in Eq.~\eqref{eq:forceOneforceTv}. 
The first Runge-Kutta stage correctors are given by 
\begin{equation}\label{eq:res-def-d}
\begin{aligned}
\energyResidual{2} & = - \dfrac{1}{4} \ROMThermodynamicMassMat^{-1}
  \energyBasis^T \forceTvObliqProjMat \left(\forceMat
  (\stateApproxt{\timeIndex}) \right )^T \cdot \velocityBasis
  \ROMKinematicMassMat^{-1} \velocityBasis^T \forceOneObliqProjMat\forceOne
  (\stateApproxt{\timeIndex}) \\
\positionResidual{2} & = -\dfrac{1}{4} \positionBasis^T  \velocityBasis
  \ROMKinematicMassMat^{-1} \velocityBasis^T \forceOneObliqProjMat \forceOne
  (\stateApproxt{\timeIndex}). 
\end{aligned}
\end{equation}
By denoting $\stateResidual{2} = [\zeroVec{\sizeKinematicFE};
\energyBasis\energyResidual{2}; \positionBasis\positionResidual{2}]^T$ and multiplying the basis matrices
to \eqref{eq:RK2-avg-rom-hr-rewrite-1}, we obtain
\begin{equation}\label{eq:RK2-avg-rom-hr-rewrite-1b}
 \stateApproxt{\timeIndex+\frac{1}{2}} = \stateApproxt{\timeIndex} +
  \frac{\timestepk{\timeIndex}}{2} \forceSysApprox(\stateApproxt{\timeIndex}) +
  (\timestepk{\timeIndex})^2 \stateResidual{2},
\end{equation} 
where $ \forceSysApprox(\stateApproxt{\timeIndex})$ is defined in
\eqref{eq:combinedApproxForce}.
Multiplying the basis matrix $\velocityBasis$ to
$\ROMvelocityt{\timeIndex+\frac{1}{2}}$ in \eqref{eq:RK2-avg-rom-hr-rewrite-1}
and $\ROMvelocityt{\timeIndex+1}$ in \eqref{eq:RK2-avg-rom-hr}, we have
\begin{equation}\label{eq:consistency-velocity-1}
\begin{aligned}
\velocityApproxt{\timeIndex+\frac{1}{2}} & =
  \velocityApproxt{\timeIndex} - \dfrac{\timestepk{\timeIndex}}{2}
  \velocityBasis \ROMKinematicMassMat^{-1} \velocityBasis^T
  \forceOneObliqProjMat \forceOne(\stateApproxt{\timeIndex}) \\
  \velocityApproxt{\timeIndex+1} & =
  \velocityApproxt{\timeIndex} - \timestepk{\timeIndex}
  \velocityBasis \ROMKinematicMassMat^{-1} \velocityBasis^T
  \forceOneObliqProjMat \forceOne(\stateApproxt{\timeIndex+\frac{1}{2}}),
\end{aligned}
\end{equation}
By the definition of $\avgvelocityApproxt{\timeIndex+\frac{1}{2}}$ below \eqref{eq:RK2-avg-rom-hr} 
and the second equation in \eqref{eq:consistency-velocity-1}, 
we observe that 
\begin{equation}\label{eq:consistency-velocity-2}
\begin{aligned}
\avgvelocityApproxt{\timeIndex+\frac{1}{2}} & =
  \velocityApproxt{\timeIndex} - \dfrac{\timestepk{\timeIndex}}{2}
  \velocityBasis \ROMKinematicMassMat^{-1} \velocityBasis^T
  \forceOneObliqProjMat \forceOne(\stateApproxt{\timeIndex+\frac{1}{2}}). 
\end{aligned}
\end{equation}
By subtracting the first equation in \eqref{eq:consistency-velocity-1} from \eqref{eq:consistency-velocity-2}, 
we obtain
\begin{equation}\label{eq:consistency-velocity}
\begin{aligned}
\avgvelocityApproxt{\timeIndex+\frac{1}{2}} & =
  \velocityApproxt{\timeIndex+\frac{1}{2}} + \dfrac{\timestepk{\timeIndex}}{2}
  \velocityBasis \ROMKinematicMassMat^{-1} \velocityBasis^T
  \forceOneObliqProjMat \left (\forceOne (\stateApproxt{\timeIndex}) - \forceOne
  (\stateApproxt{\timeIndex+\frac{1}{2}}) \right ), 
\end{aligned}
\end{equation}
which allows us to rewrite the second Runge-Kutta stage of
\eqref{eq:RK2-avg-rom-hr} as
\begin{align}\label{eq:RK2-avg-rom-hr-rewrite-2}
    \ROMvelocityt{\timeIndex+1} &= \ROMvelocityt{\timeIndex} -
    \timestepk{\timeIndex} \ROMKinematicMassMat^{-1} \velocityBasis^T
    \forceOneObliqProjMat \forceOne (\stateApproxt{\timeIndex+\frac{1}{2}}) \\
    \ROMenergyt{\timeIndex+1} &= \ROMenergyt{\timeIndex} +
    \timestepk{\timeIndex} \ROMThermodynamicMassMat^{-1}  \energyBasis^T
    \forceTvObliqProjMat  \forceTv (\stateApproxt{\timeIndex+\frac{1}{2}}) +
    (\timestepk{\timeIndex})^3 \energyResidual{3} \\
    \ROMpositiont{\timeIndex+1} &=
    \ROMpositiont{\timeIndex} + \timestepk{\timeIndex} \positionBasis^T
    \velocityApproxt{\timeIndex+\frac{1}{2}} + (\timestepk{\timeIndex})^3
    \positionResidual{3},
\end{align}
where the second Runge-Kutta correctors are given by 
\begin{equation}\label{eq:res-def-e}
\begin{aligned}
\energyResidual{3} & = \dfrac{1}{2\timestepk{\timeIndex}}
  \ROMThermodynamicMassMat^{-1}  \energyBasis^T \forceTvObliqProjMat
  \left(\forceMat (\stateApproxt{\timeIndex+\frac{1}{2}}) \right )^T \cdot
\velocityBasis \ROMKinematicMassMat^{-1} \velocityBasis^T \forceOneObliqProjMat
  \left(\forceOne (\stateApproxt{\timeIndex}) - \forceOne
  (\stateApproxt{\timeIndex+\frac{1}{2}}) \right ) \\
\positionResidual{3} & = \dfrac{1}{2\timestepk{\timeIndex}} \positionBasis^T
  \velocityBasis \ROMKinematicMassMat^{-1} \velocityBasis^T
  \forceOneObliqProjMat \left(\forceOne (\stateApproxt{\timeIndex}) - \forceOne
  (\stateApproxt{\timeIndex+\frac{1}{2}}) \right ). 
\end{aligned}
\end{equation}
Using Taylor's remainder theorem again, we have 
\begin{equation}\label{eq:taylor-RHS-1}
\begin{aligned}
\forceOne (\stateApproxt{\timeIndex+\frac{1}{2}}) & = 
\forceOne (\stateApproxt{\timeIndex}) + 
  \jacobianState \forceOne(\stateApproxt{\timeIndex}) \left(\stateApproxt{\timeIndex+\frac{1}{2}} - \stateApproxt{\timeIndex} \right) + (\timestepk{\timeIndex})^2 \forceOneResidual{4} \\
  \forceTv (\stateApproxt{\timeIndex+\frac{1}{2}}) & = 
\forceTv (\stateApproxt{\timeIndex}) + 
  \jacobianState \forceTv(\stateApproxt{\timeIndex}) \left(\stateApproxt{\timeIndex+\frac{1}{2}} - \stateApproxt{\timeIndex} \right) + (\timestepk{\timeIndex})^2 \forceTvResidual{4},
\end{aligned}
\end{equation}
where the remainders are given by 
\begin{equation}\label{eq:res-def-f}
\begin{aligned}
\forceOneResidual{4} & = \dfrac{2}{(\timestepk{\timeIndex})^2} \sum_{\vert \multiIndex \vert = 2} \dfrac{(\stateApproxt{\timeIndex+\frac{1}{2}} - \stateApproxt{\timeIndex})^{\multiIndex}}{\multiIndex !} 
  \int_0^1 (1-s) \partial^{\multiIndex} \forceOne\left((1-\dummyTimeSymbol)\stateApproxt{\timeIndex}+s \stateApproxt{\timeIndex+\frac{1}{2}}\right) d\dummyTimeSymbol \\
\forceTvResidual{4} & = \dfrac{2}{(\timestepk{\timeIndex})^2} \sum_{\vert \multiIndex \vert = 2} \dfrac{(\stateApproxt{\timeIndex+\frac{1}{2}} - \stateApproxt{\timeIndex})^{\multiIndex}}{\multiIndex !} 
  \int_0^1 (1-s) \partial^{\multiIndex} \forceTv\left((1-\dummyTimeSymbol)\stateApproxt{\timeIndex}+s \stateApproxt{\timeIndex+\frac{1}{2}}\right) d\dummyTimeSymbol.
\end{aligned}
\end{equation}
Here $\multiIndex = (\multiIndexSymbol_1, \multiIndexSymbol_2, \ldots, \multiIndexSymbol_\sizeWholeFE)$ is a multi-index with order $\vert \multiIndex \vert = \multiIndexSymbol_1 + \multiIndexSymbol_2 + \ldots + \multiIndexSymbol_\sizeWholeFE$, 
and the notations $\state^{\multiIndex}$, $\multiIndex !$ and $\partial^{\multiIndex}$ are formally defined as 
\begin{equation}
\begin{aligned}
\state^{\multiIndex} & = \stateSymbol_1^{\multiIndexSymbol_1} \stateSymbol_2^{\multiIndexSymbol_2} \ldots  \stateSymbol_\sizeWholeFE^{\multiIndexSymbol_\sizeWholeFE} \\
\multiIndex! & = \multiIndexSymbol_1 ! \multiIndexSymbol_2 ! \ldots \multiIndexSymbol_\sizeWholeFE ! \\
\partial^{\multiIndex} & = \dfrac{\partial^{\vert \multiIndex \vert}}{\partial \stateSymbol_1^{\multiIndexSymbol_1}\partial \stateSymbol_2^{\multiIndexSymbol_2}\cdots\partial \stateSymbol_\sizeWholeFE^{\multiIndexSymbol_\sizeWholeFE}}.
\end{aligned}
\end{equation}
Substituting \eqref{eq:RK2-avg-rom-hr-rewrite-1b} and \eqref{eq:taylor-RHS-1} into \eqref{eq:RK2-avg-rom-hr-rewrite-2}, we have 
\begin{align}\label{eq:RK2-avg-rom-hr-rewrite-2b}
    \ROMvelocityt{\timeIndex+1} &= \ROMvelocityt{\timeIndex} - \timestepk{\timeIndex}
    \ROMKinematicMassMat^{-1} \velocityBasis^T \forceOneObliqProjMat \left[ \forceOne (\stateApproxt{\timeIndex}) + 
  \jacobianState \forceOne(\stateApproxt{\timeIndex}) \left( \frac{\timestepk{\timeIndex}}{2} \forceSysApprox(\stateApproxt{\timeIndex}) + (\timestepk{\timeIndex})^2 \stateResidual{2} \right) + (\timestepk{\timeIndex})^2 \forceOneResidual{4} \right] \\
    \ROMenergyt{\timeIndex+1} &= \ROMenergyt{\timeIndex} + \timestepk{\timeIndex}
    \ROMThermodynamicMassMat^{-1}  \energyBasis^T \forceTvObliqProjMat  \left[ \forceTv (\stateApproxt{\timeIndex}) + 
  \jacobianState \forceTv(\stateApproxt{\timeIndex}) \left( \frac{\timestepk{\timeIndex}}{2} \forceSysApprox(\stateApproxt{\timeIndex}) + (\timestepk{\timeIndex})^2 \stateResidual{2} \right) + (\timestepk{\timeIndex})^2 \forceTvResidual{4} \right] + (\timestepk{\timeIndex})^3 \energyResidual{3} \\
    \ROMpositiont{\timeIndex+1} &=
    \ROMpositiont{\timeIndex} + \timestepk{\timeIndex}
    \positionBasis^T \left(\velocityApproxt{\timeIndex} - \frac{\timestepk{\timeIndex}}{2}
    \velocityBasis \ROMKinematicMassMat^{-1} \velocityBasis^T \forceOneObliqProjMat 
\forceOne (\stateApproxt{\timeIndex})\right) + (\timestepk{\timeIndex})^3 \positionResidual{3}.
\end{align}
Combining \eqref{eq:taylor-sol2}, \eqref{eq:RK2-avg-rom-hr-rewrite-2b} and \eqref{eq:derror-1} yields the error identities 
\begin{equation}\label{eq:derror-2}
\begin{aligned}
    \discreteVelocityErrort{\timeIndex+1} & = \discreteVelocityErrort{\timeIndex} + 
    \timestepk{\timeIndex} \rkVelocityErrorCoeff{1} + (\timestepk{\timeIndex})^2 \rkVelocityErrorCoeff{2} + (\timestepk{\timeIndex})^3 \rkVelocityErrorCoeff{3} \\
    \discreteEnergyErrort{\timeIndex+1} & = \discreteEnergyErrort{\timeIndex} + 
    \timestepk{\timeIndex} \rkEnergyErrorCoeff{1} + (\timestepk{\timeIndex})^2 \rkEnergyErrorCoeff{2} + (\timestepk{\timeIndex})^3 \rkEnergyErrorCoeff{3} \\
    \discretePositionErrort{\timeIndex+1} & = \discretePositionErrort{\timeIndex} + 
    \timestepk{\timeIndex} \rkPositionErrorCoeff{1} + (\timestepk{\timeIndex})^2 \rkPositionErrorCoeff{2} + (\timestepk{\timeIndex})^3 \rkPositionErrorCoeff{3}, 
\end{aligned}
\end{equation}
where the coefficient vectors are defined as 
\begin{equation}
\begin{aligned}
\rkVelocityErrorCoeff{1} & = 
- \velocityBasis \ROMKinematicMassMat^{-1} \velocityBasis^T \forceOneObliqProjMat 
\left( \forceOne(\stateApprox(\timek{\timeIndex})) - \forceOne (\stateApproxt{\timeIndex}) \right)\\
\rkEnergyErrorCoeff{1} & = 
\energyBasis \ROMThermodynamicMassMat^{-1} \energyBasis^T \forceTvObliqProjMat 
\left( \forceTv(\stateApprox(\timek{\timeIndex})) - \forceTv (\stateApproxt{\timeIndex}) \right)\\
\rkPositionErrorCoeff{1} & = 
\positionBasis \positionBasis^T 
\left( \velocityApprox(\timek{\timeIndex}) - \velocityApproxt{\timeIndex} \right)\\
\rkVelocityErrorCoeff{2} & = 
- \frac{1}{2} \velocityBasis \ROMKinematicMassMat^{-1} \velocityBasis^T \forceOneObliqProjMat 
\left( \jacobianState \forceOne(\stateApprox(\timek{\timeIndex})) \forceSysApprox(\stateApprox(\timek{\timeIndex})) - \jacobianState \forceOne(\stateApproxt{\timeIndex}) \forceSysApprox(\stateApproxt{\timeIndex}) \right) \\
\rkEnergyErrorCoeff{2} & = 
\frac{1}{2} \energyBasis \ROMThermodynamicMassMat^{-1} \energyBasis^T \forceTvObliqProjMat 
\left( \jacobianState \forceTv(\stateApprox(\timek{\timeIndex})) \forceSysApprox(\stateApprox(\timek{\timeIndex})) - \jacobianState \forceTv(\stateApproxt{\timeIndex}) \forceSysApprox(\stateApproxt{\timeIndex}) \right) \\
\rkPositionErrorCoeff{2} & = 
-\frac{1}{2} \positionBasis \positionBasis^T \velocityBasis\ROMKinematicMassMat^{-1} \velocityBasis^T \forceOneObliqProjMat \left( \forceOne(\stateApprox(\timek{\timeIndex})) - \forceOne (\stateApproxt{\timeIndex}) \right)\\
\rkVelocityErrorCoeff{3} & = \velocityBasis(\velocityResidual{1} + \ROMKinematicMassMat^{-1} \velocityBasis^T \forceOneObliqProjMat (\jacobianState \forceOne(\stateApproxt{\timeIndex}) \stateResidual{2} + \forceOneResidual{4})) \\
\rkEnergyErrorCoeff{3} & = \energyBasis( \energyResidual{1} - \energyResidual{3} - \ROMThermodynamicMassMat^{-1} \energyBasis^T \forceTvObliqProjMat (\jacobianState \forceTv(\stateApproxt{\timeIndex}) \stateResidual{2}  + \forceTvResidual{4})) \\
\rkPositionErrorCoeff{3} & = \positionBasis (\positionResidual{1} -  \positionResidual{3}).
\end{aligned}
\end{equation}
We remark that all the residual vectors only appear in $(\rkVelocityErrorCoeff{3}, \rkEnergyErrorCoeff{3}, 
\rkPositionErrorCoeff{3})$ and therefore scale with $(\Delta t_n)^3$.
We estimate the first order terms using a similar argument as \eqref{eq:res5a} in Theorem~\ref{thm:apriori-semi}, i.e. 
\begin{equation}
  \begin{aligned}
\euclideanNorm{\kinematicCholeskyMat^{-1} \rkVelocityErrorCoeff{1}} & \leq 
\euclideanNorm{ \kinematicMassMat^{-1/2}}  \euclideanNorm{ \forceOneObliqProjMat  
(\forceOne(\stateApprox(\timek{\timeIndex})) - \forceOne (\stateApproxt{\timeIndex}))} \\
\euclideanNorm{\thermodynamicCholeskyMat^{-1} \rkEnergyErrorCoeff{1}} & \leq
\euclideanNorm{ \thermodynamicMassMat^{-1/2}} 
\euclideanNorm{ \forceTvObliqProjMat ( \forceTv(\stateApprox(\timek{\timeIndex})) - \forceTv (\stateApproxt{\timeIndex}))}  \\
\euclideanNorm{\rkPositionErrorCoeff{1}} & \leq \euclideanNorm{ \kinematicMassMat^{-1/2}} \velocityInducedNorm
{\velocityApprox(\timek{\timeIndex}) - \velocityApproxt{\timeIndex}}.
  \end{aligned}
\end{equation}
With the assumption \eqref{eq:LipRHS-1}, we have 
\begin{equation}
\fullNorm{(\rkVelocityErrorCoeff{1}, \rkEnergyErrorCoeff{1}, \rkPositionErrorCoeff{1})}
\leq C_1 \fullNorm{(\discreteVelocityErrort{\timeIndex}, \discreteEnergyErrort{\timeIndex}, \discretePositionErrort{\timeIndex})}.
\end{equation}
where $C_1$ is the constant defined in \eqref{eq:error-eq-int5a}. 
Similarly, with the assumptions \eqref{eq:LipRHS-2} and \eqref{eq:LipRHS-1}, we 
estimate the second order terms by 
\begin{equation}
\fullNorm{(\rkVelocityErrorCoeff{2}, \rkEnergyErrorCoeff{2}, \rkPositionErrorCoeff{2})}
\leq C_3 \fullNorm{(\discreteVelocityErrort{\timeIndex}, \discreteEnergyErrort{\timeIndex}, \discretePositionErrort{\timeIndex})},
\end{equation}
with $C_3 = \frac{1}{2} \max\left\{ 
\euclideanNorm{ \kinematicMassMat^{-1/2}} \euclideanNorm{\forceOneObliqProjMat } K_3, 
\euclideanNorm{ \thermodynamicMassMat^{-1/2}} \euclideanNorm{\forceTvObliqProjMat } K_4, 
\euclideanNorm{\forceOneObliqProjMat}K_1
\right\}$. 
Using the error identities \eqref{eq:derror-2} and applying triangle inequality,
we have 
\begin{equation}\label{eq:derror-3}
\begin{aligned}
    \fullNorm{(\discreteVelocityErrort{\timeIndex+1}, \discreteEnergyErrort{\timeIndex+1}, \discretePositionErrort{\timeIndex+1})} \leq (1 + C_1 \timestepk{\timeIndex} + C_3 (\timestepk{\timeIndex})^2)
    \fullNorm{(\discreteVelocityErrort{\timeIndex}, \discreteEnergyErrort{\timeIndex}, \discretePositionErrort{\timeIndex})} + 
    (\timestepk{\timeIndex})^3 \fullNorm{(\rkVelocityErrorCoeff{3}, \rkEnergyErrorCoeff{3}, \rkPositionErrorCoeff{3})}.
\end{aligned}
\end{equation}
Next, we find a bound for the third order terms. 
By direct computation, we observe that 
\begin{equation}\label{eq:rkCoeff-3}
\begin{aligned}
\fullNorm{(\rkVelocityErrorCoeff{3}, \rkEnergyErrorCoeff{3}, \rkPositionErrorCoeff{3})}
& \leq \fullNorm{(\velocityBasis\velocityResidual{1}, \energyBasis\energyResidual{1}, \positionBasis\positionResidual{1})} + \\
& \quad \quad C_4 C_5 \euclideanNorm{ \left(\jacobianState \forceOne(\stateApproxt{\timeIndex}) \stateResidual{2}, 
\jacobianState \forceTv(\stateApproxt{\timeIndex}) \stateResidual{2} , 
\zeroVec{\sizeKinematicFE}\right) } + \\
& \quad \quad  C_4 \euclideanNorm{(\zeroVec{\sizeKinematicFE}, \energyBasis\energyResidual{3}, \positionBasis\positionResidual{3})} + 
C_4 C_5 \euclideanNorm{ \left(\forceOneResidual{4}, \forceTvResidual{4}, 
\zeroVec{\sizeKinematicFE}\right)},
\end{aligned}
\end{equation}
where $C_4 = \max\left\{ \euclideanNorm{ \kinematicMassMat^{1/2}}, 
\euclideanNorm{ \thermodynamicMassMat^{1/2}}, 1 \right\}$ and 
$C_5 = \max\left\{ \euclideanNorm{\forceOneObliqProjMat },  \euclideanNorm{\forceTvObliqProjMat } \right\}$.
We will estimate each of terms on the right hand side of \eqref{eq:rkCoeff-3}. 
We start with the last term on the right hand side of \eqref{eq:rkCoeff-3}, where 
we invoke \eqref{eq:RK2-avg-rom-hr-rewrite-1b} to estimate \eqref{eq:res-def-f} by
\begin{equation}\label{eq:res-def-f2}
\begin{aligned}
\euclideanNorm{\forceOneResidual{4}} 
& \leq \dfrac{1}{8} \sizeWholeFE M^{(4)}_{\forceOne}
\euclideanNorm{\forceSysApprox(\stateApproxt{\timeIndex}) + 2\timestepk{\timeIndex} \stateResidual{2}}^2 \\
\euclideanNorm{\forceTvResidual{4}} 
& \leq \dfrac{1}{8} \sizeWholeFE M^{(4)}_{\forceTv}
\euclideanNorm{\forceSysApprox(\stateApproxt{\timeIndex}) + 2\timestepk{\timeIndex} \stateResidual{2}}^2, 
\end{aligned}
\end{equation}
where the constants $M^{(2)}_{\forceOne}$ and $M^{(2)}_{\forceTv}$ are defined as
\begin{equation}
\begin{split}
M^{(2)}_{\forceOne} & = \max_{0 \leq \dummyTimeIndex \leq \ntimestepROM-1} 
\max_{\vert \multiIndex \vert = 2} \max_{0 \leq s \leq 1} 
\euclideanNorm{ \partial^{\multiIndex} \forceOne\left((1-\dummyTimeSymbol)\stateApproxt{\dummyTimeIndex}+ 
\dummyTimeSymbol \stateApproxt{\dummyTimeIndex+\frac{1}{2}}\right)} \\
M^{(2)}_{\forceTv} & = \max_{0 \leq \dummyTimeIndex \leq \ntimestepROM-1} 
\max_{\vert \multiIndex \vert = 2} \max_{0 \leq s \leq 1} 
\euclideanNorm{ \partial^{\multiIndex} \forceTv\left((1-\dummyTimeSymbol)\stateApproxt{\dummyTimeIndex}+ 
\dummyTimeSymbol \stateApproxt{\dummyTimeIndex+\frac{1}{2}}\right)}. 
\end{split}
\end{equation}
Therefore, the last term on the right hand side of 
\eqref{eq:rkCoeff-3} can be estimated by 
\begin{equation}
\label{eq:res-f}
\euclideanNorm{ \left(\forceOneResidual{4}, \forceTvResidual{4}, 
\zeroVec{\sizeKinematicFE}\right)} 
\leq \dfrac{1}{8} \sizeWholeFE (M_{\forceOne}^{(2)} + M_{\forceTv}^{(2)})
\euclideanNorm{\forceSysApprox(\stateApproxt{\timeIndex}) + 2\timestepk{\timeIndex} \stateResidual{2}}^2. 
\end{equation}
For the third term on the right hand side of \eqref{eq:rkCoeff-3}, substituting
\eqref{eq:taylor-RHS-1} and \eqref{eq:RK2-avg-rom-hr-rewrite-1b} into
\eqref{eq:res-def-e} gives 
\begin{equation}\label{eq:res-def-e2}
\begin{aligned}
\energyResidual{3} & = -\dfrac{1}{4} \ROMThermodynamicMassMat^{-1}  \energyBasis^T \forceTvObliqProjMat \left(\forceMat (\stateApproxt{\timeIndex+\frac{1}{2}}) \right )^T \cdot
\velocityBasis \ROMKinematicMassMat^{-1} \velocityBasis^T \forceOneObliqProjMat \left(\jacobianState \forceOne(\stateApproxt{\timeIndex}) \forceSysApprox(\stateApproxt{\timeIndex}) + 2\timestepk{\timeIndex} (\jacobianState \forceOne(\stateApproxt{\timeIndex}) \stateResidual{2} + \forceOneResidual{4})\right) \\
\positionResidual{3} & = -\dfrac{1}{4} \positionBasis^T  \velocityBasis \ROMKinematicMassMat^{-1} \velocityBasis^T \forceOneObliqProjMat \left(\jacobianState \forceOne(\stateApproxt{\timeIndex}) \forceSysApprox(\stateApproxt{\timeIndex}) + 2\timestepk{\timeIndex} (\jacobianState \forceOne(\stateApproxt{\timeIndex}) \stateResidual{2} + \forceOneResidual{4})\right).
\end{aligned}
\end{equation}
\eqref{eq:res-def-e2} and \eqref{eq:res-def-f2} together imply
\begin{equation}\label{eq:res-e}
\begin{split}
\euclideanNorm{(\zeroVec{\sizeKinematicFE}, \energyBasis\energyResidual{3}, \positionBasis\positionResidual{3})} 
& \leq \dfrac{1}{4} \left(\euclideanNorm{\ROMThermodynamicMassMat^{-1} } 
\euclideanNorm{\forceTvObliqProjMat}\euclideanNorm{\forceMat (\stateApproxt{\timeIndex+\frac{1}{2}})} + 1\right) \euclideanNorm{\ROMKinematicMassMat^{-1} } \euclideanNorm{\forceOneObliqProjMat } \\
& \quad \quad \left[\euclideanNorm{\jacobianState \forceOne(\stateApproxt{\timeIndex})} 
\euclideanNorm{\forceSysApprox(\stateApproxt{\timeIndex}) + 2\timestepk{\timeIndex} \stateResidual{2}} + 
\dfrac{1}{4}  \sizeWholeFE M^{(2)}_{\forceOne} \timestepk{\timeIndex} 
\euclideanNorm{\forceSysApprox(\stateApproxt{\timeIndex}) + 2\timestepk{\timeIndex} \stateResidual{2}}^2  \right].
\end{split}
\end{equation}
For the second term on the right hand side of \eqref{eq:rkCoeff-3}, we have 
\begin{equation}\label{eq:res-d}
\begin{aligned}
\euclideanNorm{ \left(\jacobianState \forceOne(\stateApproxt{\timeIndex}) \stateResidual{2}, 
\jacobianState \forceTv(\stateApproxt{\timeIndex}) \stateResidual{2} , 
\zeroVec{\sizeKinematicFE}\right) } 
& \leq 
\left( \euclideanNorm{\jacobianState \forceOne(\stateApproxt{\timeIndex})}^2 + 
\euclideanNorm{\jacobianState \forceTv(\stateApproxt{\timeIndex})}^2 \right)^\frac{1}{2}
\euclideanNorm{ \stateResidual{2} }.
\end{aligned}
\end{equation}
Using \eqref{eq:res-f}, \eqref{eq:res-e} and \eqref{eq:res-d}, we can rewrite \eqref{eq:rkCoeff-3} as 
\begin{equation}\label{eq:rkCoeff-3-2}
\begin{aligned}
\fullNorm{(\rkVelocityErrorCoeff{3}, \rkEnergyErrorCoeff{3}, \rkPositionErrorCoeff{3})}
& \leq \fullNorm{(\velocityBasis\velocityResidual{1}, \energyBasis\energyResidual{1}, \positionBasis\positionResidual{1})} + 
C_4 C_5 \left( M^{(1)}_{\forceOne} + M^{(1)}_{\forceTv} \right)
\euclideanNorm{ \stateResidual{2} } + \\
& \quad \quad \dfrac{C_4}{4}  
\left(\euclideanNorm{\ROMThermodynamicMassMat^{-1} } \euclideanNorm{\forceTvObliqProjMat}\euclideanNorm{\forceMat (\stateApproxt{\timeIndex+\frac{1}{2}})} + 1\right) \euclideanNorm{\ROMKinematicMassMat^{-1} }  \euclideanNorm{\forceOneObliqProjMat }M^{(1)}_{\forceOne}
\euclideanNorm{\forceSysApprox(\stateApproxt{\timeIndex}) + 2\timestepk{\timeIndex} \stateResidual{2}}+ \\
& \quad \quad  \dfrac{C_4}{16} \left[
2C_5 \sizeWholeFE (M_{\forceOne}^{(2)} + M_{\forceTv}^{(2)}) + 
\left(\euclideanNorm{\ROMThermodynamicMassMat^{-1} } \euclideanNorm{\forceTvObliqProjMat}\euclideanNorm{\forceMat (\stateApproxt{\timeIndex+\frac{1}{2}})} + 1\right) \euclideanNorm{\ROMKinematicMassMat^{-1} } \euclideanNorm{\forceOneObliqProjMat }
\sizeWholeFE M_{\forceOne}^{(2)} \timestepk{\timeIndex} \right] \\
& \quad \quad
\euclideanNorm{\forceSysApprox(\stateApproxt{\timeIndex}) + 2\timestepk{\timeIndex} \stateResidual{2}}^2,
\end{aligned}
\end{equation}
where the constants $M^{(1)}_{\forceOne}$ and $M^{(1)}_{\forceTv}$ are defined as
\begin{equation}
\begin{split}
M^{(1)}_{\forceOne} & = \max_{0 \leq \dummyTimeIndex \leq \ntimestepROM-1} 
\euclideanNorm{\jacobianState \forceOne(\stateApproxt{\dummyTimeIndex})} \\
M^{(1)}_{\forceTv} & = \max_{0 \leq \dummyTimeIndex \leq \ntimestepROM-1} 
\euclideanNorm{\jacobianState \forceTv(\stateApproxt{\dummyTimeIndex})}. 
\end{split}
\end{equation}
By the definition of $\stateResidual{2}$ before \eqref{eq:RK2-avg-rom-hr-rewrite-1b}, we obtain 
\begin{equation}\label{eq:res-def-d2}
\begin{aligned}
\euclideanNorm{ \stateResidual{2} } & \leq 
\dfrac{1}{4} 
\left(\euclideanNorm{\ROMThermodynamicMassMat^{-1} }  \euclideanNorm{\forceTvObliqProjMat } \euclideanNorm{ \forceMat(\stateApproxt{\timeIndex})} + 1 \right) \euclideanNorm{\ROMKinematicMassMat^{-1} } 
\euclideanNorm{\forceOneObliqProjMat } \euclideanNorm{ \forceMat(\stateApproxt{\timeIndex}) } \sqrt{\sizeKinematicFE}, 
\end{aligned}
\end{equation}
which further implies 
\begin{equation}\label{eq:res-def-d3}
\begin{aligned}
\euclideanNorm{ \forceSysApprox(\stateApproxt{\timeIndex}) + 2\timestepk{\timeIndex} \stateResidual{2} } & \leq 
\euclideanNorm{ \forceSysApprox(\stateApproxt{\timeIndex})} +  \dfrac{1}{2} 
\left(\euclideanNorm{\ROMThermodynamicMassMat^{-1} }  \euclideanNorm{\forceTvObliqProjMat } \euclideanNorm{ \forceMat(\stateApproxt{\timeIndex})} + 1\right) \euclideanNorm{\ROMKinematicMassMat^{-1} } 
\euclideanNorm{\forceOneObliqProjMat } \euclideanNorm{ \forceMat(\stateApproxt{\timeIndex}) } \sqrt{\sizeKinematicFE} \timestepk{\timeIndex}. 
\end{aligned}
\end{equation}
There, \eqref{eq:rkCoeff-3-2} can be further rewritten as 
\begin{equation}\label{eq:rkCoeff-3-3}
\begin{aligned}
\fullNorm{(\rkVelocityErrorCoeff{3}, \rkEnergyErrorCoeff{3}, \rkPositionErrorCoeff{3})}
& \leq \fullNorm{(\velocityBasis\velocityResidual{1}, \energyBasis\energyResidual{1}, \positionBasis\positionResidual{1})} + \\
& \quad \quad \dfrac{C_4 C_5}{4} \left( M^{(1)}_{\forceOne} + M^{(1)}_{\forceTv} \right)
\left(\euclideanNorm{\ROMThermodynamicMassMat^{-1} }  \euclideanNorm{\forceTvObliqProjMat } 
M^{(0)}_{\forceMat} + 1 \right)
\euclideanNorm{\forceOneObliqProjMat } M^{(0)}_{\forceMat} \sqrt{\sizeKinematicFE}  + \\
& \quad \quad \dfrac{C_4}{8}  
\left(\euclideanNorm{\ROMThermodynamicMassMat^{-1} } \euclideanNorm{\forceTvObliqProjMat}
M^{(0)}_{\forceMat} + 1\right) \euclideanNorm{\ROMKinematicMassMat^{-1} }  \euclideanNorm{\forceOneObliqProjMat }M^{(1)}_{\forceOne} \\
& \quad \quad \left[2\tilde{M}^{(0)}_{\forceSysApprox} + 
\left(\euclideanNorm{\ROMThermodynamicMassMat^{-1} }  \euclideanNorm{\forceTvObliqProjMat } 
M^{(0)}_{\forceMat} + 1 \right) \euclideanNorm{\ROMKinematicMassMat^{-1} } 
\euclideanNorm{\forceOneObliqProjMat } M^{(0)}_{\forceMat} \sqrt{\sizeKinematicFE} \timestepk{\timeIndex}  \right] + \\
& \quad \quad  \dfrac{C_4}{64} \left[
2C_5 \sizeWholeFE (M_{\forceOne}^{(2)} + M_{\forceTv}^{(2)}) + 
\left(\euclideanNorm{\ROMThermodynamicMassMat^{-1} } \euclideanNorm{\forceTvObliqProjMat}
M^{(0)}_{\forceMat} + 1\right) \euclideanNorm{\ROMKinematicMassMat^{-1} } \euclideanNorm{\forceOneObliqProjMat }
\sizeWholeFE M_{\forceOne}^{(2)} \timestepk{\timeIndex} \right] \\
& \quad \quad
\left[2\tilde{M}^{(0)}_{\forceSysApprox} + 
\left(\euclideanNorm{\ROMThermodynamicMassMat^{-1} }  \euclideanNorm{\forceTvObliqProjMat } 
M^{(0)}_{\forceMat} + 1 \right) \euclideanNorm{\ROMKinematicMassMat^{-1} } 
\euclideanNorm{\forceOneObliqProjMat } M^{(0)}_{\forceMat} \sqrt{\sizeKinematicFE} \timestepk{\timeIndex}  \right]^2,
\end{aligned}
\end{equation}
where the constants $M^{(0)}_{\forceMat}$ and $\tilde{M}^{(0)}_{\forceSys}$ are defined as
\begin{equation}
\begin{split}
M^{(0)}_{\forceMat} & = \max_{0 \leq \dummyTimeIndex \leq \ntimestepROM-1} 
\left\{\euclideanNorm{ \forceMat(\stateApproxt{\dummyTimeIndex})}, 
\euclideanNorm{ \forceMat(\stateApproxt{\dummyTimeIndex+\frac{1}{2}})}\right\} \\
\tilde{M}^{(0)}_{\forceSysApprox} & = \max_{0 \leq \dummyTimeIndex \leq \ntimestepROM-1} 
\euclideanNorm{\forceSysApprox(\stateApproxt{\dummyTimeIndex})}. 
\end{split}
\end{equation}
Finally, for the first term on the right hand side of \eqref{eq:rkCoeff-3}, from \eqref{eq:res-def-c}, we see that 
\begin{equation}\label{eq:res-def-c2}
\fullNorm{(\velocityBasis\velocityResidual{1}, \energyBasis\energyResidual{1}, \positionBasis\positionResidual{1})} 
\leq \dfrac{1}{6} \max_{0 \leq \dummyTimeSymbol \leq \finalTime} 
\fullNorm{\left(\dfrac{d^3 \velocityApprox}{d\timeSymbol^3}(\dummyTimeSymbol), \dfrac{d^3 \energyApprox}{d\timeSymbol^3}(\dummyTimeSymbol), \dfrac{d^3 \positionApprox}{d\timeSymbol^3}(\dummyTimeSymbol) \right)},
\end{equation}
where we have used that fact that $\stateApprox$ is $C^3$ on $[0,\finalTime]$, since 
$\forceMat$ is of class $C^2$ on $\kinematicFE \times \thermodynamicFE \times \kinematicFE$ 
and so are $\forceOne$ and $\forceTv$.
Combining \eqref{eq:rkCoeff-3-3} and \eqref{eq:res-def-c2}, we conclude that, 
\begin{equation}
\label{eq:derror-3b}
\begin{aligned}
\fullNorm{(\rkVelocityErrorCoeff{3}, \rkEnergyErrorCoeff{3}, \rkPositionErrorCoeff{3})}
& \leq A_0 + A_1 \timestep + A_2 (\timestep)^2 + A_3 (\timestep)^3, 
\end{aligned}
\end{equation}
where the coefficients $A_0, A_1, A_2, A_3$ are given by 
\begin{equation}
\begin{aligned}
A_0 
& = \dfrac{1}{6} \max_{0 \leq \dummyTimeSymbol \leq \finalTime} 
\fullNorm{\left(\dfrac{d^3 \velocityApprox}{d\timeSymbol^3}(\dummyTimeSymbol), \dfrac{d^3 \energyApprox}{d\timeSymbol^3}(\dummyTimeSymbol), \dfrac{d^3 \positionApprox}{d\timeSymbol^3}(\dummyTimeSymbol) \right)} + \\
& \quad \quad \dfrac{C_4 C_5}{4} \left( M^{(1)}_{\forceOne} + M^{(1)}_{\forceTv} \right)
\left(\euclideanNorm{\ROMThermodynamicMassMat^{-1} }  \euclideanNorm{\forceTvObliqProjMat } 
M^{(0)}_{\forceMat} + 1 \right)
\euclideanNorm{\forceOneObliqProjMat } M^{(0)}_{\forceMat} \sqrt{\sizeKinematicFE}  + \\
& \quad \quad \dfrac{C_4}{4}  
\left(\euclideanNorm{\ROMThermodynamicMassMat^{-1} } \euclideanNorm{\forceTvObliqProjMat}
M^{(0)}_{\forceMat} + 1\right) \euclideanNorm{\ROMKinematicMassMat^{-1} }  \euclideanNorm{\forceOneObliqProjMat }M^{(1)}_{\forceOne} \tilde{M}^{(0)}_{\forceSysApprox} + 
\dfrac{C_4 C_5}{8} \sizeWholeFE (M_{\forceOne}^{(2)} + M_{\forceTv}^{(2)}) \left(\tilde{M}^{(0)}_{\forceSysApprox}\right)^2 \\
A_1 & = \dfrac{C_4}{8}  
\left(\euclideanNorm{\ROMThermodynamicMassMat^{-1} } \euclideanNorm{\forceTvObliqProjMat}
M^{(0)}_{\forceMat} + 1\right)^2 \euclideanNorm{\ROMKinematicMassMat^{-1} }^2  \euclideanNorm{\forceOneObliqProjMat }^2 
M^{(1)}_{\forceOne} M^{(0)}_{\forceMat} \sqrt{\sizeKinematicFE} + \\ 
& \quad \quad  \dfrac{C_4}{32} 
\left(\euclideanNorm{\ROMThermodynamicMassMat^{-1} } \euclideanNorm{\forceTvObliqProjMat}
M^{(0)}_{\forceMat} + 1\right) \euclideanNorm{\ROMKinematicMassMat^{-1} } \euclideanNorm{\forceOneObliqProjMat }
\sizeWholeFE M_{\forceOne}^{(2)}\tilde{M}^{(0)}_{\forceSysApprox} + \\
& \quad \quad \dfrac{C_4 C_5}{8} \sizeWholeFE (M_{\forceOne}^{(2)} + M_{\forceTv}^{(2)})
\tilde{M}^{(0)}_{\forceSysApprox} 
\left(\euclideanNorm{\ROMThermodynamicMassMat^{-1} }  \euclideanNorm{\forceTvObliqProjMat } 
M^{(0)}_{\forceMat} + 1 \right) \euclideanNorm{\ROMKinematicMassMat^{-1} } 
\euclideanNorm{\forceOneObliqProjMat } M^{(0)}_{\forceMat} \sqrt{\sizeKinematicFE} \\
A_2 & = \dfrac{C_4 C_5}{32} \sizeWholeFE (M_{\forceOne}^{(2)} + M_{\forceTv}^{(2)}) 
\left(\euclideanNorm{\ROMThermodynamicMassMat^{-1} }  \euclideanNorm{\forceTvObliqProjMat } 
M^{(0)}_{\forceMat} + 1 \right)^2 \euclideanNorm{\ROMKinematicMassMat^{-1} }^2 
\euclideanNorm{\forceOneObliqProjMat }^2 \left(M^{(0)}_{\forceMat}\right)^2 \sizeKinematicFE + \\
& \quad \quad \dfrac{C_4}{16} 
\left(\euclideanNorm{\ROMThermodynamicMassMat^{-1} }  \euclideanNorm{\forceTvObliqProjMat } 
M^{(0)}_{\forceMat} + 1 \right)^2 \euclideanNorm{\ROMKinematicMassMat^{-1} }^2 
\euclideanNorm{\forceOneObliqProjMat }^2 \sizeWholeFE M_{\forceOne}^{(2)} 
\tilde{M}^{(0)}_{\forceSysApprox} M^{(0)}_{\forceMat} \sqrt{\sizeKinematicFE} \\
A_3 & = \dfrac{C_4}{64}
\left(\euclideanNorm{\ROMThermodynamicMassMat^{-1} } \euclideanNorm{\forceTvObliqProjMat}
M^{(0)}_{\forceMat} + 1\right)^3 \euclideanNorm{\ROMKinematicMassMat^{-1} }^3 \euclideanNorm{\forceOneObliqProjMat }^3
\sizeWholeFE M_{\forceOne}^{(2)} \timestepk{\timeIndex} \left(M^{(0)}_{\forceMat}\right)^2 \sizeKinematicFE.
\end{aligned}
\end{equation}
Substituting \eqref{eq:derror-3b} back to \eqref{eq:derror-3}, we see that 
\begin{equation}\label{eq:derror-3c}
\begin{aligned}
    \fullNorm{(\discreteVelocityErrort{\timeIndex+1}, \discreteEnergyErrort{\timeIndex+1}, \discretePositionErrort{\timeIndex+1})} 
    + A (\timestep)^2 & \leq (1 + C_1 \timestepk{\timeIndex} + C_3 (\timestepk{\timeIndex})^2)
    \left(\fullNorm{(\discreteVelocityErrort{\timeIndex}, \discreteEnergyErrort{\timeIndex}, \discretePositionErrort{\timeIndex})} 
    +  A (\timestep)^2 \right) \\
    & \leq e^{(C_1+ C_3 \timestep)  \timestepk{\timeIndex}}
    \left(\fullNorm{(\discreteVelocityErrort{\timeIndex}, \discreteEnergyErrort{\timeIndex}, \discretePositionErrort{\timeIndex})} 
    +  A (\timestep)^2 \right), 
\end{aligned}
\end{equation}
where $A = C_1^{-1} (A_0 + A_1 \timestep + A_2(\timestep)^2 + A_3 (\timestep)^3)$.
By induction, we have 
\begin{equation}\label{eq:derror-4}
\begin{aligned}
    \fullNorm{(\discreteVelocityErrort{\timeIndex}, \discreteEnergyErrort{\timeIndex}, \discretePositionErrort{\timeIndex})} + A (\timestep)^2
    \leq e^{(C_1+ C_3 \timestep) \timek{\timeIndex}}
    \left(\fullNorm{(\discreteVelocityErrort{0}, \discreteEnergyErrort{0}, \discretePositionErrort{0})} +  A (\timestep)^2 \right), 
\end{aligned}
\end{equation}
which implies 
\begin{equation}\label{eq:derror-4b}
\begin{aligned}
    \fullNorm{(\discreteVelocityErrort{\timeIndex}, \discreteEnergyErrort{\timeIndex}, \discretePositionErrort{\timeIndex})}
    & \leq e^{(C_1+ C_3 \timestep) \timek{\timeIndex}}
    \fullNorm{(\discreteVelocityErrort{0}, \discreteEnergyErrort{0}, \discretePositionErrort{0})} + 
    \left(e^{(C_1+ C_3 \timestep) \timek{\timeIndex}} - 1\right) A (\timestep)^2 \\
    & \leq e^{(C_1+ C_3 \timestep) \timek{\timeIndex}} 
    \left(\fullNorm{(\discreteVelocityErrort{0}, \discreteEnergyErrort{0}, \discretePositionErrort{0})} + 
     A (\timestep)^2 (C_1 + C_3 \timestep) \timek{\timeIndex} \right), 
\end{aligned}
\end{equation}
which yields the desired result. 
\end{proof}
\end{theorem}

\section{Numerical experiments}\label{sec:numericalresults}
In this section, we present some numerical results to test the performance of
our proposed method.  Our ROM is applied to several Lagrangian hydrodynamics
problems that can be simulated with Laghos\footnote{GitHub page, {\it
https://github.com/CEED/Laghos/tree/rom}.} and libROM\footnote{GitHub page, {\it
https://github.com/LLNL/libROM}.}. Simple command-line options for simulating
these problems are provided in Appendix~\ref{sec:appendix} for the purpose of
reproducible research.  The relative error for each ROM field is measured
against the corresponding FOM solution at the final time $\finalTime$, which is
defined as:
\begin{align}\label{eq:relerrors}
  \relerrorVelocityt{\finalTime} &= \frac{\| \velocityt{\ntimestep} -
  \velocityApproxt{\ntimestepROM} \|_2 }{\| \velocityt{\ntimestep} \|_2},& 
  \relerrorEnergyt{\finalTime} &= \frac{\| \energyt{\ntimestep} - 
  \energyApproxt{\ntimestepROM} \|_2 }{\| \energyt{\ntimestep} \|_2},& 
  \relerrorPositiont{\finalTime} &= \frac{\| \positiont{\ntimestep} -
  \positionApproxt{\ntimestepROM} \|_2 }{\| \positiont{\ntimestep} \|_2}.
\end{align}
The speed-up of each ROM simulation is measured by dividing the wall-clock time
for the ROM time loop by the wall-clock time for the corresponding FOM time
loop.  A visualization tool, VisIt \cite{ahern2013visit}, is used to plot the
meshes and solution fields both for FOM and ROM.  All the simulations in this
numerical section use Quartz in Livermore Computing Center\footnote{High
performance computing at LLNL, https://hpc.llnl.gov/hardware/platforms/Quartz},
on Intel Xeon CPUs with 128 GB memory, peak TFLOPS of 3251.4, and peak single
CPU memory bandwidth of 77 GB/s. For simplicity, we only report serial tests,
but we have also observed good speed-up for parallel FOM simulation and serial ROM
simulation.

\subsection{Problem setting} 
In this subsection, we introduce the settings of several benchmark
problems, including the Gresho vortex problem, the Sedov Blast problem, the
Taylor--Green vortex problem, and the triple-point problem.

\subsubsection{Gresho vortex problem}\label{sec:gresho}
The Gresho vortex problem is a two-dimensional standard benchmark test 
for the incompressible inviscid Navier--Stokes equations \cite{gresho1990theory}. 
In this problem, we consider a manufactured smooth solution from extending the 
steady state Gresho vortex solution to the compressible Euler equations.  The
computational domain is the unit square $\initialDomain = [-0.5,0.5]^2$ with wall
boundary conditions on all surfaces, i.e. $\velocitySymbol \cdot \normalSymbol =
0$.  Let $(r, \phi)$ denote the polar coordinates of a particle
$\initialPosition \in \initialDomain$.  The initial angular velocity is given by
\begin{equation}
\velocitySymbol_\phi = 
\begin{cases}
5r & \text{ for } 0 \leq r < 0.2 \\
2-5r & \text{ for } 0.2 \leq r < 0.4 \\
0 & \text{ for } r \geq 0.4.
\end{cases}
\end{equation}
The initial density is given by $\densitySymbol = 1$. 
The initial thermodynamic pressure is given by 
\begin{equation}
p = \begin{cases}
5 + \frac{25}{2} r^2 & \text{ for } 0 \leq r < 0.2 \\
9 - 4 \log(0.2) + \frac{25}{2} - 20r + 4 \log(r) &  \text{ for } 0.2 \leq r < 0.4 \\
3 + 4\log(2) & \text{ for } r \geq 0.4. 
\end{cases}
\end{equation}
The initial energy is related to the pressure and the density by \eqref{eq:EOS}.
The adiabatic index in the ideal gas equations of state is set to be a constant
$\adiabaticIndexSymbol = 5/3$.  The initial mesh is a uniform Cartesian
hexahedral mesh, which deforms over time.  No artificial viscosity stress is
added.  We compute the resultant source terms for driving the time-dependent
simulation up to some point in time, and perform normed error analysis on the
final computational mesh. The visualized solution of the Gresho vortex problem
is given in the first column of Fig.~\ref{fig:long-fom}.

\subsubsection{Sedov Blast problem}\label{sec:sedov}
The Sedov blast problem is a three-dimensional standard shock hydrodynamic benchmark test
\cite{sedov1993similarity}. In this problem, we consider a delta source of
internal energy initially deposited at the origin of a three-dimensional cube,
analogous to the two-dimensional test in \cite{dobrev2012high}, to which we
refer for further details.  The computational domain is the unit cube
$\initialDomain = [0,1]^3$ with wall boundary conditions on all surfaces, i.e.
$\velocitySymbol \cdot \normalSymbol = 0$.  The initial velocity is given by
$\velocitySymbol = 0$.  The initial density is given by $\densitySymbol = 1$.
The initial energy is given by a delta function at the origin.  In our
implementation, the delta function energy source is approximated by setting the
internal energy to zero in all degrees of freedom except at the origin.  The
default value of the initial internal energy at the origin is
$\energySymbol(0,0,0) = 0.25$.  The adiabatic index in the ideal gas equations
of state is set to be a constant $\adiabaticIndexSymbol = 1.4$.  The initial
mesh is a uniform Cartesian hexahedral mesh, which deforms over time.
Artificial viscosity stress is added.  We compute the resultant source terms for
driving the time-dependent simulation up to some point in time, and perform
normed error analysis on the final computational mesh. 
The visualized solution of the Sedov blast problem
is given in the second column of Fig.~\ref{fig:long-fom}. 
It can be seen that the radial symmetry is maintained 
in the shock wave propagation, thanks to a high-order artificial viscosity formulation.

\subsubsection{Taylor--Green vortex problem}\label{sec:taylorGreen}

The Taylor--Green vortex problem is a three-dimensional standard flow benchmark test 
for the incompressible inviscid Navier--Stokes equations \cite{taylor1937mechanism}. 
In this problem, we consider a manufactured smooth solution from extending the 
steady state Taylor--Green vortex solution to the compressible Euler equations, 
analogous to the two-dimensional test in \cite{dobrev2012high}, to which we
refer for further details.  The computational domain is the unit cube
$\initialDomain = [0,1]^3$ with wall boundary conditions on all surfaces, i.e.
$\velocitySymbol \cdot \normalSymbol = 0$.  The initial velocity is given by 
\begin{equation}
\velocitySymbol = ( \sin(\pi x) \cos(\pi y) \cos(\pi z), -\cos(\pi x) \sin(\pi y) \cos(\pi z), 0 ).
\end{equation}
The initial density is given by $\densitySymbol = 1$. 
The initial thermodynamic pressure is given by 
\begin{equation}
\pressureSymbol = 100 + \dfrac{(\cos(2\pi x) + \cos(2\pi y))(\cos(2 \pi z) + 2) - 2}{16}. 
\end{equation}
The initial energy is related to the pressure and the density by \eqref{eq:EOS}.
The adiabatic index in the ideal gas equations of state is set to be a constant
$\adiabaticIndexSymbol = 5/3$.  The initial mesh is a uniform Cartesian
hexahedral mesh, which deforms over time.  Artificial viscosity stress is added.
We compute the resultant source terms for driving the time-dependent simulation
up to some point in time, and perform normed error analysis on the final
computational mesh. 
The visualized solution of the Taylor--Green vortex problem is given in the
third column of Fig.~\ref{fig:long-fom}.

\subsubsection{Triple--point problem}\label{sec:triple}

\begin{figure}[ht!]
\centering
\includegraphics[width=0.24\linewidth]{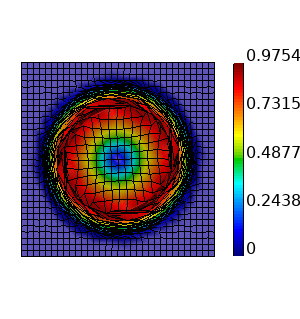}
\includegraphics[width=0.24\linewidth]{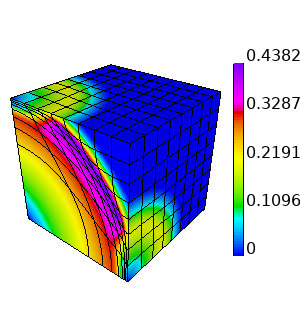}
\includegraphics[width=0.24\linewidth]{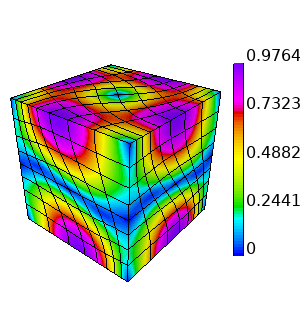}
\includegraphics[width=0.24\linewidth]{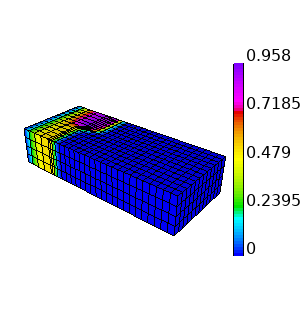}\\
\includegraphics[width=0.24\linewidth]{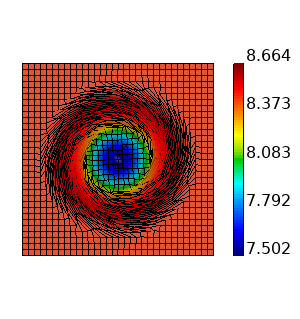}
\includegraphics[width=0.24\linewidth]{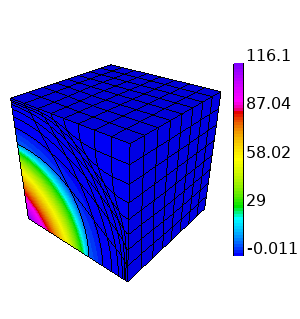}
\includegraphics[width=0.24\linewidth]{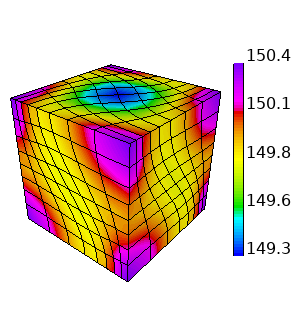}
\includegraphics[width=0.24\linewidth]{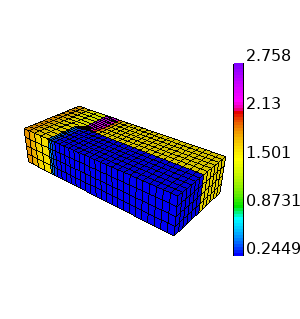}\\
\caption{Visualization of final-time velocity (top) and energy (bottom) on the
  Lagrangian frame for long-time FOM simulation in various benchmark problems:
  the Gresho vortex problem (first column), the Sedov Blast problem (second
  column), the Taylor--Green vortex problem (third column) and the triple--point
  problem (fourth column).}
\label{fig:long-fom}
\end{figure}

The triple-point problem is a three-dimensional shock test with two materials in
three states \cite{galera2010two}.  In this problem, we consider a three-state,
two-material, 2D Riemann problem which generates vorticity, analogous to the
two-dimensional test in \cite{dobrev2012high}, to which we refer for further
details.  The computational domain is $\initialDomain = [0,7] \times
[0,3] \times [0,1.5]$ with wall boundary conditions on all surfaces, i.e.
$\velocitySymbol \cdot \normalSymbol = 0$.  The initial velocity is given by
$\velocitySymbol = 0$.  The initial density is given by 
\begin{equation}
\densitySymbol = 
\begin{cases} 
1 & \text{ if } x \leq 1 \text{ or } y \leq 1.5, \\
1/8 & \text{ if } x > 1 \text{ and } y > 1.5.
\end{cases}
\end{equation}
The initial thermodynamic pressure is given by 
\begin{equation}
\pressureSymbol = 
\begin{cases}
1 & \text{ if } x \leq 1, \\
0.1 & \text{ if } x > 1. 
\end{cases}
\end{equation}
The initial energy is related to the pressure and the density by \eqref{eq:EOS}. 
The adiabatic index in the ideal gas equations of state is set to be 
\begin{equation}
\adiabaticIndexSymbol = 
\begin{cases}
1.5 & \text{ if } x \leq 1 \text{ or } y > 1.5, \\
1.4 & \text{ if } x > 1 \text{ and } y \leq 1.5.  
\end{cases}
\end{equation}
The initial mesh is a uniform Cartesian hexahedral mesh, which deforms over
time.  Artificial viscosity stress is added.  We compute the resultant source
terms for driving the time-dependent simulation up to some point in time, and
perform normed error analysis on the final computational mesh. 
The visualized solution of the triple point problem is given in the
fourth column of Fig.~\ref{fig:long-fom}.

We remark that the FOM simulation for all tests is characterized
by several FOM user-defined input values, namely the CFL constant $\CFLconst$ in
\eqref{eq:adaptive-timestep}, the mesh refinement level $\meshLevel$ which
controls the mesh size $\meshSize = 2^{-\meshLevel} \meshSize_0$ where
$\meshSize_0$ is the coarsest mesh size, and the polynomial order $\feOrder$
used in the $Q_\feOrder$--$Q_{\feOrder-1}$ pair for finite element
discretization in \cite{dobrev2012high}.  In particular, the dimensions
$\sizeKinematicFE$ and $\sizeThermodynamicFE$ of the FOM state space increase
with the mesh size $\meshSize$ and the polynomial order $\feOrder$.
Table~\ref{tab:fom} summarizes the FOM user-defined input values and
discretization parameters in various benchmark problems, which characterizes the
FOM simulation in the later sections.
The resultant system of semi-discrete Lagrangian hydrodynamics equations 
is numerically solved by the RK2-average scheme 
introduced in Section~\ref{sec:RK2-avg}.
\begin{table}[ht!]
\centering
\begin{tabular}{|c||c|c|c|c|}
\hline
Problem & Gresho vortex & Sedov Blast & Taylor--Green vortex & Triple--point  \\ 
\hline
$\CFLconst$ & 0.5 & 0.5 & 0.1 & 0.5 \\
$\meshLevel$ & 4 & 2 & 2 & 2 \\
$\feOrder$ & 3 & 2 & 2 & 2 \\ 
\hline
$\sizeKinematicFE$ & 18818 & 14739 & 14739 & 38475 \\
$\sizeThermodynamicFE$ & 9216 & 4096 & 4096 & 10752 \\ 
\hline
\end{tabular}
\caption{FOM user-defined input values and discretization parameters in various
  benchmark problems.}
\label{tab:fom}
\end{table}

\begin{figure}[ht!]
\centering
\includegraphics[width=0.24\linewidth]{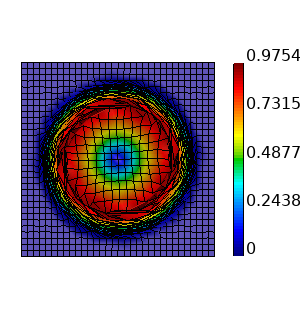}
\includegraphics[width=0.24\linewidth]{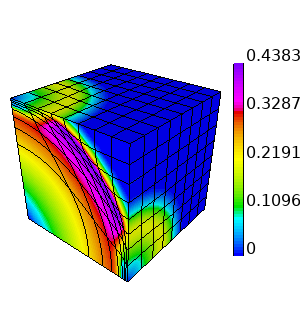}
\includegraphics[width=0.24\linewidth]{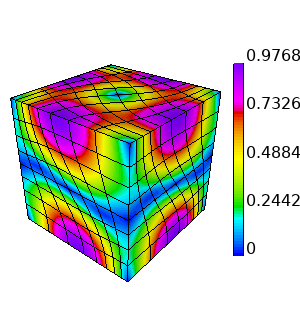}
\includegraphics[width=0.24\linewidth]{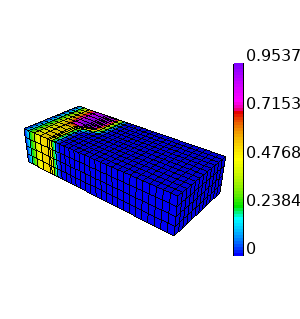}\\
\includegraphics[width=0.24\linewidth]{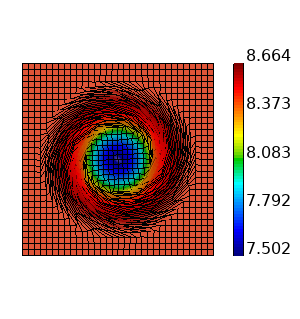}
\includegraphics[width=0.24\linewidth]{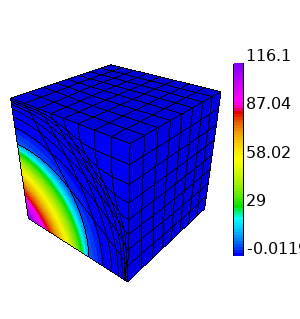}
\includegraphics[width=0.24\linewidth]{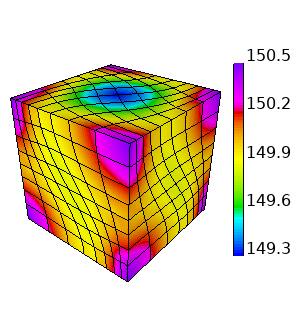}
\includegraphics[width=0.24\linewidth]{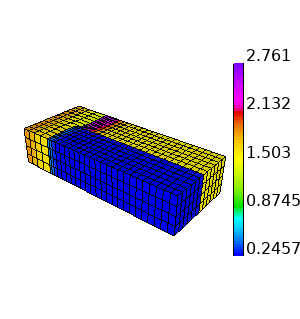}\\
\caption{Visualization of final-time velocity (top) and energy (bottom) on the
  Lagrangian frame for long-time ROM simulation in various benchmark problems:
  the Gresho vortex problem (first column), the Sedov Blast problem (second
  column), the Taylor--Green vortex problem (third column) and the triple--point
  problem (fourth column).}
\label{fig:long-twrom}
\end{figure}

At the end of this section, we demonstrate some representative 
simulation results on the long-time simulation in various benchmark problems, 
which will be further discussed in Section~\ref{sec:exp-long}.
The final time $\finalTime$ of the simulation is taken to be 0.62 for the Gresho
vortex problem, 0.8 for the Sedov Blast problem, 0.25 for the Taylor--Green
vortex problem and 0.8 for the triple--point problem.  Figure~\ref{fig:long-fom}
shows the visualization of final-time velocity and energy on the Lagrangian
frame for long-time FOM simulation in various benchmark problems. It can be
seen that the Lagrangian mesh is extensively distorted and compressed in the
Gresho vortex problem and the Sedov Blast problem, which suggests that sharp
gradients are developed in these problems. Figure~\ref{fig:long-twrom}
shows the visualization of final-time velocity and energy on the Lagrangian
frame for long-time ROM simulation in various benchmark problems using the
time-windowing approach. It can be seen that the ROM simulation results are in
good agreement with the FOM simulation results in Figure~\ref{fig:long-fom}.
The ROM simulation results successfully capture the extreme mesh distortion.

\subsection{Short-time ROM simulation}\label{sec:exp-short}

In this section, we use the spatial ROM introduced in Section~\ref{sec:ROM} to
perform short-time simulation in various benchmark problems.  Only reproductive
cases are considered, i.e.  the problem setting in the ROM is identical to that
in the FOM.  First, in the offline phase, the fully discrete FOM scheme in
Section~\ref{sec:time_integrator} is used to compute snapshot solution for
performing POD discussed in Section~\ref{sec:POD}.  The reduced basis matrices
for the solution variables are then used to formulate the reduced order model
\eqref{eq:rom} in the online phase. In the case of hyper-reduction, we also
obtain the reduced basis matrices for the nonlinear terms using SVD or the SNS
relation, and then follow the DEIM procedure discussed in Section~\ref{sec:DEIM}
to obtain sampling matrices and formulate the hyper-reduced system
\eqref{eq:rom-hr}.  Using a time integrator again, the fully discrete
(hyper-)reduced order model is used to compute an approximate solution.  The ROM
solution is compared with the FOM solution at the final time $\finalTime$ with
respect to the relative errors \eqref{eq:relerrors}.  Also, the wall time of ROM
simulation is compared with that of the FOM simulation.  We will study the
dependence of the accuracy and speed-up of ROM on the ratios of dimension
reduction 
\begin{equation}
\ratioROMvelocity = \dfrac{\sizeROMvelocity}{\sizeKinematicFE}, 
\quad 
\ratioROMenergy = \dfrac{\sizeROMenergy}{\sizeThermodynamicFE},
\quad 
\ratioROMposition = \dfrac{\sizeROMposition}{\sizeKinematicFE},
\quad
\ratioROMforceOneSample = \dfrac{\sizeROMforceOneSample}{\sizeKinematicFE}, 
\quad 
\ratioROMforceTvSample = \dfrac{\sizeROMforceTvSample}{\sizeThermodynamicFE} . 
\end{equation}
On the other hand, for hyper-reduction, the numbers of sampling indices 
$\{\sizeROMforceOneSample, \sizeROMforceTvSample\}$ are controlled by the product of 
the over-sampling factors $\{\factorROMforceOneSample, \factorROMforceTvSample\} \in [1, \infty)$ 
and the reduced basis dimensions $\{\sizeROMforceOne, \sizeROMforceTv\}$, i.e. 
\begin{equation}
\sizeROMforceOneSample = \min\left\{ \sizeKinematicFE,  \, \factorROMforceOneSample\sizeROMforceOne \right\}, 
\quad 
\sizeROMforceTvSample = \min\left\{ \sizeThermodynamicFE, \, \factorROMforceTvSample\sizeROMforceTv \right\}. 
\end{equation}
The resultant reduced system of semi-discrete Lagrangian hydrodynamics equations 
is numerically solved by the RK2-average scheme 
introduced in Section~\ref{sec:rom_time_integrator}.

\begin{figure}[ht!]
\centering
\includegraphics[width=0.45\linewidth]{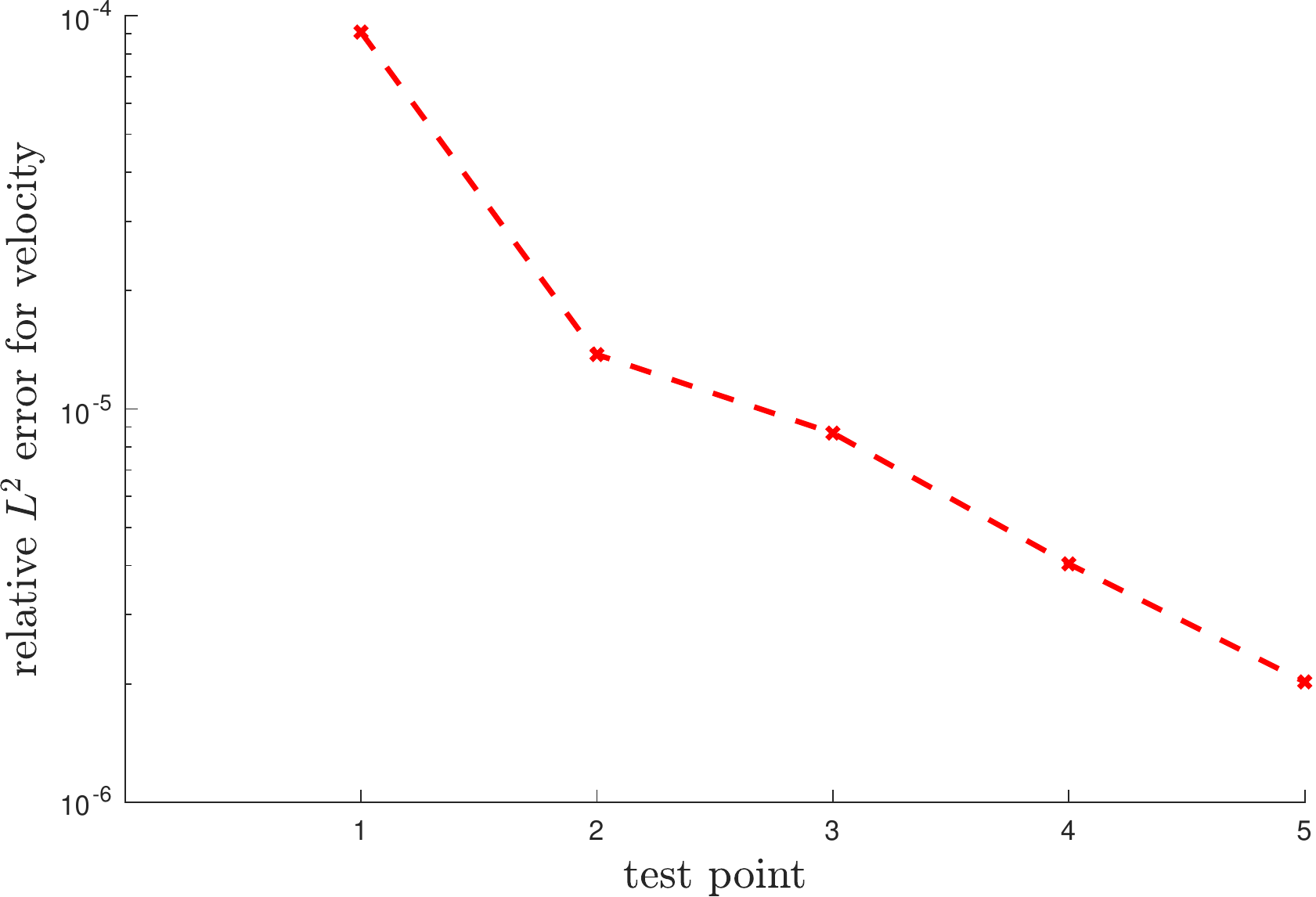}
\hspace{0.05\linewidth}
\includegraphics[width=0.45\linewidth]{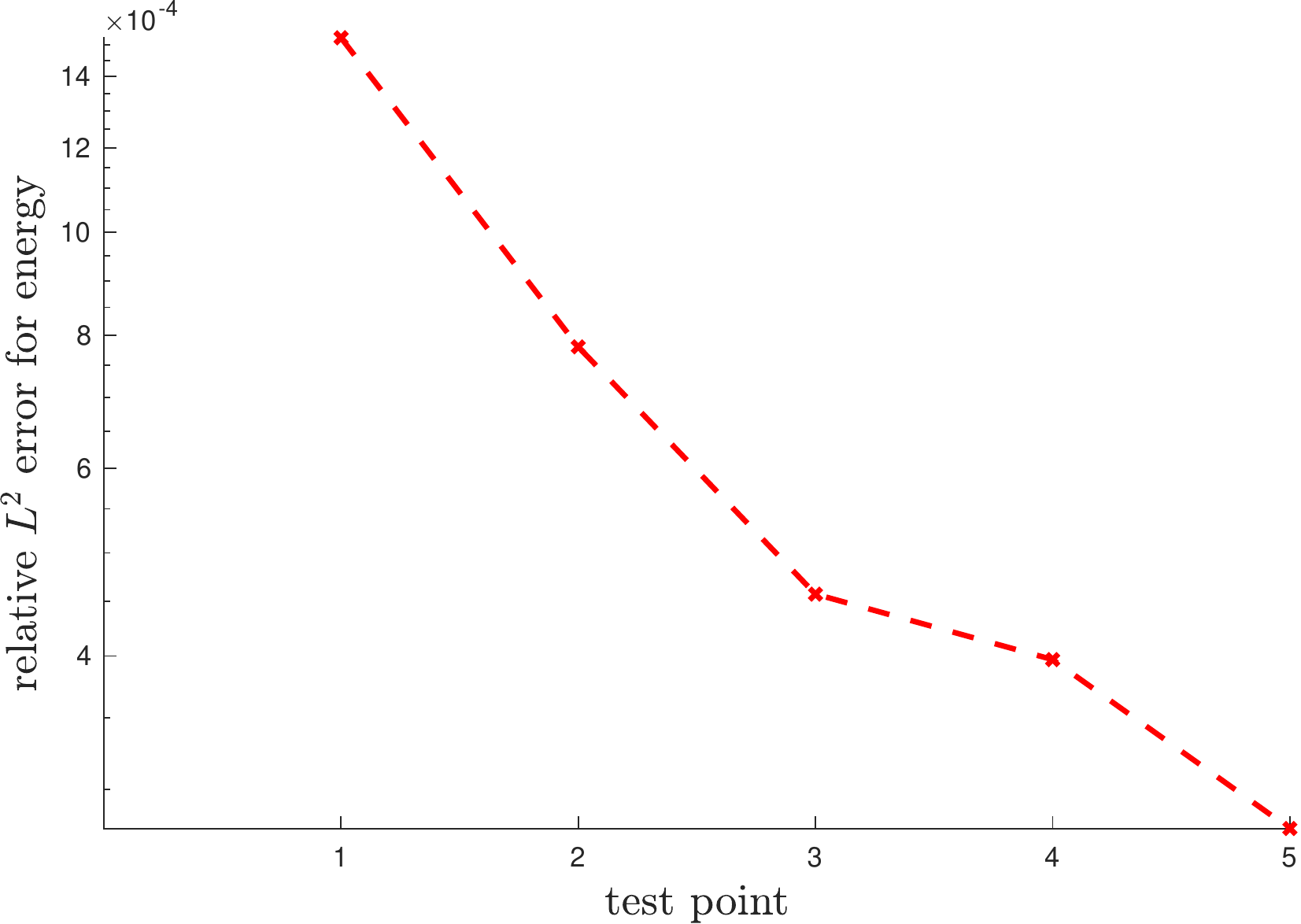}\\
\vspace{0.2in}
\includegraphics[width=0.45\linewidth]{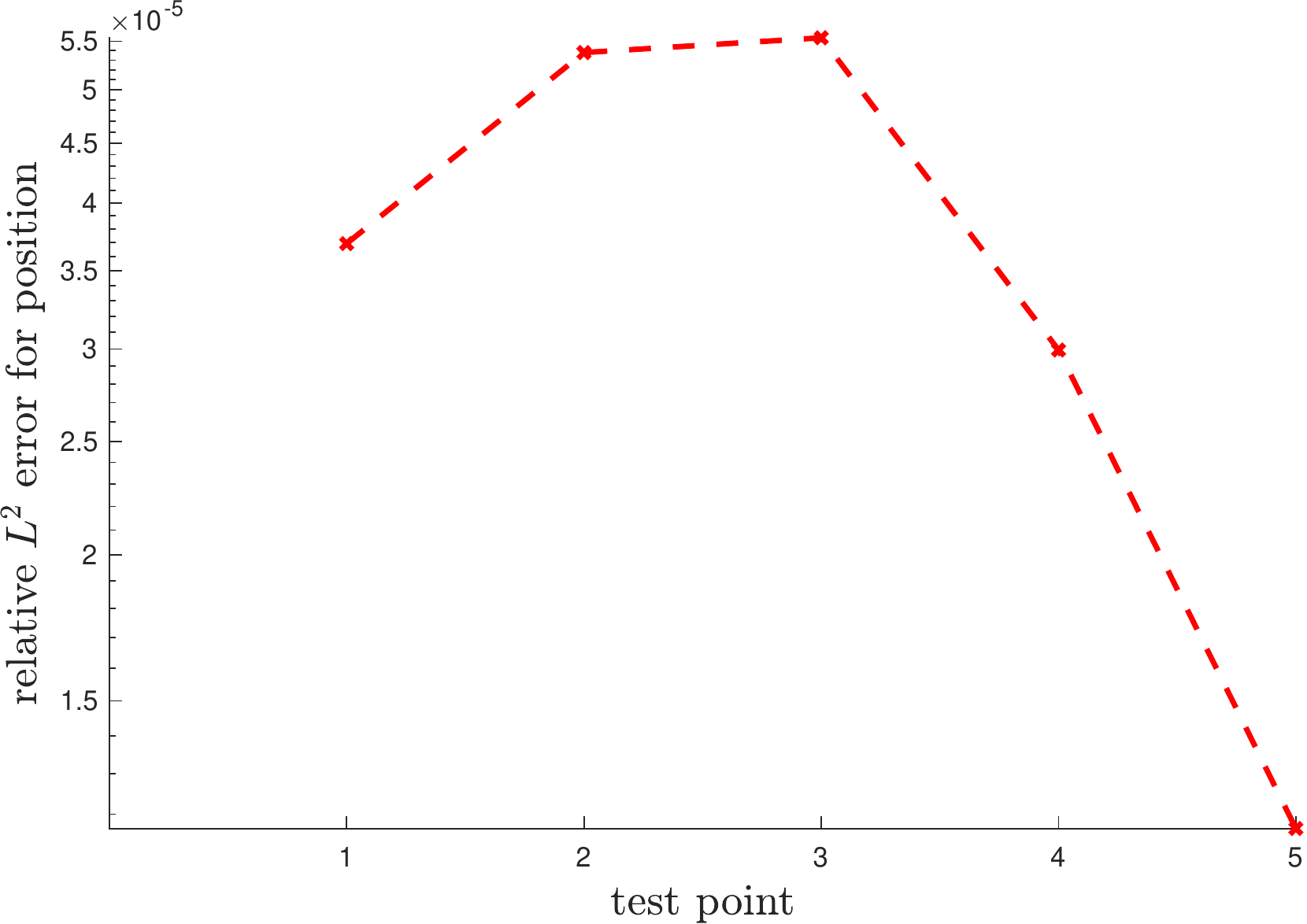}
\hspace{0.05\linewidth}
\includegraphics[width=0.45\linewidth]{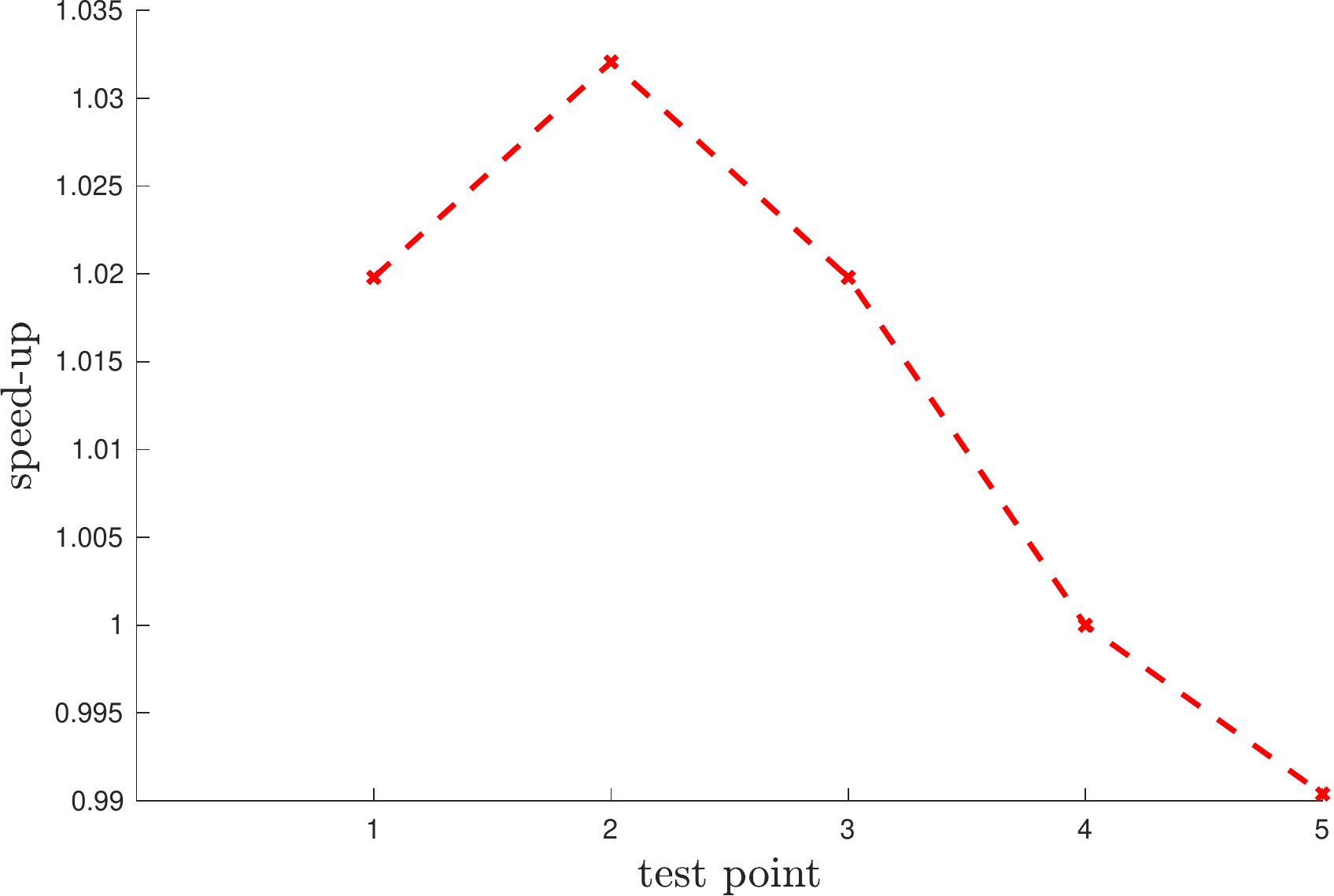} \\
\caption{ROM performance comparison for short-time simulation in 
the Gresho vortex problem with varying ratio of reduced basis dimensions: 
relative $L^2$ error for velocity (top-left), 
relative $L^2$ error for energy (top-right), 
relative $L^2$ error for position (bottom-left) and 
speed-up (bottom-right). Table~\ref{tab:rdim_gresho_short} provides test points. 
  }
\label{fig:rdim_gresho_short}
\end{figure}

In Figure~\ref{fig:rdim_gresho_short}, we compare the ROM performance on the
short-time simulation in the Gresho vortex problem against reduced basis
dimensions.  This comparison examines how the solution representability of the
reduced order model depends on the dimensions of the POD reduced bases. To this
end, we fix all other ROM parameters and do not perform hyper-reduction, to avoid
complicating the results.  Table~\ref{tab:rdim_gresho_short} shows the
reduced basis dimensions being tested and their ratios to the dimensions of the
corresponding FOM finite element spaces.  From
Figure~\ref{fig:rdim_gresho_short}, we see that the relative errors are overall
decreasing with the reduced basis dimensions, and the accuracy is very
outstanding even when low dimensional subspaces are used. However, since the
nonlinear terms scale with the FOM size, the speed-up achieved is only around
1.4 for all the cases tested.  This highlights the need of hyper-reduction
for the sake of remarkable speed-up.

\begin{table}[ht!]
\centering
\begin{tabular}{|c||c|c|c|c|c|}
\hline
Test point & 1 & 2 & 3 & 4 & 5  \\ 
\hline
$\sizeROMvelocity$ & 6 & 11 & 17 & 23 & 28 \\ 
$\sizeROMenergy$ & 9 & 18 & 27 & 36 & 44 \\ 
$\sizeROMposition$ & 2 & 4 & 6 & 8 & 10 \\ 
\hline
$\ratioROMvelocity$ 
& 0.000318844 & 0.000584547 & 0.00090339 & 0.00122223 & 0.00148794 \\ 
$\ratioROMenergy$ 
& 0.000976562 & 0.00195312 & 0.00292969 & 0.00390625 & 0.00477431 \\ 
$\ratioROMposition$ 
& 0.000106281 & 0.000212562 & 0.000318844 & 0.000425125 & 0.000531406 \\ 
\hline
\end{tabular}
\caption{List of reduced basis dimensions being tested in
  Figure~\ref{fig:rdim_gresho_short} and their ratio to the dimensions of the
  corresponding FOM finite element spaces for short-time simulation in the
  Gresho vortex problem.}\label{tab:rdim_gresho_short}
\end{table}

\begin{figure}[ht!]
\centering
\includegraphics[width=0.45\linewidth]{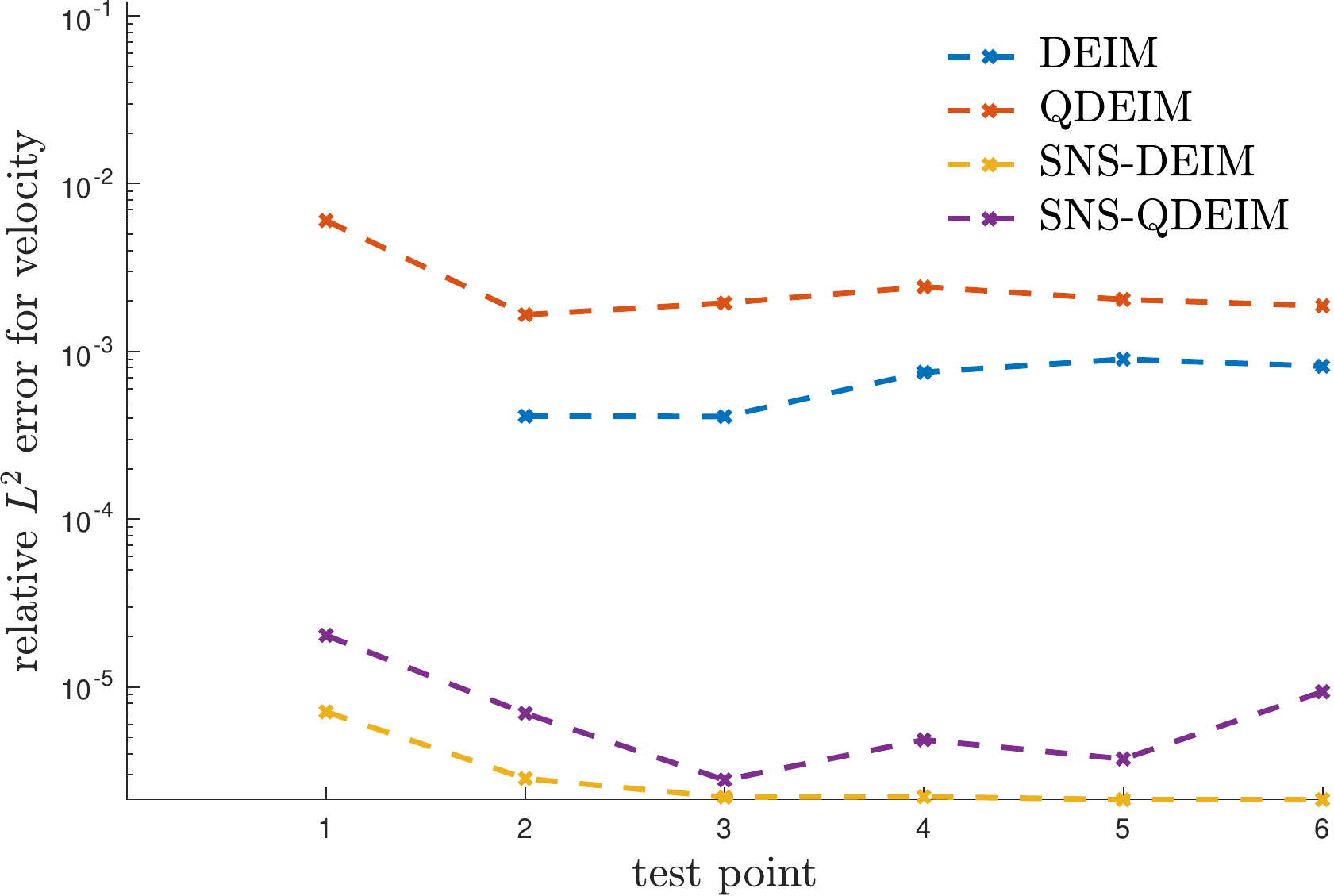}
\hspace{0.05\linewidth}
\includegraphics[width=0.45\linewidth]{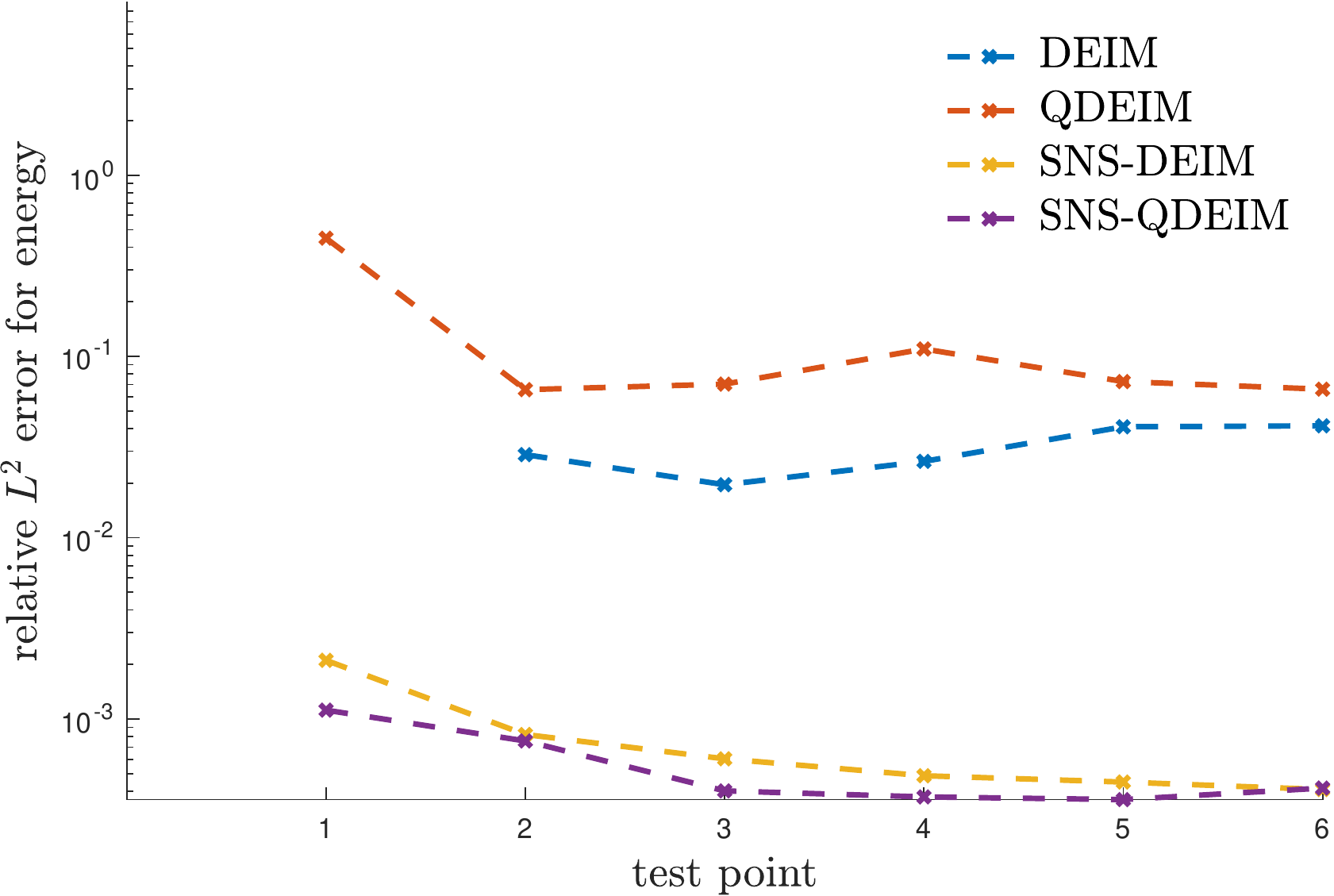}\\
\vspace{0.2in}
\includegraphics[width=0.45\linewidth]{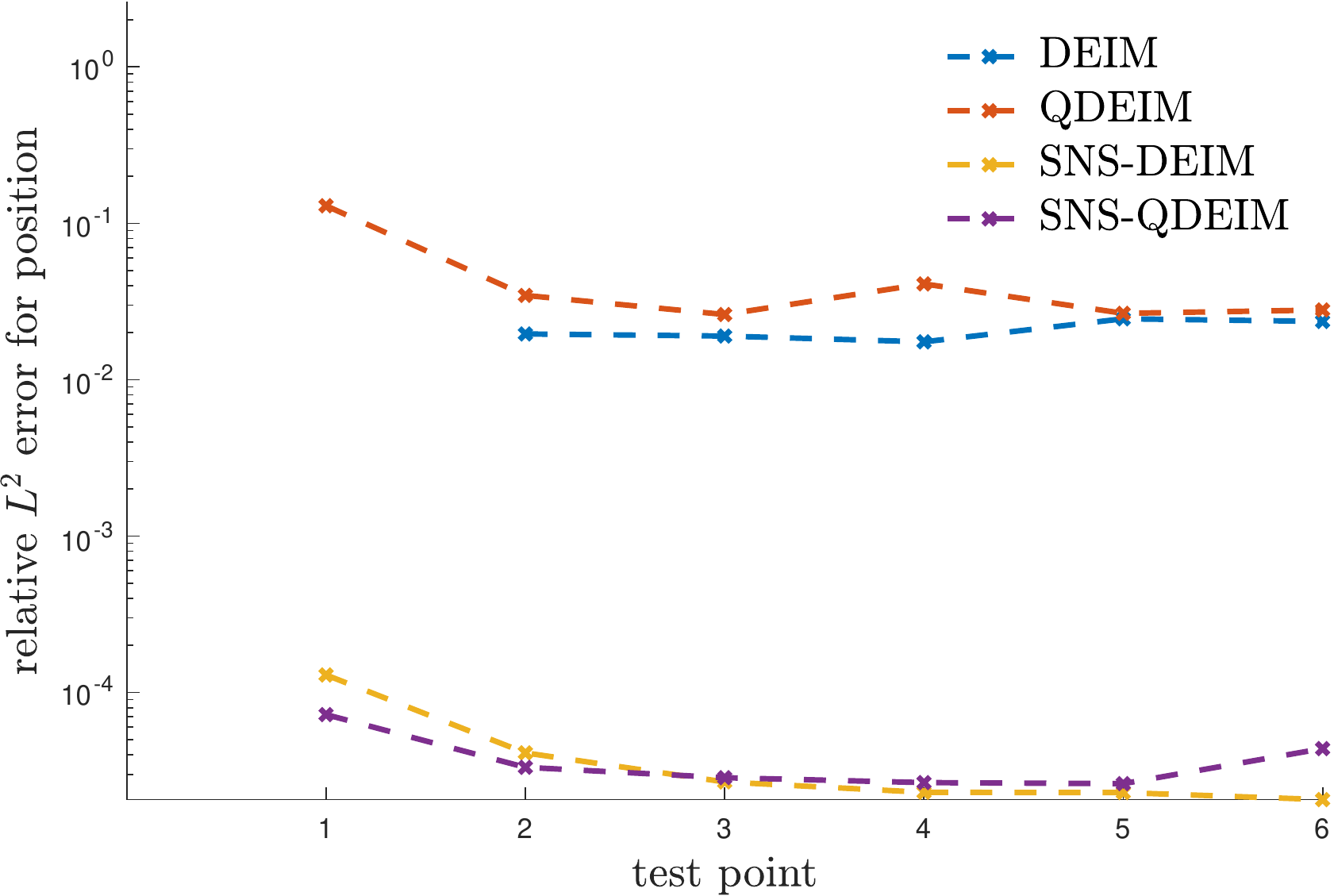}
\hspace{0.05\linewidth}
\includegraphics[width=0.45\linewidth]{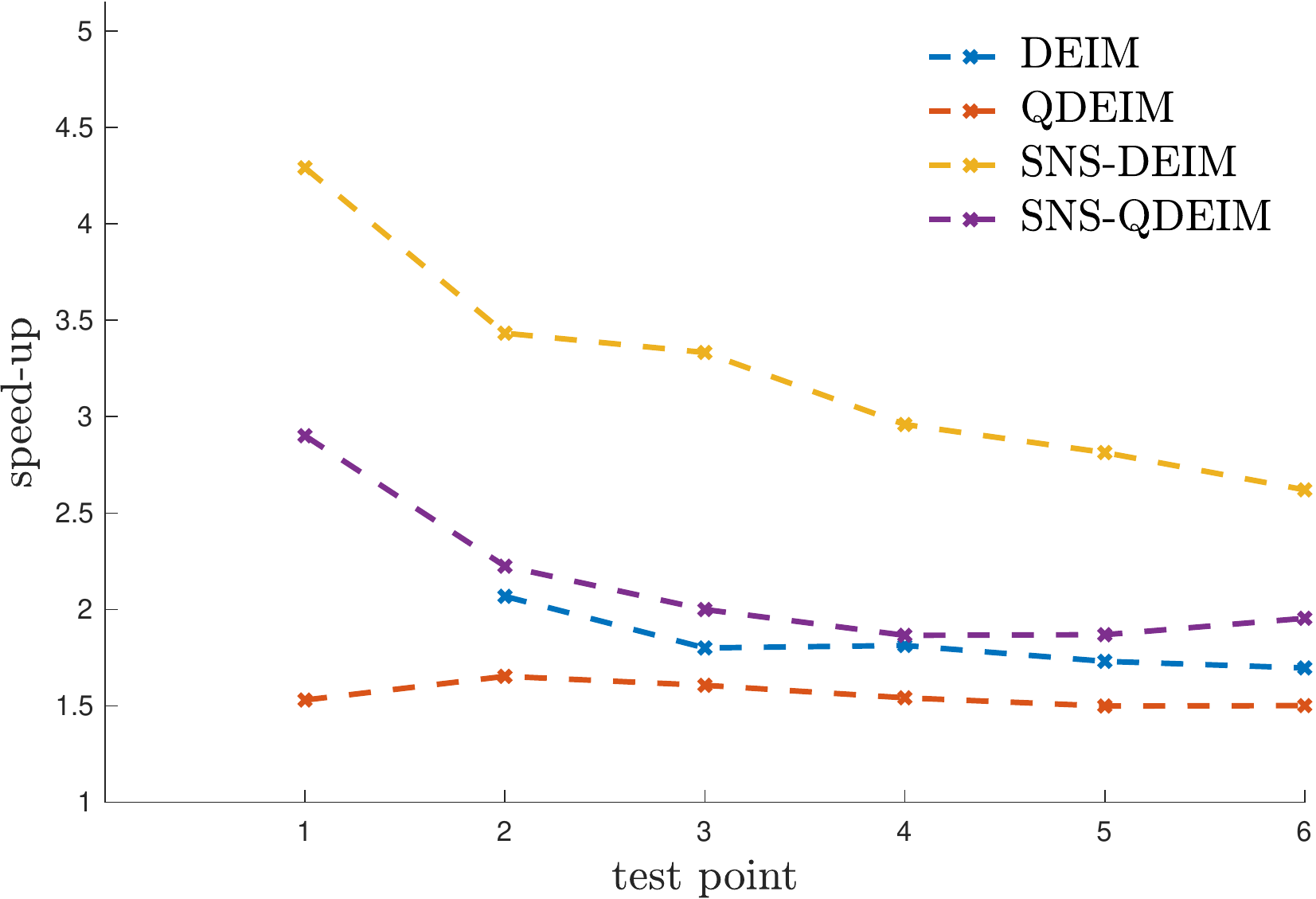} \\
\caption{ROM performance comparison for short-time simulation in 
the Gresho vortex problem with varying ratio of sampling indices: 
relative $L^2$ error for velocity (top-left), 
relative $L^2$ error for energy (top-right), 
relative $L^2$ error for position (bottom-left) and 
speed-up (bottom-right). Table~\ref{tab:sfac_gresho_short} provides test points. 
  }
\label{fig:sfac_gresho_short}
\end{figure}

In Figure~\ref{fig:sfac_gresho_short}, we compare the ROM performance on the
short-time simulation in the Gresho vortex problem against the number of
sampling indices.  This comparison examines how the solution accuracy depends on
the projection error of the nonlinear terms in the DEIM nonlinear model
reduction.  We compare different approaches for obtaining the nonlinear term
bases, i.e.  directly applying snapshot SVD and using the SNS relation
\cite{choi2020sns}.  To avoid complicating the results with respect to
different ROM parameters, we fix the reduced basis dimensions as
$(\sizeROMvelocity, \sizeROMenergy, \sizeROMposition, \sizeROMforceOne,
\sizeROMforceTv) = (28, 44, 10, 28, 44)$.  We also compare different algorithms
for obtaining the sampling indices, namely the over-sampling DEIM (see Algorithm
3 of \cite{carlberg2013gnat} and Algorithm 5 of \cite{carlberg2011efficient})
and QDEIM (see \cite{drmac2016new}). Table~\ref{tab:sfac_gresho_short} shows
the number of sampling indices being tested and their ratio to the dimensions of
the corresponding FOM finite element spaces. 
From Figure~\ref{fig:sfac_gresho_short}, we see that the relative errors and speed-up
are overall decreasing with the number of sampling indices, while the accuracy
is still remarkable for each of the nonlinear model reduction technique. 
We remark that the first test point with DEIM has a numerical instability 
as the time step size eventually vanishes.
In this experiment, we observe that SNS-DEIM has the
highest accuracy and best speed-up compared with the other approaches. In the
rest of the paper, unless otherwise specified, we will mainly focus on the
SNS-DEIM approach as suggested by this experiment. 

\begin{table}[ht!]
\centering
\begin{tabular}{|c||c|c|c|c|c|c|}
\hline
Test point & 1 & 2 & 3 & 4 & 5 & 6 \\ 
\hline
$\factorROMforceOneSample$ & 8 & 15 & 23 & 30 & 37 & 45 \\ 
$\factorROMforceTvSample$ & 7 & 14 & 21 & 28 & 35 & 41 \\ 
\hline

$\sizeROMforceOneSample$ & 616 & 1232 & 1848 & 2464 & 3080 & 3696 \\ 
$\sizeROMforceTvSample$ & 396 & 792 & 1188 & 1584 & 1980 & 2376 \\ 
\hline
$\ratioROMforceOneSample$ 
& 0.0327346 & 0.0654692 & 0.0982038 & 0.130938 & 0.163673 & 0.196408 \\ 
$\ratioROMforceTvSample$ 
& 0.0429688 & 0.0859375 & 0.128906 & 0.171875 & 0.214844 & 0.257812 \\ 
\hline
\end{tabular}
\caption{List of number of sampling indices being tested in
  Figure~\ref{fig:sfac_gresho_short} and their ratio to the dimensions of the
  corresponding FOM finite element spaces for short-time simulation in the
  Gresho vortex problem.}\label{tab:sfac_gresho_short}
\end{table}

Table~\ref{tab:short} shows the short-time numerical results using SNS-DEIM
spatial ROM for various benchmark problems.  The quantities reported in
Table~\ref{tab:short} include 
\begin{itemize}
\item the final time $\finalTime$ and the number of time steps $\ntimestep$ in
  the FOM simulations, 
\item the ROM user-defined input values, namely the threshold
  $\singularValueThreshold$ and the over-sampling factors
    $\{\factorROMforceOneSample, \factorROMforceTvSample\}$, 
\item the ROM parameters, namely the reduced dimensions $\{\sizeROMvelocity,
  \sizeROMenergy, \sizeROMposition\}$ (note that with our SNS setting, we have
    $\sizeROMforceOne = \sizeROMvelocity$ and $\sizeROMforceTv =
    \sizeROMenergy$) and the number of sampling indices
    $\{\sizeROMforceOneSample, \sizeROMforceTvSample\}$, and 
\item the ROM results, namely the number of time steps $\ntimestepROM$ in the
  ROM simulations, the relative errors $\{\relerrorVelocityt{\finalTime},
    \relerrorEnergyt{\finalTime}, \relerrorPositiont{\finalTime}\}$ at the final
    time, and the relative speed-up of the ROM simulations to FOM simulations.
\end{itemize}
In each of the benchmark problems, it can be observed that the number of time
steps $\ntimestepROM$ in the ROM simulations is close to but slightly differ
from $\ntimestep$ in the FOM simulations.  The ROM simulation benefits from the
relatively cheap reduced dimension computation in each time step to achieve an
overall relative speed-up.  The relative error is small for each variable,
indicating that the ROM solutions provide accurate approximations to the FOM
solutions.

\begin{table}[ht!]
\centering
\begin{tabular}{|c||c|c|c|c|}
\hline
Problem & Gresho vortex & Sedov Blast & Taylor--Green vortex & Triple--point  \\ 
\hline
$\finalTime$ & 0.1 & 0.1 & 0.05 & 0.2 \\ 
$\ntimestep$ & 87 & 245 & 122 & {28} \\ 
\hline
$\singularValueThreshold$ & 0.9999 & 0.9999 & 0.9999 & 0.9999 \\
$\factorROMforceOneSample$ & {8} & {2} & {2} & {2} \\
$\factorROMforceTvSample$ & {8} & {2} & {2} & {2} \\
\hline
$\sizeROMvelocity$ & {32} & {53} & 9 & {10} \\
$\sizeROMenergy$ & {72} & 13 & {16} & {10} \\
$\sizeROMposition$ & 10 & {17} & 3 & {7} \\
$\sizeROMforceOneSample$ & {256} & {106} & {18} & {20} \\ 
$\sizeROMforceTvSample$ & {576} & {26} & {32} & {20} \\ 
\hline
$\ntimestepROM$ & 87 & {225} & 122 & 29 \\ 
$\relerrorVelocityt{\finalTime}$ & {1.385366e-02} & {2.74144e-03} & {9.05004e-05} & {9.41694e-03} \\ 
$\relerrorEnergyt{\finalTime}$ & {2.02558e-03} & {8.39586e-04} & {3.69690e-06} & {4.91346e-03} \\ 
$\relerrorPositiont{\finalTime}$ & {1.76761e-04} & {2.04953e-05} & {6.81088e-07} & {4.03934e-05} \\ 
speed-up & {8.27130} & {26.55553} & {14.75786} & {53.02794} \\  \hline
\end{tabular}
\caption{FOM results, ROM user-defined input values, parameters and results 
for short-time simulation in various benchmark problems.}
\label{tab:short}
\end{table}

\subsection{Long-time ROM simulation}\label{sec:exp-long}

In this section, we use the spatial ROM introduced in Section~\ref{sec:ROM} and the
time-windowing ROM introduced in Section~\ref{sec:timewindowing} to perform
long-time simulation in various benchmark problems.  Again, only reproductive
cases are considered.  For the spatial ROM, we use SNS-DEIM for hyper-reduction
and the offline-online procedure splitting is the same as in Section~\ref{sec:exp-short}.
Table~\ref{tab:long} shows the long-time numerical results using SNS-DEIM
spatial ROM for various benchmark problems.  It can be observed that the
relative errors are still reasonable, but  the relative speed-up is lower than
that in the short-time simulations.  In these simulations, large over-sampling
factors $\{\factorROMforceOneSample, \factorROMforceTvSample\}$ are used to
maintain sufficient numbers of sampling indices $\{\sizeROMforceOneSample,
\sizeROMforceTvSample\}$, for otherwise the time step size would eventually vanish
due to the adaptive time stepping control and the CFL constraints.  However,
larger numbers of sampling indices lead to more expensive overdetermined
systems \eqref{eq:sol-ls-hyper} and \eqref{eq:sol-ls-hyper-2} for the
hyper-reduction and hurt the overall speed-up in ROM simulations.  It is
especially important to point out that the spatial ROM does not achieve a
speed-up in the Gresho vortex problem since more reduced basis vectors are
required to ensure the solution representability of the reduced linear space, and
almost full sampling has to be used for the sake of numerical stability. 

\begin{table}[ht!]
\centering
\begin{tabular}{|c||c|c|c|c|}
\hline
Problem & Gresho vortex & Sedov Blast & Taylor--Green vortex & Triple--point  \\ 
\hline
$\finalTime$ & 0.62 & 0.8 & 0.25 & 0.8 \\ 
$\ntimestep$ & 1672 & {719} & 897 & 193 \\ 
\hline
$\singularValueThreshold$ & 0.9999 & 0.9999 & 0.9999 & 0.9999 \\
$\factorROMforceOneSample$ & {100} & {6} & {80} & {6} \\ 
$\factorROMforceTvSample$ & {100} & {6} & {40} & {6} \\ 
\hline
$\sizeROMvelocity$ & {165} & {225} & {40} & {20} \\
$\sizeROMenergy$ & {331} & {27} & {62} & {16} \\
$\sizeROMposition$ & 34 & 29 & 6 & 10 \\
$\sizeROMforceOneSample$ & {16500} & {1350} & {3200} & {120} \\ 
$\sizeROMforceTvSample$ & 9216 & {162} & {2480} & {96} \\ 
\hline 
$\ntimestepROM$ & {1670} & {742} & {898} & 193 \\ 
$\relerrorVelocityt{\finalTime}$ & {4.29468e-02} & {9.72044e-03} & {3.60272e-02} & {2.52978e-03} \\ 
$\relerrorEnergyt{\finalTime}$ & {6.82899e-03} & {5.95163e-04} & {3.28395e-03} & {2.36835e-03} \\ 
$\relerrorPositiont{\finalTime}$ & {7.1045e-05} & {3.67286e-04} & {2.77386e-04} & {4.72467e-05} \\ 
speed-up & {0.79646} & {3.35523} & {1.96808} & {15.29550} \\  \hline
\end{tabular}
\caption{FOM results, ROM user-defined input values, parameters and results 
for long-time simulation in various benchmark problems.}
\label{tab:long}
\end{table}

Next, we examine the performance of the time windowing ROM approach.  In our
implementation, the end points of the time windows
$\{\windowk{\windowIndex}\}_{\windowIndex=1}^{\nwindow}$ can either be
determined by physical time as in Section~\ref{sec:physical_tw}, or by the
number of samples as in Section~\ref{sec:sample_tw}.  We focus our discussion on 
the SNS-DEIM time windowing ROM approach with the end points of the time windows
determined by the number of samples. In the offline phase, 
the fully discrete RK2-average FOM scheme is first used to compute solution snapshots.
The time windows are determined according to
\eqref{eq:sample_window}.  Following Section~\ref{sec:POD_tw}, the collected
solution data are clustered into different windows for performing POD, and the
resultant reduced solution bases for velocity and energy are left-multiplied by
the corresponding mass matrix to obtain the nonlinear term bases, which in turn
determine the sampling indices by DEIM. In the online phase, the reduced basis matrices
and the sampling matrices are then used to formulate the
hyper-reduced system \eqref{eq:rom-hr-tw} in each time window. The fully
discrete hyper-reduced system follows from applying the RK2-average scheme as
discussed in Section~\ref{sec:rom_time_integrator}. 

Table~\ref{tab:tw} shows the long-time numerical results for various benchmark
problems with time windowing ROM.  In the case of decomposing time windows by
number of samples, the number of time windows $\nwindow$ is an additional ROM
parameter controlled by the number of samples per window $\ntimestepWindow$ as
an additional user-defined input.  We use $\ntimestepWindow = 10$ in each
of the benchmark problems.  Moreover, in our implementation, the over-sampling
factors $\{\factorROMforceOneSample, \factorROMforceTvSample\}$ are user-defined
constants over all time windows, which control the number of sampling indices in
the window $\windowIndex$ by 
\begin{equation}
\sizeROMforceOneSampleWindow{\windowIndex} = \min\left\{ \sizeKinematicFE,  \, 
\factorROMforceOneSample\sizeROMforceOneWindow{\windowIndex} \right\}, 
\quad 
\sizeROMforceTvSampleWindow{\windowIndex} = \min\left\{ \sizeThermodynamicFE, \, 
\factorROMforceTvSample\sizeROMforceTvWindow{\windowIndex} \right\}. 
\end{equation}
In Table~\ref{tab:tw}, the reported values of reduced basis dimensions and the
number of sampling indices are taken as the greatest among all windows, i.e. 
\begin{equation}
\sizeROMvelocity = \max_{1 \leq \windowIndex \leq \nwindow} 
\sizeROMvelocityWindow{\windowIndex}, 
\quad
\sizeROMenergy = \max_{1 \leq \windowIndex \leq \nwindow} 
\sizeROMenergyWindow{\windowIndex}, 
\quad
\sizeROMposition = \max_{1 \leq \windowIndex \leq \nwindow} 
\sizeROMpositionWindow{\windowIndex}, 
\quad
\sizeROMforceOneSample = \max_{1 \leq \windowIndex \leq \nwindow} 
\sizeROMforceOneSampleWindow{\windowIndex},  
\quad 
\sizeROMforceTvSample = \max_{1 \leq \windowIndex \leq \nwindow} 
\sizeROMforceTvSampleWindow{\windowIndex}. 
\end{equation}
It can be observed from comparing Table~\ref{tab:tw} with Table~\ref{tab:long}
that it is possible to achieve a higher relative speed-up and
improve solution accuracy using time windowing ROM. 
In particular, for the Gresho vortex problem, when the spatial ROM 
fails to achieve any relative speed-up, the time windowing approach
yields relative error or order $10^{-6}$ in all solution variables and a $25$
times speed-up. This suggests that the time windowing approach is more useful
for the advection-dominated cases, in which the intrinsic dimensions of reduced
subspaces in the spatial ROM are large. 

\begin{table}[ht!]
\centering
\begin{tabular}{|c||c|c|c|c|}
\hline
Problem & Gresho vortex & Sedov Blast & Taylor--Green vortex & Triple--point  \\ 
\hline
$\finalTime$ & 0.62 & 0.8 & 0.25 & 0.8 \\ 
$\ntimestep$ & 1672 & {719} & 897 & 193 \\ 
\hline
$\ntimestepWindow$ & 10 & 10 & 10 & 10 \\
$\singularValueThreshold$ & 0.9999 & 0.9999 & 0.9999 & 0.9999 \\
$\factorROMforceOneSample$ & {2} & {2} & {2} & {2} \\ 
$\factorROMforceTvSample$ & {2} & {2} & {2} &{2} \\ 
\hline
$\nwindow$ & {335} & {144} & {180} & {39} \\ 
$\sizeROMvelocity$ & {8} & 10 & {3} & {7} \\
$\sizeROMenergy$ & 10 & {7} & 5 & {8} \\
$\sizeROMposition$ & {6} & {7} & 2 & 5 \\
$\sizeROMforceOneSample$ & {16} & {20} & {6} & {14} \\ 
$\sizeROMforceTvSample$ & {20} & {14} & {10} & {16} \\ 
\hline
$\ntimestepROM$ & {1599} & {697} & {897} & {155} \\ 
$\relerrorVelocityt{\finalTime}$ & {3.50163e-06} & {4.58111e-04} & {1.22769e-06} & {9.31072e-04}  \\ 
$\relerrorEnergyt{\finalTime}$ & {3.51595e-07} & {2.95541e-05} & {5.35795e-08} & {6.99656e-04} \\ 
$\relerrorPositiont{\finalTime}$ & {4.28729e-07} & {3.86835e-05} & {6.95479e-09} & {3.82897e-05} \\ 
speed-up & {24.59098} & {26.47649} & {23.46642} & {75.03020} \\ 
\hline
\end{tabular}
\caption{FOM results, time windowing ROM user-defined input values, parameters
  and results for long-time simulation in various benchmark problems.}
\label{tab:tw}
\end{table}

At the end of this subsection, we compare the performance of the ROM approach 
with SNS-DEIM across different window size on
long-time simulation in the Gresho vortex problem and the Sedov Blast problem.
We will consider the spatial ROM in Section~\ref{sec:ROM} 
and the time windowing ROM in Section~\ref{sec:timewindowing}.  
The endpoints of the time windows
$\{\windowk{\windowIndex}\}_{\windowIndex=1}^{\nwindow}$ are determined by 
\eqref{eq:sample_window} with the number of samples $\ntimestepWindow$. 
More precisely, the cases under consideration in the Gresho vortex problem are: 
\begin{itemize}
\item $\nwindow = 1$ corresponds to serial ROM, 
\item $\ntimestepWindow = 10$ corresponds to time windowing ROM which generates $\nwindow = 335$ time windows, and  
\item $\ntimestepWindow = 20$ corresponds to time windowing ROM which generates $\nwindow = 168$ time windows, and  
\item $\ntimestepWindow = 40$ corresponds to time windowing ROM which generates $\nwindow = 84$ time windows. 
\end{itemize}
Meanwhile, the cases under consideration in the Sedov Blast problem are: 
\begin{itemize}
\item $\nwindow = 1$ corresponds to serial ROM, 
\item $\ntimestepWindow = 10$ corresponds to time windowing ROM which generates $\nwindow = 144$ time windows, and  
\item $\ntimestepWindow = 20$ corresponds to time windowing ROM which generates $\nwindow = 72$ time windows, and  
\item $\ntimestepWindow = 40$ corresponds to time windowing ROM which generates $\nwindow = 36$ time windows. 
\end{itemize}
Depending on the reduced dimensions $\{\sizeROMforceOne, \sizeROMforceTv\}$ of
the nonlinear term bases, each case is tested with various combinations of
over-sampling factors $(\factorROMforceOneSample, \factorROMforceTvSample)$ 
to obtain adjustable solution accuracy and speed-up.

\begin{figure}[ht!]
\centering
\includegraphics[width=0.6\linewidth]{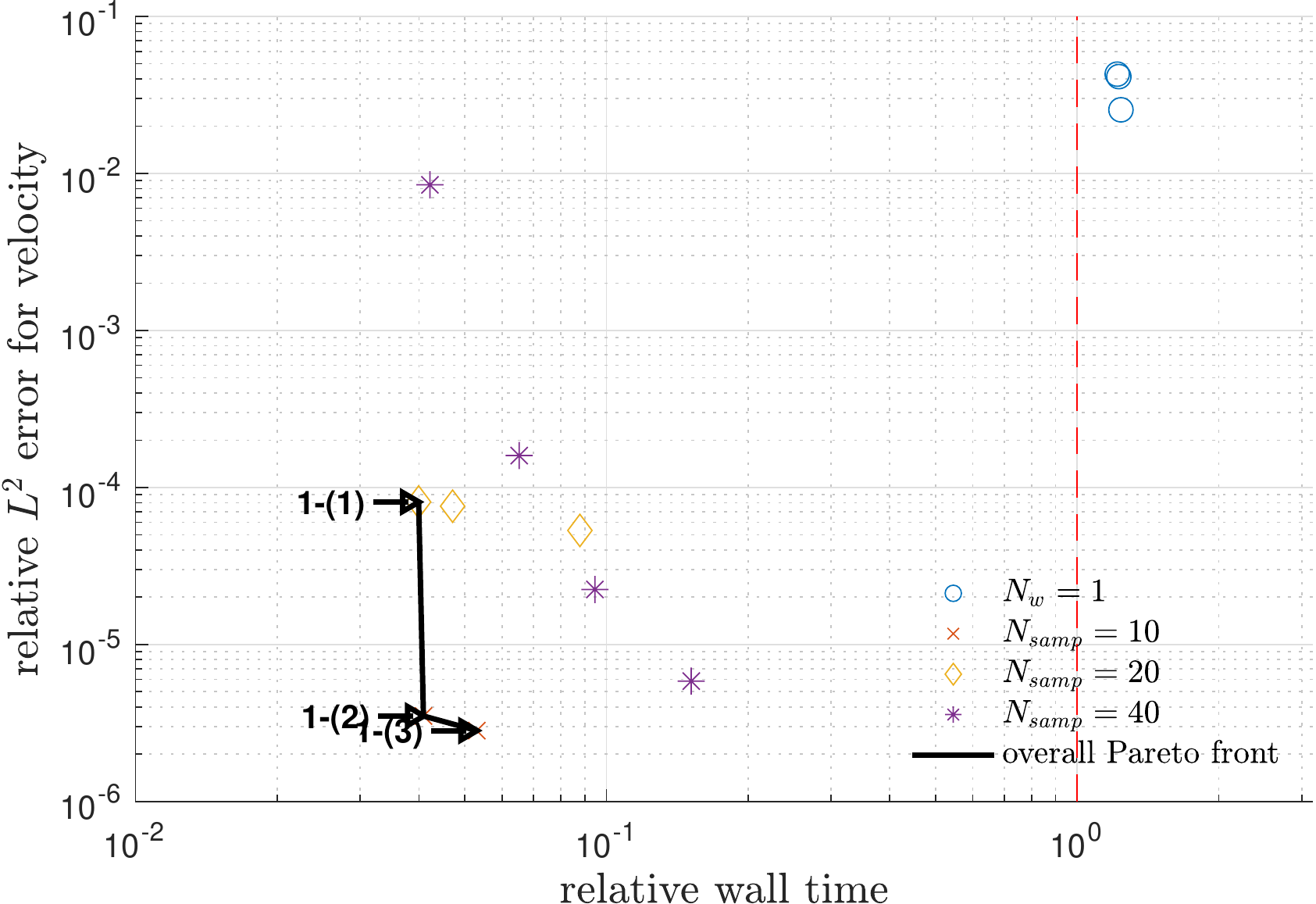} \\ 
\vspace{0.2in}
\includegraphics[width=0.6\linewidth]{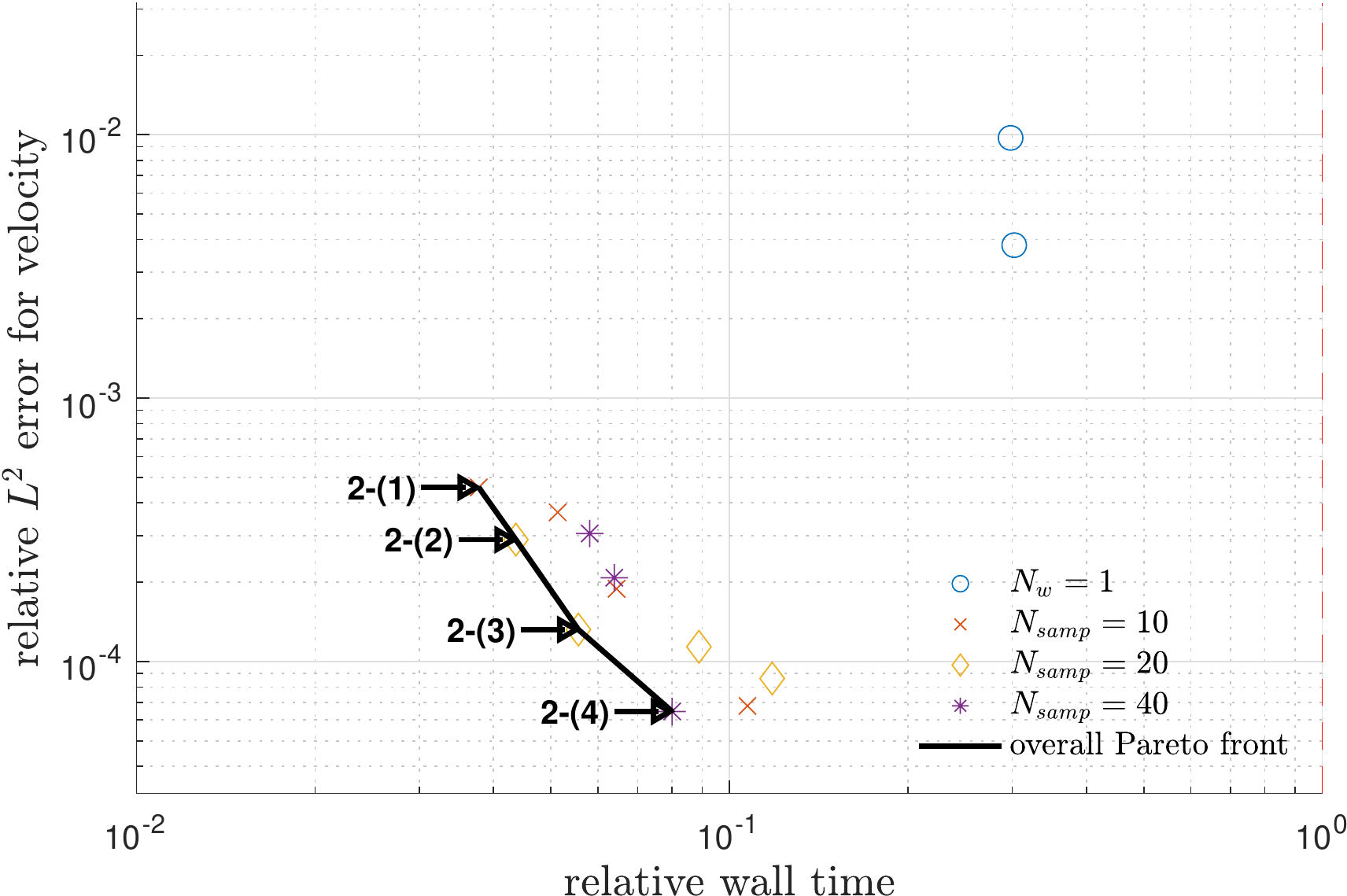} \\
\caption{ROM performance comparison for long-time simulation in the Gresho vortex problem 
(top) and the Sedov Blast problem (bottom). 
Relative $L^2$ error for velocity versus relative wall time for varying ROM parameters.
  }
\label{fig:pareto_repro}
\end{figure}

For each of the problems, we construct a Pareto front in
Figure~\ref{fig:pareto_repro} for each of the cases, which is characterized by
the ROM user-defined input values that minimize the competing objectives of
relative $L^2$ error for velocity and relative wall time.  An overall Pareto
front that selects the ROM parameters that are Pareto-optimal across all groups
is also illustrated for each problem. Table~\ref{tab:pareto_repro}
reports the ROM user-defined  input values that yielded the results on the
overall Pareto front. The results also suggest that the time windowing 
approaches can produce a remarkable
speed-up as well as accurate approximations with the appropriate use of
hyper-reduction. 

\begin{table}[ht!]
\centering
\begin{tabular}{|c||c|c|c||c|c|c|c|}
\hline
Label & 1-(1) & 1-(2) & 1-(3) & 2-(1) & 2-(2) & 2-(3) & 2-(4) \\
\hline
$\nwindow$ & -- & -- & -- & -- & -- & -- & -- \\
$\ntimestepWindow$ & 20 & 20 & 10 & 10 & 20 & 20 & 40 \\ 
$\singularValueThreshold$ & 0.9999 & 0.9999 & 0.9999 & 0.9999 & 0.9999 & 0.9999 & 0.9999  \\
$\factorROMforceOneSample$ & 2 & 2 & 4 & 2 & 2 & 4 & 8 \\ 
$\factorROMforceTvSample$ & 2 & 2 & 4 & 2 & 2 & 4 & 8 \\ 
\hline
\end{tabular}
\caption{ROM user-defined input values yielding Pareto-optimal performance for
  long-time simulation in the Gresho vortex problem and the Sedov Blast problem.
  Figure~\ref{fig:pareto_repro} provides labels.}
\label{tab:pareto_repro}
\end{table}

\subsection{Parametric ROM simulation}

In this section, we use the time-windowing ROM introduced in
Section~\ref{sec:timewindowing} to perform long-time simulation in the Sedov
Blast problem in a parametric problem setting.  In our numerical experiments,
the initial internal energy deposited at the origin is parametrized by the
problem parameter $\param \in \paramDomain = \RR{+}$, where we set
$\energySymbol(0,0,0; \param) = 0.25 \param$.  First, in the offline phase, the fully discrete FOM
scheme is used to compute snapshot solution with several problem
parameters $\{ \param_{\paramIndex} \}_{\paramIndex=1}^{\nparam}$.  Then, we
consider the SNS-DEIM time windowing ROM approach with the end points of the
time windows determined by number of samples as in Section~\ref{sec:sample_tw}.

\begin{figure}[ht!]
\centering
\includegraphics[width=0.45\linewidth]{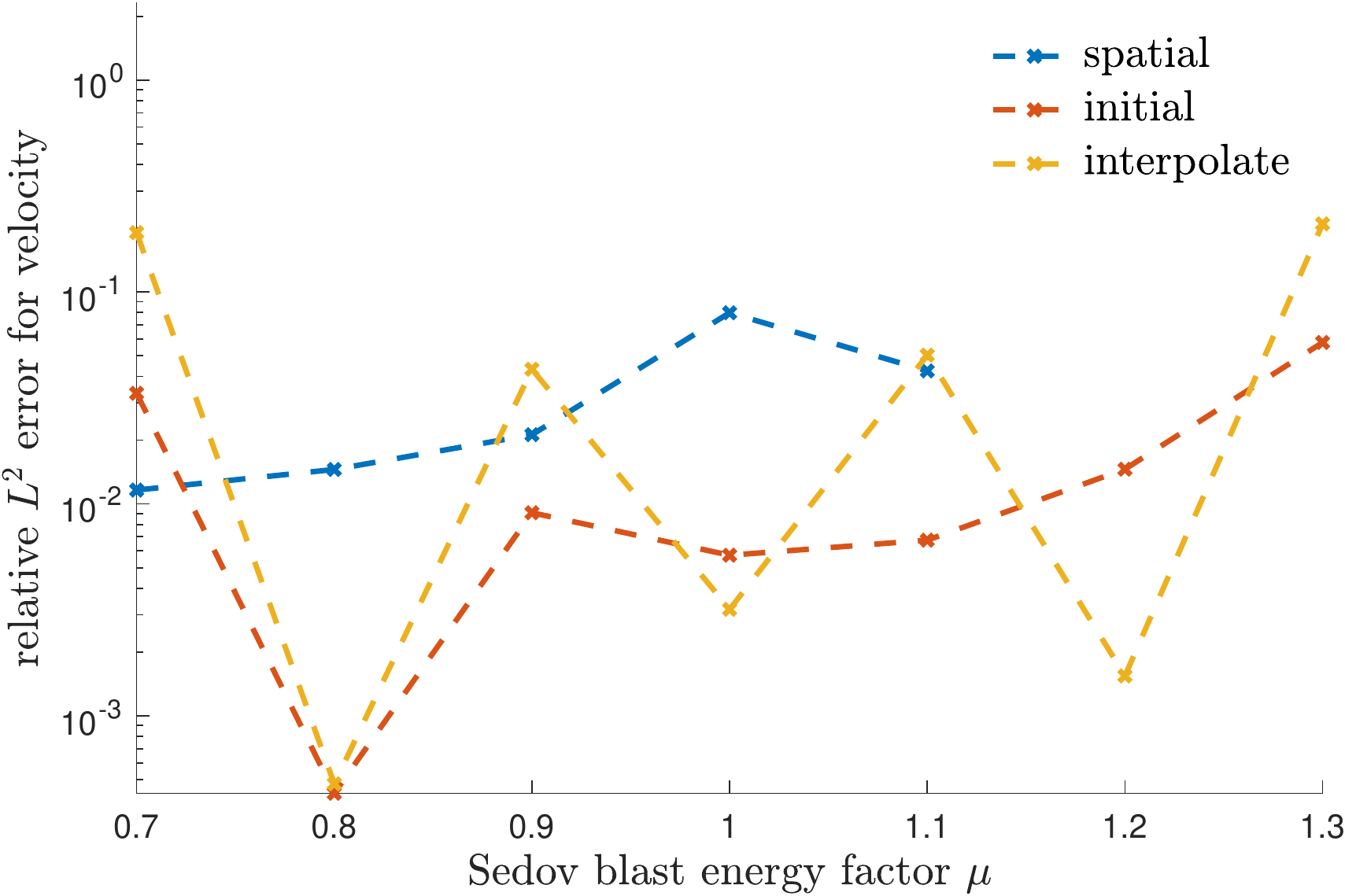}
\hspace{0.05\linewidth}
\includegraphics[width=0.45\linewidth]{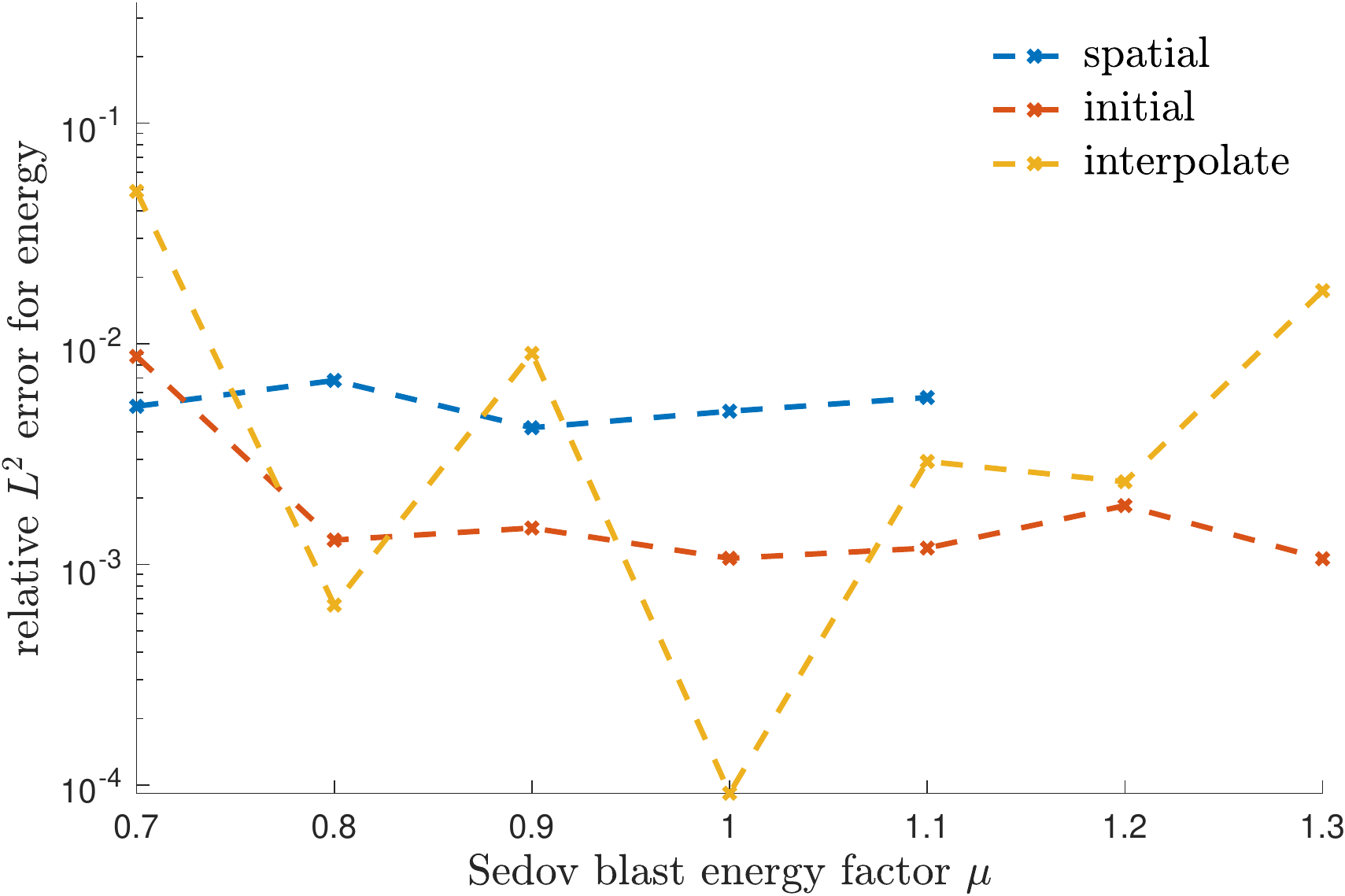}\\
\vspace{0.2in}
\includegraphics[width=0.45\linewidth]{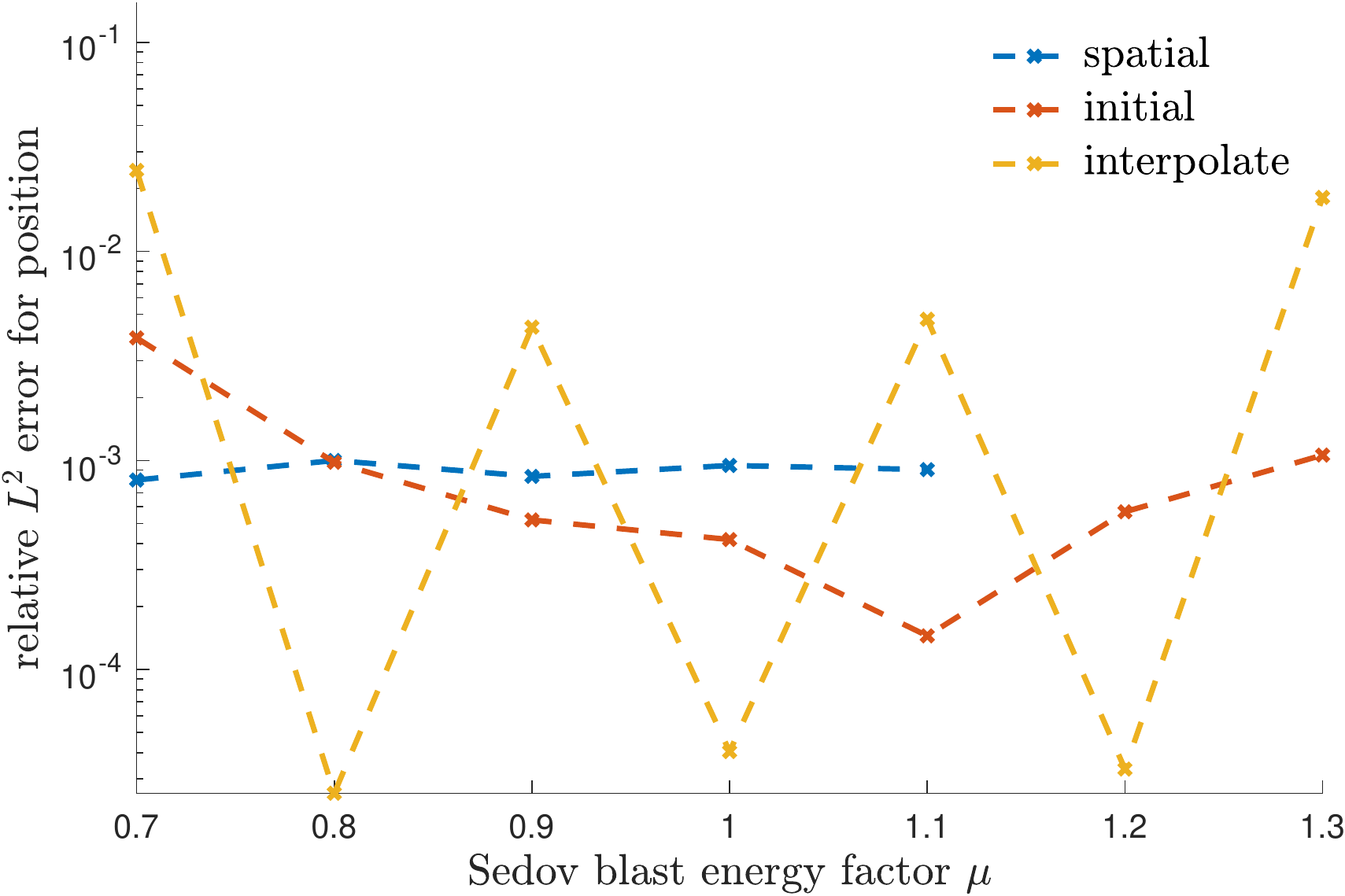}
\hspace{0.05\linewidth}
\includegraphics[width=0.45\linewidth]{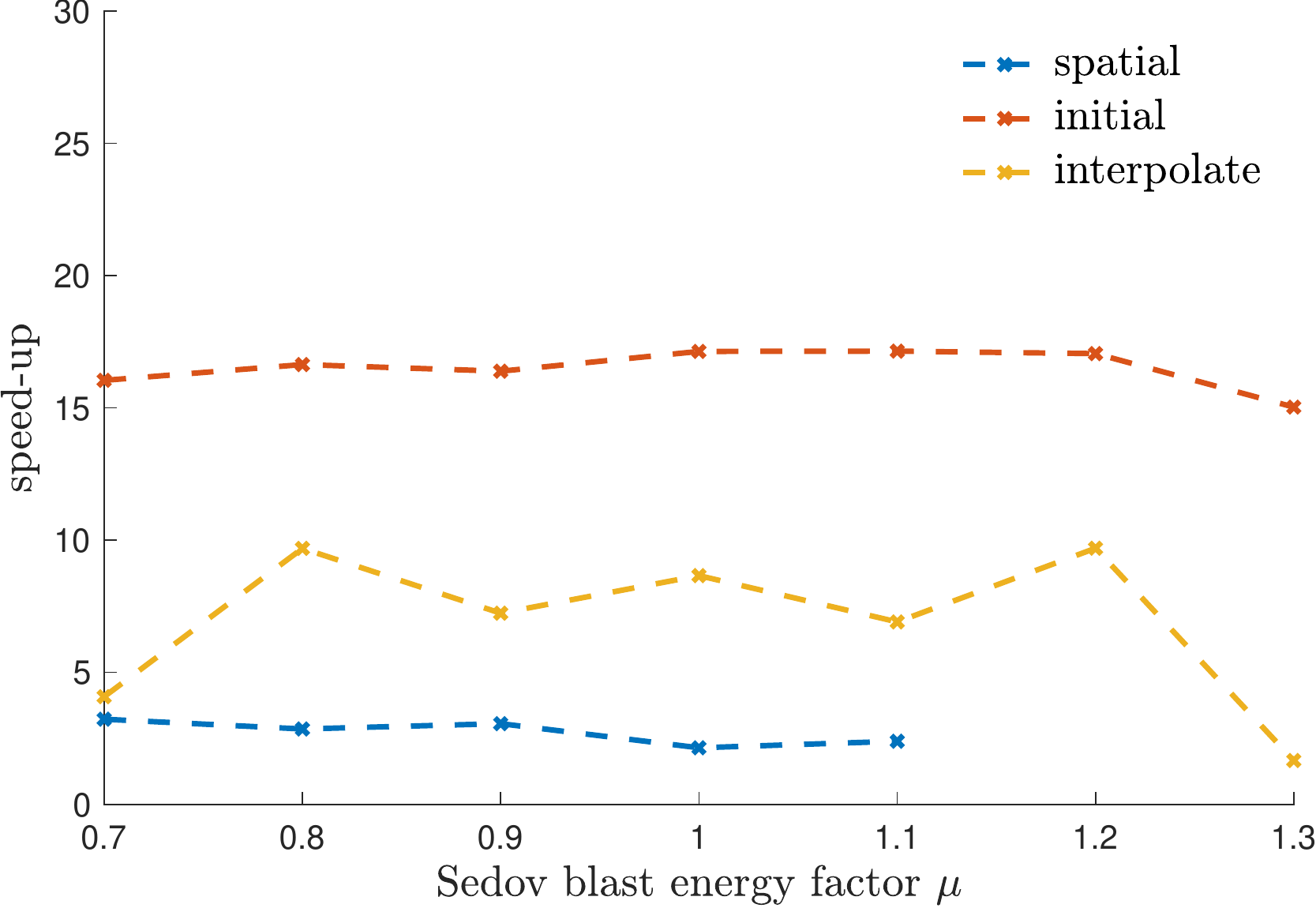} \\
\caption{Time-windowing ROM performance comparison for long-time simulation 
in the parametric Sedov Blast problem with varying problem parameter: 
relative $L^2$ error for velocity (top-left), 
relative $L^2$ error for energy (top-right), 
relative $L^2$ error for position (bottom-left) and 
speed-up (bottom-right).
  }
\label{fig:param_sedov}
\end{figure}

The FOM solutions with $\nparam = 3$ problem parameters,
namely $\{\param_1, \param_2, \param_3\} = \{0.8, 1.0, 1.2\}$, are computed and
collected as snapshots for determining the time windows. Then POD and SNS-DEIM are
performed in each window using each choice of the offset vectors as discussed
in Section~\ref{sec:offset_tw}.  The time windows are determined according to
\eqref{eq:sample_window} with $\ntimestepWindow = 20$. Then the ROM is tested
on a generic problem parameter $\param$, where the offset vectors
$\velocityOSWindow{\windowIndex}(\param), \energyOSWindow{\windowIndex}(\param),
\positionOSWindow{\windowIndex}(\param)$ are computed in each of the time
windows $\windowIndex$ according to the choice of offset vectors used in
constructing the solution bases.  In Figure~\ref{fig:param_sedov}, we compare
the ROM performance on the long-time simulation in the Sedov Blast problem
against the generic problem parameter $\param \in [0.7, 1.3]$.  This comparison
examines how the solution representability of the reduced order model depends on
the problem parameter and the choice of offset vectors.  To this end, we fix all other ROM parameters and
do not perform hyper-reduction to avoid complicating the results.  
The comparison is also made against the spatial ROM. 
As shown in Figure~\ref{fig:param_sedov}, the time windowing approach 
provides a much better speed-up than the serial ROM. 
We also remark that the serial ROM has a numerical instability 
in $\param \in \{1.2, 1.3\}$ as the time step size eventually vanishes, 
which we do not see in the time windowing ROM.
We can observe from Figure~\ref{fig:param_sedov} that using the interpolation scheme in
Section~\ref{sec:offset_interpolate} outperforms the other approaches of
computing offset vectors in terms of the relative errors in the reproductive
cases $\param \in \{0.8, 1.0, 1.2\}$, at which the offset vectors are exactly
the interpolating values for each time window.  However, 
this approach produces larger error at the predictive cases. 
This suggests that the error of offset determination by the interpolation scheme 
plays a crucial part in the final-time error. On the other hand, 
using the the initial state as offset vectors, which is described in Section~\ref{sec:offset_initial}, 
provides uniformly improved accuracy when compared with serial ROM. 
Moreover, it provides a speed-up of around 15 times at all the tested parameter values, 
which outperforms the interpolation scheme. 

At the end of this subsection, we compare the performance of various parametric
time windowing ROM approaches across a wide variation of all ROM parameters on
long-time simulation in the Sedov Blast problem in the setting of varying
the energy factor.  The FOM solutions with $\nparam = 3$ problem
parameters, namely $\{\param_1, \param_2, \param_3\} = \{0.8, 1.0, 1.2\}$, are
computed and collected as snapshots for determining the time windows, and then
performing POD and SNS-DEIM in each window using each choice of the offset
vectors as discussed in Section~\ref{sec:offset_tw}.  The time windows are
determined according to \eqref{eq:sample_window} with $\ntimestepWindow = 20$.
The ROM performance is tested on the problem parameters $\param \in \{0.7, 0.8,
0.9\}$.  Again, we compare different choices of the offset vectors as discussed
in Section~\ref{sec:offset_tw}, and the comparison is also made against the
spatial ROM.  Depending on the reduced dimensions $\{\sizeROMforceOne,
\sizeROMforceTv\}$ of the nonlinear term bases, each case is tested with various
combinations of over-sampling factors $(\factorROMforceOneSample,
\factorROMforceTvSample)$ to obtain adjustable solution accuracy and speed-up.

\begin{figure}[ht!]
\centering
\includegraphics[width=0.6\linewidth]{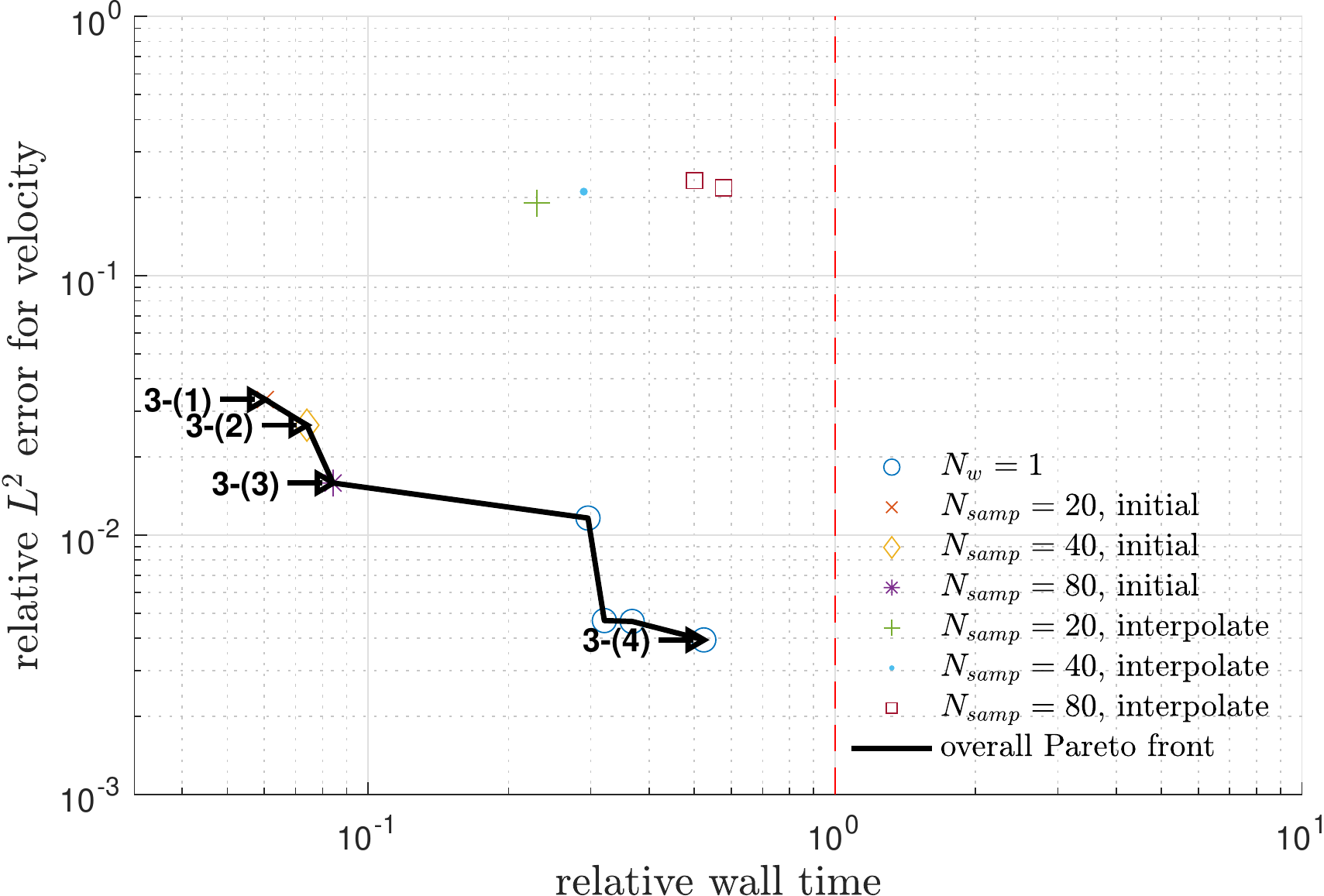} \\ 
\vspace{0.2in}
\includegraphics[width=0.6\linewidth]{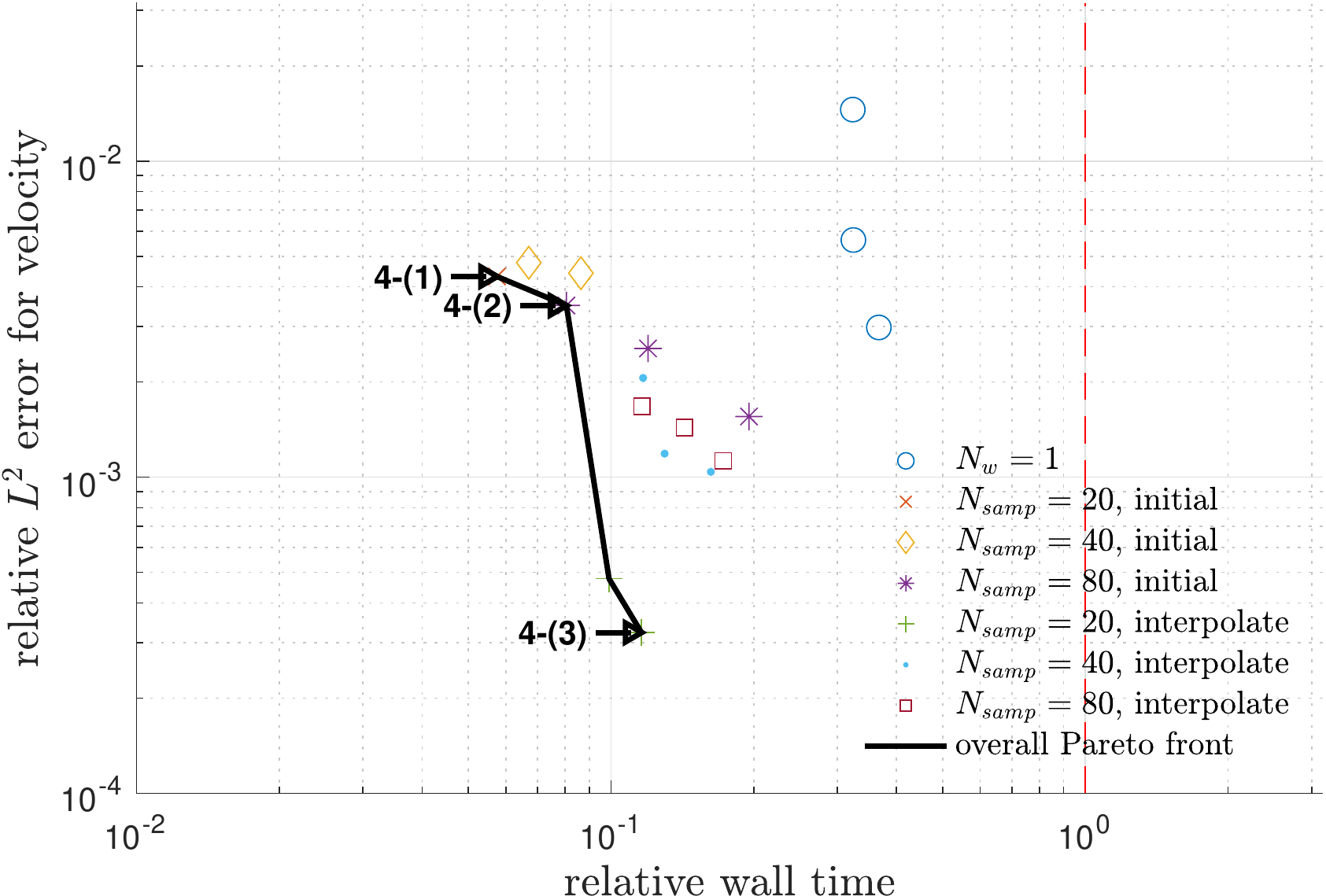} \\
\vspace{0.2in}
\includegraphics[width=0.6\linewidth]{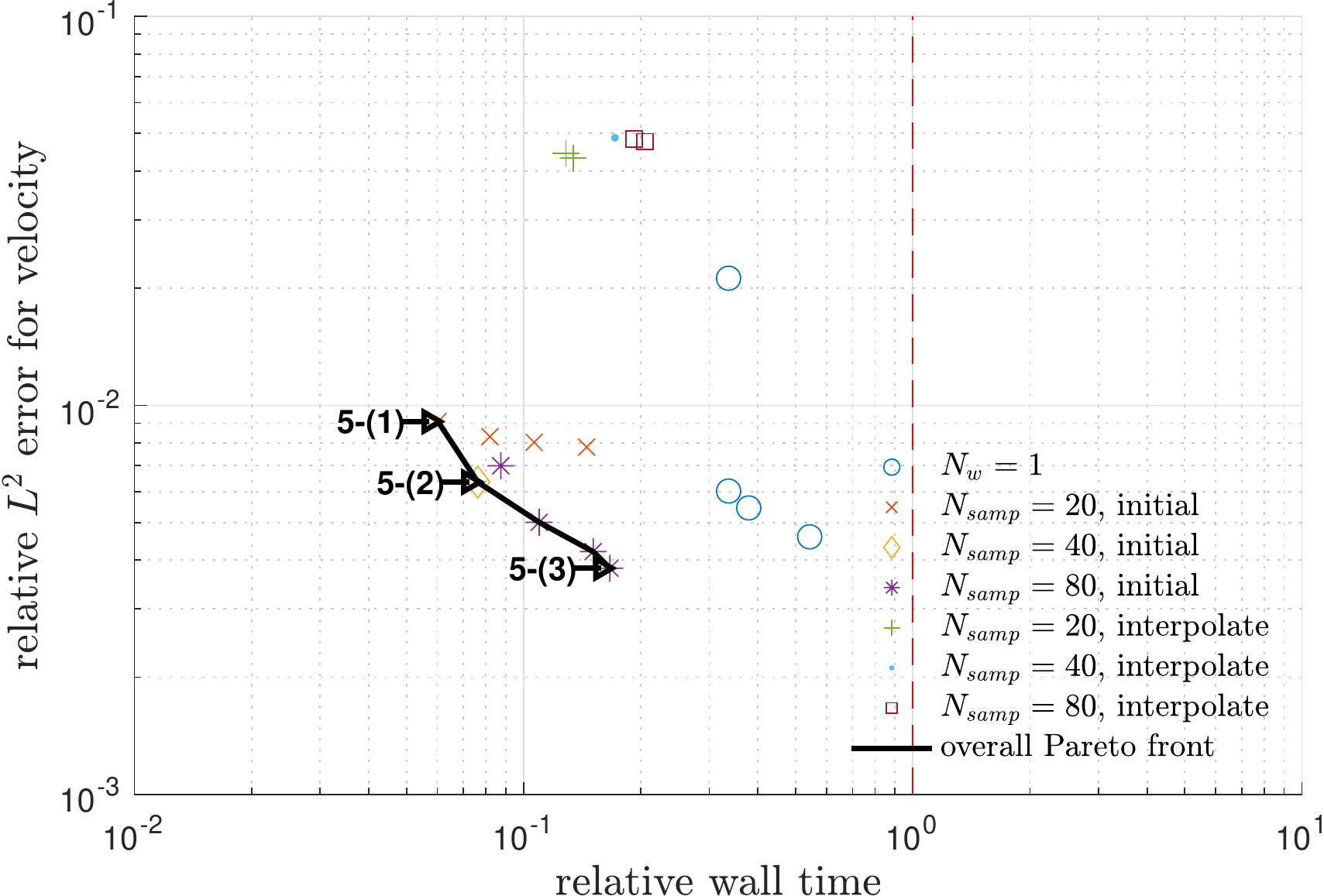} \\
\caption{ROM performance comparison for long-time simulation in the parametric Sedov Blast problem
with the tested problem parameters $\param = 0.7$ (top), $\param = 0.8$ (mid) and $\param = 0.9$ (bottom). 
Relative $L^2$ error for velocity versus relative wall time for varying ROM
  parameters. 
  }
\label{fig:pareto_param_sedov}
\end{figure}

We remark that the tested parameter $\param = 0.8$ is a reproductive case, the
tested parameter $\param = 0.9$ is a predictive case in the interpolation
regime, while the case $\nparam = 0.7$ is predictive in the
extrapolation regime.  For each of the tested problem parameters $\param \in
\{0.7, 0.8, 0.9\}$, we construct a Pareto front in
Figure~\ref{fig:pareto_param_sedov} for each of the cases,   which is
characterized by the ROM user-defined input values that minimize the competing
objectives of relative $L^2$ error for velocity and relative wall time.  An
overall Pareto front that selects the ROM parameters that are Pareto-optimal
across all groups is also illustrated for each problem.
Table~\ref{tab:pareto_param_sedov} reports the ROM user-defined  input values
that yielded the results on the overall Pareto front.  The numerical results
show that for the predictive cases $\param \in \{0.7, 0.9\}$, using the initial state as
offset vectors described in Section~\ref{sec:offset_initial} is Pareto optimal
in any region. It is also Pareto optimal for smaller relative wall time for the
reproductive case $\param = 0.8$. In larger relative wall time, using
the interpolation scheme in Section~\ref{sec:offset_interpolate} is Pareto
optimal for $\param = 0.8$.

\begin{table}[ht!]
\centering
\begin{tabular}{|c||c|c|c|c||c|c|c||c|c|c|}
\hline
Label & 3-(1) & 3-(2) & 3-(3) & 3-(4) & 4-(1) & 4-(2) & 4-(3) & 5-(1) & 5-(2) & 5-(3) \\
\hline
$\nwindow$ & -- & -- & -- & 1 & -- & -- & -- & -- & -- & -- \\
$\ntimestepWindow$ & 20 & 40 & 80 & -- & 20 & 80 & 20 & 20 & 40 & 80 \\ 
$\singularValueThreshold$ & 0.9999 & 0.9999 & 0.9999 & 0.9999 & 0.9999 & 0.9999 & 0.9999 & 0.9999 & 0.9999 & 0.9999 \\
Offset & Initial & Initial & Initial & -- & Initial & Initial & Interpolate & Initial & Initial & Initial \\
$\factorROMforceOneSample$ & 4 & 4 & 4 & 32 & 4 & 4 & 8 & 4 & 4 & 32 \\ 
$\factorROMforceTvSample$ & 4 & 4 & 4 & 32 & 4 & 4 & 8 & 4 & 4 & 32 \\ 
\hline
\end{tabular}
\caption{ROM user-defined input values yielding Pareto-optimal performance for
  long-time simulation in the parametric Sedov Blast problem.
  Figure~\ref{fig:pareto_param_sedov} provides labels.}
\label{tab:pareto_param_sedov}
\end{table}

\section{Conclusion}\label{sec:conclusion} 
In this paper, we develop an efficient reduced order modeling approach for
Lagrangian hydrodynamics simulations. A time-windowing approach is introduced to
construct temporally-local ROM spaces which are small but accurate within a
short period in advection-dominated problems. Different techniques of window
division and construction of offset variables are discussed. Over-sampling
hyper-reduction on the nonlinear terms are introduced to ensure adequate
speed-up. Error bounds for the fully discrete reduced order model with adaptive
time-step control are developed using the continuous-in-time full order model
solution as the reference.  The error bounds are controlled by several
quantities, including the oblique projection error of the solution and the
nonlinear terms onto reduced subspaces, the mismatch of initial condition, an
exponential factor of the time of evaluation, and the maximum time step size.
Numerical examples are shown to verify the capability of our method.  For the
long-time simulation in benchmark advection-dominated problems, the time
windowing reduced order model successfully captures the extreme mesh distortion
exhibited in the Lagrangian full order model, and attains remarkable speed-up of 
25 to 75 times and outstanding solution accuracy with the order of -9
to -4 in the reproductive cases.  As suggested by the numerical results, the
subspace relation is useful for extracting quality nonlinear term bases and
improving the solution accuracy and overall speed-up.  Various offset
approaches in the time windowing reduced order model are compared.
In the reproductive case, using the interpolation scheme is Pareto optimal
in the region of higher solution accuracy, while using the initial states is
Pareto optimal in the region of smaller relative clock time.  In the predictive
cases, using the initial states provides uniformly improved accuracy 
when compared with serial ROM and a speed-up of around 15 times at all the
tested parameter values.

\section*{Acknowledgments}
This work was performed at Lawrence Livermore National Laboratory and partially
funded by LDRD (21-FS-042).  Lawrence Livermore National Laboratory is operated
by Lawrence Livermore National Security, LLC, for the U.S. Department of
Energy, National Nuclear Security Administration under Contract
DE-AC52-07NA27344 and LLNL-JRNL-820660.

\section*{Disclaimer}
This document was prepared as an account of work sponsored by an agency of the
United States government.  Neither the United States government nor Lawrence
Livermore National Security, LLC, nor any of their employees makes any warranty,
expressed or implied, or assumes any legal liability or responsibility for the
accuracy, completeness, or usefulness of any information, apparatus, product, or
process disclosed, or represents that its use would not infringe privately owned
rights.  Reference herein to any specific commercial product, process, or
service by trade name, trademark, manufacturer, or otherwise does not
necessarily constitute or imply its endorsement, recommendation, or favoring by
the United States government or Lawrence Livermore National Security, LLC.  The
views and opinions of authors expressed herein do not necessarily state or
reflect those of the United States government or Lawrence Livermore National
Security, LLC, and shall not be used for advertising or product endorsement
purposes.

\appendix
\section{Command line options of Laghos}\label{sec:appendix} In this section,
we present some examples of the command line options of Laghos simulation for
the purpose of reproducible research.  Due to rapid software development in the
repository, the commands lines are subject to change. However, we try our best
to present the command lines compatible with recent versions of different
dependent softwares and maintain a simple usage of the program. The following
command lines are compatible with the recent commits of 
the \texttt{master} branch of MFEM
\footnote{GitHub page, {\it https://github.com/mfem/mfem}, commit 4fb1a54}, 
the \texttt{master} branch of libROM
\footnote{GitHub page, {\it https://github.com/LLNL/libROM}, commit 143d4f6}, 
and the \texttt{rom} branch of Laghos
\footnote{GitHub page, {\it https://github.com/CEED/Laghos/tree/rom}, commit 934badf}. 
In order to use the ROM capability of Laghos, a user has to 
navigate to the \texttt{rom} subdirectory. 

\subsection{Problem specification}\label{sec:cmd_prob}
First, we present the commands lines for the FOM user-defined input values 
reported in Table~\ref{tab:fom} using the executable \texttt{laghos}. 
These command lines can also be used along with 
command line options for ROM user-defined input values. 
We use the long-time simulation for illustration. 
The command line options for specifying the Gresho vortex problem are 
\begin{verbatim}
./laghos -p 4 -m data/square_gresho.mesh -rs 4 -ok 3 -ot 2 -tf 0.62 -s 7
\end{verbatim}
The command line options for specifying the Sedov Blast problem are 
\begin{verbatim}
./laghos -p 1 -m data/cube01_hex.mesh -pt 211 -tf 0.8 -s 7
\end{verbatim}
The command line options for specifying the Taylor--Green vortex problem are 
\begin{verbatim}
./laghos -p 0 -m data/cube01_hex.mesh -cfl 0.1 -tf 0.25 -s 7
\end{verbatim}
The command line options for specifying the triple--point problem are 
\begin{verbatim}
./laghos -p 3 -m data/box01_hex.mesh -tf 0.8 -s 7
\end{verbatim}

\subsection{Long-time serial ROM simulation}\label{sec:cmd_serial}
Next, we present the command line options for the long-time serial ROM
simulation in the Sedov Blast problem reported in the second column of
Table~\ref{tab:long}.  The options are appended to the command lines for
specifying the Sedov Blast problem in Section~\ref{sec:cmd_prob}.  The other
problems in Table~\ref{tab:long} can be simulated using similar options.  For
computing FOM reference solution, snapshot sampling and basis generation in the
offline phase, one appends 
\begin{verbatim}
-offline -writesol -romsns -ef 0.9999 -sample-stages
\end{verbatim}
For hyper-reduction preprocessing in the online phase, one appends
\begin{verbatim}
-online -romhrprep -romsns -rdimv 225 -rdime 27 -rdimx 29 -sfacv 6 -sface 6
\end{verbatim}
For ROM simulation in the online phase, one appends
\begin{verbatim}
-online -romhr -romsns -rdimv 225 -rdime 27 -rdimx 29 -sfacv 6 -sface 6
\end{verbatim}
Finally, for solution postprocessing and calculating the relative error, one appends 
\begin{verbatim}
-restore -soldiff -romsns -rdimv 225 -rdime 27 -rdimx 29
\end{verbatim}

\subsection{Long-time time windowing ROM simulation}\label{sec:cmd_tw}
Next, we present the command line options for the long-time time windowing ROM
simulation in the Sedov Blast problem reported in the second column of
Table~\ref{tab:tw}.  The options are appended to the command lines for
specifying the Sedov Blast problem in Section~\ref{sec:cmd_prob}.  The other
problems in Table~\ref{tab:tw} can be simulated using similar options.  For
computing FOM reference solution, snapshot sampling and basis generation in the
offline phase, one appends 
\begin{verbatim}
-offline -writesol -romsns -ef 0.9999 -sample-stages -nwinsamp 10 -rostype load
\end{verbatim}
For hyper-reduction preprocessing in the online phase, one appends
\begin{verbatim}
-online -romhrprep -romsns -sfacv 2 -sface 2 -nwin 144 -rostype load
\end{verbatim}
For ROM simulation in the online phase, one appends
\begin{verbatim}
-online -romhr -romsns -sfacv 2 -sface 2 -nwin 144 -rostype load
\end{verbatim}
Finally, for solution postprocessing and calculating the relative error, one
appends 
\begin{verbatim}
-restore -soldiff -romsns -nwin 144 -rostype load
\end{verbatim}

\subsection{Parametric time windowing ROM simulation}\label{sec:cmd_par}
Finally, we present the command line options for the long-time time windowing
ROM simulation in the parametric Sedov Blast problem reported in the Label 4-(1)
of Table~\ref{tab:pareto_param_sedov}.  The options are appended to the command
lines for specifying the Sedov Blast problem in Section~\ref{sec:cmd_prob},
except that for the basis generation a separate executable \texttt{merge} is
used.  The other labels in Table~\ref{tab:pareto_param_sedov} can be simulated
using similar options.  For computing FOM reference solution with $\param=0.8$
and snapshot sampling in the offline phase, one appends 
\begin{verbatim}
-offline -romsns -bef 1.0 -rpar 0 -sample-stages -rostype initial
-offline -romsns -bef 1.2 -rpar 1 -sample-stages -rostype initial
-offline -romsns -bef 0.8 -rpar 2 -sample-stages -rostype initial -writesol
\end{verbatim}
For basis generation in the offline phase of the parametric problem, 
one uses a separate executable \texttt{merge}
\begin{verbatim}
./merge -nset 3 -romsns -ef 0.9999 -nwinsamp 20 -rostype initial
\end{verbatim}
For hyper-reduction preprocessing in the online phase, one appends
\begin{verbatim}
-online -romhrprep -romsns -bef 0.8 -sfacv 4 -sface 4 -nwin 80 -rostype initial
\end{verbatim}
For ROM simulation in the online phase, one appends
\begin{verbatim}
-online -romhr -romsns -bef 0.8 -sfacv 4 -sface 4 -nwin 80 -rostype initial
\end{verbatim}
Finally, for solution postprocessing and calculating the relative error, one appends 
\begin{verbatim}
-restore -soldiff -romsns -nwin 80 -rostype initial
\end{verbatim}

\bibliographystyle{unsrt}
\bibliography{references}

\end{document}

%% file: commands.tex
\newcommand{\bmat}[1]{\begin{bmatrix}#1\end{bmatrix}} 
\newcommand{\pmat}[1]{\begin{pmatrix}#1\end{pmatrix}} 

\newcommand{\dummyMatSymbol}{A}
\newcommand{\dummyElemSymbol}{a}
\newcommand{\dummyVecSymbol}{\dummyElemSymbol}
\newcommand{\dummyTimeSymbol}{s}
\newcommand{\dummyTimeIndex}{j}

\newcommand{\dimensionSymbol}{d}
\newcommand{\solDomainSymbol}{\Omega}
\newcommand{\gradientSymbol}{\nabla}
\newcommand{\energySymbol}{e}
\newcommand{\velocitySymbol}{v}
\newcommand{\positionSymbol}{x}
\newcommand{\densitySymbol}{\rho}
\newcommand{\pressureSymbol}{p}
\newcommand{\stressSymbol}{\sigma}
\newcommand{\artificialStressSymbol}{\stressSymbol_a}
\newcommand{\adiabaticIndexSymbol}{\gamma}
\newcommand{\normalSymbol}{n}
\newcommand{\neumannSymbol}{g}
\newcommand{\identitySymbol}{I}
\newcommand{\zeroSymbol}{0}
\newcommand{\initialSymbol}[1]{\widetilde{#1}}
\newcommand{\initialDomain}{\initialSymbol{\solDomainSymbol}}
\newcommand{\initialPosition}{\initialSymbol{\positionSymbol}}
\newcommand{\meshSize}{h}
\newcommand{\meshLevel}{m}
\newcommand{\feOrder}{k}

\newcommand{\snapshotSymb}{E}
\newcommand{\snapshotSymbS}{S}

\newcommand{\subspaceSymbol}{\mathcal{S}}
\newcommand{\timeSymbol}{t}
\newcommand{\jacobianSymbol}{J}
\newcommand{\massMatSymbol}{M}
\newcommand{\choleskyMatSymbol}{L}
\newcommand{\oneSymbol}{1}
\newcommand{\forceSymbol}{F}

\newcommand{\kinematicSymbol}{V}
\newcommand{\thermodynamicSymbol}{E}
\newcommand{\feSpace}[1]{\mathcal{#1}}
\newcommand{\domain}[1]{\mathsf{#1}}
\newcommand{\paramSymbol}{\mu}
\newcommand{\mapto}{\rightarrow}
\newcommand{\sizeFOMsymbol}{N}
\newcommand{\sizeROMsymbol}{n}
\newcommand{\sizeSnapshot}{N_{\snapshotSymb}}
\newcommand{\sizeSnapshotS}{N_{\snapshotSymbS}}
\newcommand{\nsmall}{\sizeROMsymbol}
\newcommand{\Span}[1]{\text{Span}\{#1\}}
\newcommand{\RR}[1]{\mathbb{R}^{#1}}

\newcommand{\timeIndex}{n}
\newcommand{\basisIndex}{i}
\newcommand{\paramIndex}{k}
\newcommand{\windowIndex}{w}
\newcommand{\paramDomainSymbol}{D}
\newcommand{\finalSymbol}{f}

\newcommand{\ROMBasisSymbol}{\mathbf{\Phi}}
\newcommand{\ROMBasisVecSymbol}{\mathbf{\phi}}
\newcommand{\ROMSamplingMatSymbol}{\mathbf{Z}}
\newcommand{\ROMSolutionSymbol}[1]{\widehat{#1}}
\newcommand{\ProjectionSymbol}{\mathcal{P}}
\newcommand{\IdentitySymbol}{\mathbb{I}}
\newcommand{\continuousErrorSymbol}{\varepsilon}
\newcommand{\discreteErrorSymbol}{\delta}
\newcommand{\relErrorSymbol}{\mathcal{\epsilon}}
\newcommand{\differenceSymbol}{\theta}
\newcommand{\initialErrorSymbol}{\zeta}
\newcommand{\residualSymbol}{\mathcal{R}}
\newcommand{\norm}[1]{\left\lVert #1 \right\rVert}
\newcommand{\fullNorm}[1]{{\left\vert\kern-0.25ex\left\vert\kern-0.25ex\left\vert #1 
    \right\vert\kern-0.25ex\right\vert\kern-0.25ex\right\vert}}
\newcommand{\stateSymbol}{w}
\newcommand{\rkcoeffSymbol}{Y}
\newcommand{\rkcoeffk}[1]{\mathsf{\rkcoeffSymbol}_{#1}}
\newcommand{\rkcoeffI}{\rkcoeffk{1}}
\newcommand{\rkcoeffII}{\rkcoeffk{2}}
\newcommand{\rkcoeffIII}{\rkcoeffk{3}}
\newcommand{\rkcoeffIV}{\rkcoeffk{4}}
\newcommand{\rkErrorCoeffSymbol}{a}
\newcommand{\CFLconst}{\alpha}

\newcommand{\forceOne}{\forceSys^{\mathfrak{1}}}
\newcommand{\forceTv}{\forceSys^{t\velocitySymbol}}
\newcommand{\avgforceTv}{\bar{\forceSys}^{t\velocitySymbol}}
\newcommand{\forceOnek}[1]{\forceOne_{#1}}
\newcommand{\forceTvk}[1]{\forceTv_{#1}}
\newcommand{\avgforceTvk}[1]{\avgforceTv_{#1}}
\newcommand{\forceOneApprox}{\tilde{\forceSys}^{\mathfrak{1}}}
\newcommand{\forceTvApprox}{\tilde{\forceSys}^{t\velocitySymbol}}
\newcommand{\avgforceTvApprox}{\bar{\tilde{\forceSys}}^{t\velocitySymbol}}
\newcommand{\forceOneApproxk}[1]{\forceOneApprox_{#1}}
\newcommand{\forceTvApproxk}[1]{\forceTvApprox_{#1}}
\newcommand{\avgforceTvApproxk}[1]{\avgforceTvApprox_{#1}}
\newcommand{\forceSysApprox}{\tilde{\forceSys}}

\newcommand{\energySubspace}{\subspaceSymbol_{\energySymbol}}
\newcommand{\velocitySubspace}{\subspaceSymbol_{\velocitySymbol}}
\newcommand{\positionSubspace}{\subspaceSymbol_{\positionSymbol}}
\newcommand{\forceOneSubspace}{\subspaceSymbol_{\forceOne}}
\newcommand{\forceTvSubspace}{\subspaceSymbol_{\forceTv}}

\newcommand{\paramDomain}{\domain{\paramDomainSymbol}}

\newcommand{\energy}{\boldsymbol{\energySymbol}}
\newcommand{\velocity}{\boldsymbol{\velocitySymbol}}
\newcommand{\position}{\boldsymbol{\positionSymbol}}
\newcommand{\state}{\boldsymbol{\stateSymbol}}

\newcommand{\velocityApprox}{\tilde{\velocity}}
\newcommand{\energyApprox}{\tilde{\energy}}
\newcommand{\positionApprox}{\tilde{\position}}
\newcommand{\stateApprox}{\tilde{\state}}

\newcommand{\continuousError}{\boldsymbol{\continuousErrorSymbol}}
\newcommand{\discreteError}{\boldsymbol{\discreteErrorSymbol}}
\newcommand{\relError}{\relErrorSymbol}
\newcommand{\difference}{\boldsymbol{\differenceSymbol}}
\newcommand{\initialError}{\boldsymbol{\initialErrorSymbol}}
\newcommand{\residual}{\boldsymbol{\residualSymbol}}

\newcommand{\massMat}{\boldsymbol{\massMatSymbol}}
\newcommand{\choleskyMat}{\boldsymbol{\choleskyMatSymbol}}
\newcommand{\param}{\boldsymbol{\paramSymbol}}
\newcommand{\oneVec}{\boldsymbol{\oneSymbol}}
\newcommand{\kinematicFE}{\feSpace{\kinematicSymbol}}
\newcommand{\thermodynamicFE}{\feSpace{\thermodynamicSymbol}}
\newcommand{\kinematicMassMat}{\massMat_{\kinematicFE}}
\newcommand{\thermodynamicMassMat}{\massMat_{\thermodynamicFE}}
\newcommand{\kinematicCholeskyMat}{\choleskyMat_{\kinematicFE}}
\newcommand{\thermodynamicCholeskyMat}{\choleskyMat_{\thermodynamicFE}}
\newcommand{\jacobian}{\boldsymbol{\jacobianSymbol}}

\newcommand{\forceMat}{\boldsymbol{\forceSymbol}}

\newcommand{\forceSys}{\mathsf{\forceSymbol}}
\newcommand{\sizeKinematicFE}{\sizeFOMsymbol_{\kinematicFE}}
\newcommand{\sizeThermodynamicFE}{\sizeFOMsymbol_{\thermodynamicFE}}
\newcommand{\sizeROMforceOne}{\sizeROMsymbol_{\forceOne}}
\newcommand{\sizeROMforceTv}{\sizeROMsymbol_{\forceTv}}
\newcommand{\samplingOp}[1]{\breve{#1}}
\newcommand{\sizeROMforceOneSample}{\sizeROMsymbol_{\samplingOp{\forceOne}}}
\newcommand{\sizeROMforceTvSample}{\sizeROMsymbol_{\samplingOp{\forceTv}}}
\newcommand{\sizeROMvelocity}{\sizeROMsymbol_{\velocitySymbol}}
\newcommand{\sizeROMenergy}{\sizeROMsymbol_{\energySymbol}}
\newcommand{\sizeROMposition}{\sizeROMsymbol_{\positionSymbol}}
\newcommand{\ratioROMvelocity}{\zeta_{\velocitySymbol}}
\newcommand{\ratioROMenergy}{\zeta_{\energySymbol}}
\newcommand{\ratioROMposition}{\zeta_{\positionSymbol}}
\newcommand{\ratioROMforceOneSample}{\zeta_{\samplingOp{\forceOne}}}
\newcommand{\ratioROMforceTvSample}{\zeta_{\samplingOp{\forceTv}}}
\newcommand{\factorROMforceOneSample}{\lambda_{\samplingOp{\forceOne}}}
\newcommand{\factorROMforceTvSample}{\lambda_{\samplingOp{\forceTv}}}
\newcommand{\sizeWholeFE}{\sizeFOMsymbol}
\newcommand{\sizeWholeROM}{\sizeROMsymbol}
\newcommand{\timek}[1]{\timeSymbol_{#1}}
\newcommand{\finalTime}{\timek{\finalSymbol}}
\newcommand{\ntimestep}{N_{\timeSymbol}}
\newcommand{\timestep}{\Delta\timeSymbol}
\newcommand{\timestepk}[1]{\timestep_{#1}}
\newcommand{\timestepestimateSymbol}{\tau}
\newcommand{\timestepestimatek}[1]{\timestepestimateSymbol_{#1}}

\newcommand{\avgvelocityt}[1]{\bar{\velocity}_{#1}}
\newcommand{\statet}[1]{\state_{#1}}
\newcommand{\energyt}[1]{\energy_{#1}}
\newcommand{\velocityt}[1]{\velocity_{#1}}
\newcommand{\positiont}[1]{\position_{#1}}

\newcommand{\avgvelocityApproxt}[1]{\bar{\velocityApprox}_{#1}}
\newcommand{\stateApproxt}[1]{\stateApprox_{#1}}
\newcommand{\energyApproxt}[1]{\energyApprox_{#1}}
\newcommand{\velocityApproxt}[1]{\velocityApprox_{#1}}
\newcommand{\positionApproxt}[1]{\positionApprox_{#1}}
\newcommand{\ROMenergyt}[1]{\ROMenergy_{#1}}
\newcommand{\ROMvelocityt}[1]{\ROMvelocity_{#1}}
\newcommand{\ROMpositiont}[1]{\ROMposition_{#1}}

\newcommand{\nat}[1]{\mathbb{N}(#1)}
\newcommand{\identityMat}[1]{\boldsymbol{\identitySymbol}_{#1}}
\newcommand{\zeroVec}[1]{\boldsymbol{\zeroSymbol}_{#1}}
\newcommand{\velocityBasis}{\ROMBasisSymbol_{\velocitySymbol}}
\newcommand{\energyBasis}{\ROMBasisSymbol_{\energySymbol}}
\newcommand{\positionBasis}{\ROMBasisSymbol_{\positionSymbol}}
\newcommand{\velocityBasisVec}{\ROMBasisVecSymbol_{\velocitySymbol}}
\newcommand{\energyBasisVec}{\ROMBasisVecSymbol_{\energySymbol}}
\newcommand{\positionBasisVec}{\ROMBasisVecSymbol_{\positionSymbol}}
\newcommand{\velocityBasisVeck}[1]{\velocityBasisVec^{#1}}
\newcommand{\energyBasisVeck}[1]{\energyBasisVec^{#1}}
\newcommand{\positionBasisVeck}[1]{\positionBasisVec^{#1}}
\newcommand{\forceSysEnergy}{\forceSys_{\energySymbol}}
\newcommand{\forceSysVelocity}{\forceSys_{\velocitySymbol}}
\newcommand{\forceSysPosition}{\forceSys_{\positionSymbol}}
\newcommand{\forceTvBasis}{\ROMBasisSymbol_{\forceTv}}
\newcommand{\forceOneBasis}{\ROMBasisSymbol_{\forceOne}}
\newcommand{\forceTvBasisVec}{\ROMBasisVecSymbol_{\forceTv}}
\newcommand{\forceOneBasisVec}{\ROMBasisVecSymbol_{\forceOne}}
\newcommand{\forceTvBasisVeck}[1]{\forceTvBasisVec^{#1}}
\newcommand{\forceOneBasisVeck}[1]{\forceOneBasisVec^{#1}}
\newcommand{\forceTvSamplingMat}{\ROMSamplingMatSymbol_{\forceTv}}
\newcommand{\forceOneSamplingMat}{\ROMSamplingMatSymbol_{\forceOne}}

\newcommand{\ROMstate}{\ROMSolutionSymbol{\state}}
\newcommand{\ROMenergy}{\ROMSolutionSymbol{\energy}}
\newcommand{\ROMvelocity}{\ROMSolutionSymbol{\velocity}}
\newcommand{\ROMposition}{\ROMSolutionSymbol{\position}}
\newcommand{\ROMforceOne}{\ROMSolutionSymbol{\forceOne}}
\newcommand{\ROMforceTv}{\ROMSolutionSymbol{\forceTv}}
\newcommand{\ROMforceSys}{\ROMSolutionSymbol{\forceSys}}
\newcommand{\ROMforceSysEnergy}{\ROMforceSys_{\energySymbol}}
\newcommand{\ROMforceSysVelocity}{\ROMforceSys_{\velocitySymbol}}
\newcommand{\ROMforceSysPosition}{\ROMforceSys_{\positionSymbol}}
\newcommand{\offsetSymbol}{\text{os}}
\newcommand{\energyOS}{\energy_{\offsetSymbol}}
\newcommand{\velocityOS}{\velocity_{\offsetSymbol}}
\newcommand{\positionOS}{\position_{\offsetSymbol}}

\newcommand{\ntimestepROM}{\widetilde{N}_{\timeSymbol}}
\newcommand{\ROMReducedMassMat}{\widehat{\massMat}}
\newcommand{\ROMKinematicMassMat}{\ROMReducedMassMat_{\kinematicFE}}
\newcommand{\ROMThermodynamicMassMat}{\ROMReducedMassMat_{\thermodynamicFE}}
\newcommand{\forceOneObliqProjMat}{\ProjectionSymbol_{\forceOne}}
\newcommand{\forceTvObliqProjMat}{\ProjectionSymbol_{\forceTv}}
\newcommand{\velocityIdentity}{\IdentitySymbol_\kinematicFE}
\newcommand{\energyIdentity}{\IdentitySymbol_\thermodynamicFE}
\newcommand{\positionIdentity}{\velocityIdentity}

\newcommand{\windowk}[1]{T_{#1}}
\newcommand{\nwindow}{N_w}
\newcommand{\energySubspaceWindow}[1]{\subspaceSymbol_{\energySymbol}^{#1}}
\newcommand{\velocitySubspaceWindow}[1]{\subspaceSymbol_{\velocitySymbol}^{#1}}
\newcommand{\positionSubspaceWindow}[1]{\subspaceSymbol_{\positionSymbol}^{#1}}
\newcommand{\forceOneSubspaceWindow}[1]{\subspaceSymbol_{\forceOne}^{#1}}
\newcommand{\forceTvSubspaceWindow}[1]{\subspaceSymbol_{\forceTv}^{#1}}
\newcommand{\sizeROMvelocityWindow}[1]{\sizeROMsymbol_{\velocitySymbol}^{#1}}
\newcommand{\sizeROMenergyWindow}[1]{\sizeROMsymbol_{\energySymbol}^{#1}}
\newcommand{\sizeROMpositionWindow}[1]{\sizeROMsymbol_{\positionSymbol}^{#1}}
\newcommand{\sizeROMforceOneWindow}[1]{\sizeROMsymbol_{\forceOne}^{#1}}
\newcommand{\sizeROMforceTvWindow}[1]{\sizeROMsymbol_{\forceTv}^{#1}}
\newcommand{\sizeROMforceOneSampleWindow}[1]{\sizeROMsymbol_{\samplingOp{\forceOne}}^{#1}}
\newcommand{\sizeROMforceTvSampleWindow}[1]{\sizeROMsymbol_{\samplingOp{\forceTv}}^{#1}}
\newcommand{\energyOSWindow}[1]{\energy_{\offsetSymbol}^{#1}}
\newcommand{\velocityOSWindow}[1]{\velocity_{\offsetSymbol}^{#1}}
\newcommand{\positionOSWindow}[1]{\position_{\offsetSymbol}^{#1}}
\newcommand{\velocityBasisWindow}[1]{\ROMBasisSymbol_{\velocitySymbol}^{#1}}
\newcommand{\energyBasisWindow}[1]{\ROMBasisSymbol_{\energySymbol}^{#1}}
\newcommand{\positionBasisWindow}[1]{\ROMBasisSymbol_{\positionSymbol}^{#1}}
\newcommand{\forceOneBasisWindow}[1]{\ROMBasisSymbol_{\forceOne}^{#1}}
\newcommand{\forceTvBasisWindow}[1]{\ROMBasisSymbol_{\forceTv}^{#1}}
\newcommand{\ROMenergyWindow}[1]{\ROMenergy^{#1}}
\newcommand{\ROMvelocityWindow}[1]{\ROMvelocity^{#1}}
\newcommand{\ROMpositionWindow}[1]{\ROMposition^{#1}}
\newcommand{\ROMforceOneWindow}[1]{\ROMforceOne^{#1}}
\newcommand{\ROMforceTvWindow}[1]{\ROMforceTv^{#1}}
\newcommand{\ROMKinematicMassMatWindow}[1]{\ROMKinematicMassMat^{#1}}
\newcommand{\ROMThermodynamicMassMatWindow}[1]{\ROMThermodynamicMassMat^{#1}}
\newcommand{\forceOneObliqProjMatWindow}[1]{\forceOneObliqProjMat^{#1}}
\newcommand{\forceTvObliqProjMatWindow}[1]{\forceOneObliqProjMat^{#1}}
\newcommand{\forceTvSamplingMatWindow}[1]{\forceTvSamplingMat^{#1}}
\newcommand{\forceOneSamplingMatWindow}[1]{\forceOneSamplingMat^{#1}}
\newcommand{\ntimestepWindow}{N_{\text{sample}}}
\newcommand{\timeIndexWindow}[1]{\timeIndex_{#1}}
\newcommand{\paramMetric}{d}
\newcommand{\interpolationWeightSymbol}{q}
\newcommand{\interpolationPowerSymbol}{r}

\newcommand{\continuousVelocityError}{\continuousError_\velocitySymbol}
\newcommand{\continuousEnergyError}{\continuousError_\energySymbol}
\newcommand{\continuousPositionError}{\continuousError_\positionSymbol}
\newcommand{\discreteVelocityErrort}[1]{\discreteError_{\velocitySymbol,#1}}
\newcommand{\discreteEnergyErrort}[1]{\discreteError_{\energySymbol,#1}}
\newcommand{\discretePositionErrort}[1]{\discreteError_{\positionSymbol,#1}}
\newcommand{\relerrorVelocityt}[1]{\relError_{\velocitySymbol,#1}}
\newcommand{\relerrorEnergyt}[1]{\relError_{\energySymbol,#1}}
\newcommand{\relerrorPositiont}[1]{\relError_{\positionSymbol,#1}}
\newcommand{\velocityDifference}{\difference_\velocitySymbol}
\newcommand{\energyDifference}{\difference_\energySymbol}
\newcommand{\positionDifference}{\difference_\positionSymbol}
\newcommand{\velocityInitialError}{\initialError_\velocitySymbol}
\newcommand{\energyInitialError}{\initialError_\energySymbol}
\newcommand{\positionInitialError}{\initialError_\positionSymbol}
\newcommand{\stateResidual}[1]{\residual_\stateSymbol^{(#1)}}
\newcommand{\velocityResidual}[1]{\residual_\velocitySymbol^{(#1)}}
\newcommand{\energyResidual}[1]{\residual_\energySymbol^{(#1)}}
\newcommand{\positionResidual}[1]{\residual_\positionSymbol^{(#1)}}
\newcommand{\forceOneResidual}[1]{\residual_{\forceOne}^{(#1)}}
\newcommand{\forceTvResidual}[1]{\residual_{\forceTv}^{(#1)}}
\newcommand{\velocityInducedNorm}[1]{\norm{#1}_{\kinematicMassMat}}
\newcommand{\energyInducedNorm}[1]{\norm{#1}_{\thermodynamicMassMat}}

\newcommand{\euclideanNorm}[1]{\norm{#1}_2}
\newcommand{\snapshots}{\boldsymbol \snapshotSymb}
\newcommand{\rkStage}{r}
\newcommand{\rkErrorCoeff}{\mathbf{\rkErrorCoeffSymbol}}
\newcommand{\rkVelocityErrorCoeff}[1]{\rkErrorCoeff^{(#1)}_{\velocitySymbol,\timeIndex}}
\newcommand{\rkEnergyErrorCoeff}[1]{\rkErrorCoeff^{(#1)}_{\energySymbol,\timeIndex}}
\newcommand{\rkPositionErrorCoeff}[1]{\rkErrorCoeff^{(#1)}_{\positionSymbol,\timeIndex}}

\newcommand{\jacobianState}{\jacobian_{\stateSymbol}}

\newcommand{\dummyMat}{\boldsymbol{\dummyMatSymbol}}
\newcommand{\dummyVec}{\boldsymbol{\dummyVecSymbol}}
\newcommand{\leftSingularMat}{\boldsymbol{U}}
\newcommand{\leftSingularVec}{\boldsymbol{u}}
\newcommand{\leftSingularVeck}[1]{\leftSingularVec_{#1}}
\newcommand{\rightSingularMat}{\boldsymbol{V}}
\newcommand{\singularValueMat}{\boldsymbol{\Sigma}}
\newcommand{\singularValue}{\sigma}
\newcommand{\singularValueThreshold}{\epsilon_{\singularValue}}
\newcommand{\nparam}{\nsmall_\paramSymbol}

\newcommand{\conditionnumber}{\kappa}
\newcommand{\unitvecSymbol}{\mathfrak{e}}
\newcommand{\unitveck}[1]{\unitvecSymbol_{#1}}
\newcommand{\tuningparam}{\eta}
\DeclareMathOperator*{\argmin}{arg\,min}
\newcommand{\multiIndexSymbol}{I}
\newcommand{\multiIndex}{\boldsymbol{\multiIndexSymbol}}

%% file: pap_rev.bbl
\begin{thebibliography}{10}

\bibitem{wang2007large}
Shun Wang, Eric~de Sturler, and Glaucio~H Paulino.
\newblock Large-scale topology optimization using preconditioned {K}rylov
  subspace methods with recycling.
\newblock {\em International journal for numerical methods in engineering},
  69(12):2441--2468, 2007.

\bibitem{de2020three}
Miguel A~Salazar de~Troya and Daniel~A Tortorelli.
\newblock Three-dimensional adaptive mesh refinement in stress-constrained
  topology optimization.
\newblock {\em Structural and Multidisciplinary Optimization},
  62(5):2467--2479, 2020.

\bibitem{de2018adaptive}
Miguel A~Salazar De~Troya and Daniel~A Tortorelli.
\newblock Adaptive mesh refinement in stress-constrained topology optimization.
\newblock {\em Structural and Multidisciplinary Optimization},
  58(6):2369--2386, 2018.

\bibitem{white2020dual}
Daniel~A White, Youngsoo Choi, and Jun Kudo.
\newblock A dual mesh method with adaptivity for stress-constrained topology
  optimization.
\newblock {\em Structural and Multidisciplinary Optimization}, 61(2):749--762,
  2020.

\bibitem{choi2015practical}
Youngsoo Choi, Charbel Farhat, Walter Murray, and Michael Saunders.
\newblock A practical factorization of a {S}chur complement for
  {PDE}-constrained distributed optimal control.
\newblock {\em Journal of Scientific Computing}, 65(2):576--597, 2015.

\bibitem{choi2012simultaneous}
Youngsoo Choi.
\newblock {\em Simultaneous analysis and design in {PDE}-constrained
  optimization}.
\newblock PhD thesis, Stanford University, 2012.

\bibitem{smith2013uncertainty}
Ralph~C Smith.
\newblock {\em Uncertainty quantification: theory, implementation, and
  applications}, volume~12.
\newblock Siam, 2013.

\bibitem{biegler2011large}
Lorenz Biegler, George Biros, Omar Ghattas, Matthias Heinkenschloss, David
  Keyes, Bani Mallick, Luis Tenorio, Bart van Bloemen~Waanders, Karen Willcox,
  and Youssef Marzouk.
\newblock {\em Large-scale inverse problems and quantification of uncertainty},
  volume 712.
\newblock John Wiley \& Sons, 2011.

\bibitem{galbally2010non}
David Galbally, Krzysztof Fidkowski, Karen Willcox, and Omar Ghattas.
\newblock Non-linear model reduction for uncertainty quantification in
  large-scale inverse problems.
\newblock {\em International journal for numerical methods in engineering},
  81(12):1581--1608, 2010.

\bibitem{hoang2020domain}
Chi Hoang, Youngsoo Choi, and Kevin Carlberg.
\newblock Domain-decomposition least-squares {P}etrov-{G}alerkin ({DD-LSPG})
  nonlinear model reduction.
\newblock {\em arXiv preprint arXiv:2007.11835}, 2020.

\bibitem{fritzen2018algorithmic}
Felix Fritzen, Bernard Haasdonk, David Ryckelynck, and Sebastian Sch{\"o}ps.
\newblock An algorithmic comparison of the hyper-reduction and the discrete
  empirical interpolation method for a nonlinear thermal problem.
\newblock {\em Mathematical and computational applications}, 23(1):8, 2018.

\bibitem{choi2019space}
Youngsoo Choi and Kevin Carlberg.
\newblock Space--time least-squares {P}etrov--{G}alerkin projection for
  nonlinear model reduction.
\newblock {\em SIAM Journal on Scientific Computing}, 41(1):A26--A58, 2019.

\bibitem{choi2020sns}
Youngsoo Choi, Deshawn Coombs, and Robert Anderson.
\newblock {SNS}: a solution-based nonlinear subspace method for time-dependent
  model order reduction.
\newblock {\em SIAM Journal on Scientific Computing}, 42(2):A1116--A1146, 2020.

\bibitem{carlberg2018conservative}
Kevin Carlberg, Youngsoo Choi, and Syuzanna Sargsyan.
\newblock Conservative model reduction for finite-volume models.
\newblock {\em Journal of Computational Physics}, 371:280--314, 2018.

\bibitem{mojgani2017lagrangian}
Rambod Mojgani and Maciej Balajewicz.
\newblock {L}agrangian basis method for dimensionality reduction of convection
  dominated nonlinear flows.
\newblock {\em arXiv preprint arXiv:1701.04343}, 2017.

\bibitem{kim2020efficientII}
Youngkyu Kim, Karen~May Wang, and Youngsoo Choi.
\newblock Efficient space-time reduced order model for linear dynamical systems
  in {P}ython using less than 120 lines of code.
\newblock {\em arXiv preprint arXiv:2011.10648}, 2020.

\bibitem{xiao2014non}
Dunhui Xiao, Fangxin Fang, Andrew~G Buchan, Christopher~C Pain, Ionel~Michael
  Navon, Juan Du, and G~Hu.
\newblock Non-linear model reduction for the navier--{S}tokes equations using
  residual deim method.
\newblock {\em Journal of Computational Physics}, 263:1--18, 2014.

\bibitem{burkardt2006pod}
John Burkardt, Max Gunzburger, and Hyung-Chun Lee.
\newblock {POD} and {CVT}-based reduced-order modeling of {N}avier--{S}tokes
  flows.
\newblock {\em Computer methods in applied mechanics and engineering},
  196(1-3):337--355, 2006.

\bibitem{amsallem2015design}
David Amsallem, Matthew Zahr, Youngsoo Choi, and Charbel Farhat.
\newblock Design optimization using hyper-reduced-order models.
\newblock {\em Structural and Multidisciplinary Optimization}, 51(4):919--940,
  2015.

\bibitem{choi2020gradient}
Youngsoo Choi, Gabriele Boncoraglio, Spenser Anderson, David Amsallem, and
  Charbel Farhat.
\newblock Gradient-based constrained optimization using a database of linear
  reduced-order models.
\newblock {\em Journal of Computational Physics}, 423:109787, 2020.

\bibitem{choi2019accelerating}
Youngsoo Choi, Geoffrey Oxberry, Daniel White, and Trenton Kirchdoerfer.
\newblock Accelerating design optimization using reduced order models.
\newblock {\em arXiv preprint arXiv:1909.11320}, 2019.

\bibitem{mcbane2020component}
Sean McBane and Youngsoo Choi.
\newblock Component-wise reduced order model lattice-type structure design.
\newblock {\em arXiv preprint arXiv:2010.10770}, 2020.

\bibitem{ghasemi2015localized}
Mohamadreza Ghasemi and Eduardo Gildin.
\newblock Localized model reduction in porous media flow.
\newblock {\em IFAC-PapersOnLine}, 48(6):242--247, 2015.

\bibitem{jiang2019implementation}
Rui Jiang and Louis~J Durlofsky.
\newblock Implementation and detailed assessment of a {GNAT} reduced-order
  model for subsurface flow simulation.
\newblock {\em Journal of Computational Physics}, 379:192--213, 2019.

\bibitem{yang2016fast}
Yanfang Yang, Mohammadreza Ghasemi, Eduardo Gildin, Yalchin Efendiev, Victor
  Calo, et~al.
\newblock Fast multiscale reservoir simulations with pod-deim model reduction.
\newblock {\em SPE Journal}, 21(06):2--141, 2016.

\bibitem{wang2020generalized}
Min Wang, Siu~Wun Cheung, Eric~T. Chung, Maria Vasilyeva, and Yuhe Wang.
\newblock Generalized multiscale multicontinuum model for fractured vuggy
  carbonate reservoirs.
\newblock {\em Journal of Computational and Applied Mathematics}, 366:112370,
  2020.

\bibitem{yang2017efficient}
Huanhuan Yang and Alessandro Veneziani.
\newblock Efficient estimation of cardiac conductivities via {POD-DEIM} model
  order reduction.
\newblock {\em Applied Numerical Mathematics}, 115:180--199, 2017.

\bibitem{fu2018pod}
Hongfei Fu, Hong Wang, and Zhu Wang.
\newblock {POD/DEIM} reduced-order modeling of time-fractional partial
  differential equations with applications in parameter identification.
\newblock {\em Journal of Scientific Computing}, 74(1):220--243, 2018.

\bibitem{zhao2014pod}
Pengfei Zhao, Cai Liu, and Xuan Feng.
\newblock {POD-DEIM} based model order reduction for the spherical shallow
  water equations with {T}urkel-{Z}was finite difference discretization.
\newblock {\em Journal of Applied Mathematics}, 2014, 2014.

\bibitem{cstefuanescu2013pod}
R~{\c{S}}tef{\u{a}}nescu and Ionel~Michael Navon.
\newblock {POD/DEIM} nonlinear model order reduction of an {ADI} implicit
  shallow water equations model.
\newblock {\em Journal of Computational Physics}, 237:95--114, 2013.

\bibitem{choi2021space}
Youngsoo Choi, Peter Brown, Bill Arrighi, Roberti Anderson, and Kevin Huynh.
\newblock Space-time reduced order model for large-scale linear dynamical
  systems with application to {B}oltzmann transport problems.
\newblock {\em Journal of Computational Physics}, 424:109845, 2021.

\bibitem{mordhorst2017pod}
M~Mordhorst, Timm Strecker, D~Wirtz, Thomas Heidlauf, and Oliver R{\"o}hrle.
\newblock {POD-DEIM} reduction of computational {EMG} models.
\newblock {\em Journal of Computational Science}, 19:86--96, 2017.

\bibitem{dimitriu2013application}
Gabriel Dimitriu, Ionel~M Navon, and R{\u{a}}zvan {\c{S}}tef{\u{a}}nescu.
\newblock Application of {POD-DEIM} approach for dimension reduction of a
  diffusive predator-prey system with allee effect.
\newblock In {\em International conference on large-scale scientific
  computing}, pages 373--381. Springer, 2013.

\bibitem{antil2012reduced}
Harbir Antil, Matthias Heinkenschloss, Ronald~HW Hoppe, Christopher Linsenmann,
  and Achim Wixforth.
\newblock Reduced order modeling based shape optimization of surface acoustic
  wave driven microfluidic biochips.
\newblock {\em Mathematics and Computers in Simulation}, 82(10):1986--2003,
  2012.

\bibitem{cheng2016reduced}
Ming-C Cheng.
\newblock A reduced-order representation of the {S}chr{\"o}dinger equation.
\newblock {\em AIP Advances}, 6(9):095121, 2016.

\bibitem{gugercin2004survey}
Serkan Gugercin and Athanasios~C Antoulas.
\newblock A survey of model reduction by balanced truncation and some new
  results.
\newblock {\em International Journal of Control}, 77(8):748--766, 2004.

\bibitem{benner2015survey}
Peter Benner, Serkan Gugercin, and Karen Willcox.
\newblock A survey of projection-based model reduction methods for parametric
  dynamical systems.
\newblock {\em SIAM review}, 57(4):483--531, 2015.

\bibitem{mou2020data}
Changhong Mou, Birgul Koc, Omer San, Leo Rebholz, and Traian Iliescu.
\newblock Data-driven variational multiscale reduced order models.
\newblock {\em arXiv preprint arXiv:2002.06457}, 2020.

\bibitem{parish2017non}
Eric~J Parish and Karthik Duraisamy.
\newblock Non-{M}arkovian closure models for large eddy simulations using the
  {M}ori-{Z}wanzig formalism.
\newblock {\em Physical Review Fluids}, 2(1):014604, 2017.

\bibitem{gadalla2020comparison}
Mahmoud Gadalla, Marta Cianferra, Marco Tezzele, Giovanni Stabile, Andrea Mola,
  and Gianluigi Rozza.
\newblock On the comparison of {LES} data-driven reduced order approaches for
  hydroacoustic analysis.
\newblock {\em Computers \& Fluids}, page 104819, 2020.

\bibitem{bergmann2009enablers}
Michel Bergmann, C-H Bruneau, and Angelo Iollo.
\newblock Enablers for robust {POD} models.
\newblock {\em Journal of Computational Physics}, 228(2):516--538, 2009.

\bibitem{osth2014need}
Jan {\"O}sth, Bernd~R Noack, Sini{\v{s}}a Krajnovi{\'c}, Diogo Barros, and
  Jacques Bor{\'e}e.
\newblock On the need for a nonlinear subscale turbulence term in {POD} models
  as exemplified for a high-{R}eynolds-number flow over an {A}hmed body.
\newblock {\em Journal of Fluid Mechanics}, 747:518--544, 2014.

\bibitem{baiges2015reduced}
Joan Baiges, Ramon Codina, and Sergio Idelsohn.
\newblock Reduced-order subscales for {POD} models.
\newblock {\em Computer Methods in Applied Mechanics and Engineering},
  291:173--196, 2015.

\bibitem{san2018extreme}
Omer San and Romit Maulik.
\newblock Extreme learning machine for reduced order modeling of turbulent
  geophysical flows.
\newblock {\em Physical Review E}, 97(4):042322, 2018.

\bibitem{lu2020lagrangian}
Hannah Lu and Daniel~M Tartakovsky.
\newblock {L}agrangian dynamic mode decomposition for construction of
  reduced-order models of advection-dominated phenomena.
\newblock {\em Journal of Computational Physics}, 407:109229, 2020.

\bibitem{abgrall2016robust}
R{\'e}mi Abgrall, David Amsallem, and Roxana Crisovan.
\newblock Robust model reduction by ${L}^{1}$-norm minimization and
  approximation via dictionaries: application to nonlinear hyperbolic problems.
\newblock {\em Advanced Modeling and Simulation in Engineering Sciences},
  3(1):1--16, 2016.

\bibitem{carlberg2015adaptive}
Kevin Carlberg.
\newblock Adaptive $h$-refinement for reduced-order models.
\newblock {\em International Journal for Numerical Methods in Engineering},
  102(5):1192--1210, 2015.

\bibitem{parish2019windowed}
Eric~J Parish and Kevin~T Carlberg.
\newblock Windowed least-squares model reduction for dynamical systems.
\newblock {\em arXiv preprint arXiv:1910.11388}, 2019.

\bibitem{shimizu2020windowed}
Yukiko~S Shimizu and Eric~J Parish.
\newblock Windowed space-time least-squares {P}etrov-{G}alerkin method for
  nonlinear model order reduction.
\newblock {\em arXiv preprint arXiv:2012.06073}, 2020.

\bibitem{peherstorfer2018model}
Benjamin Peherstorfer.
\newblock Model reduction for transport-dominated problems via online adaptive
  bases and adaptive sampling.
\newblock {\em arXiv preprint arXiv:1812.02094}, 2018.

\bibitem{constantine2012reduced}
PG~Constantine and G~Iaccarino.
\newblock Reduced order models for parameterized hyperbolic conservations laws
  with shock reconstruction.
\newblock {\em Center for Turbulence Research Annual Brief}, 2012.

\bibitem{taddei2020space}
Tommaso Taddei and Lei Zhang.
\newblock Space-time registration-based model reduction of parameterized
  one-dimensional hyperbolic {PDE}s.
\newblock {\em arXiv preprint arXiv:2004.06693}, 2020.

\bibitem{reiss2018shifted}
Julius Reiss, Philipp Schulze, J{\"o}rn Sesterhenn, and Volker Mehrmann.
\newblock The shifted proper orthogonal decomposition: A mode decomposition for
  multiple transport phenomena.
\newblock {\em SIAM Journal on Scientific Computing}, 40(3):A1322--A1344, 2018.

\bibitem{rim2018transport}
Donsub Rim, Scott Moe, and Randall~J LeVeque.
\newblock Transport reversal for model reduction of hyperbolic partial
  differential equations.
\newblock {\em SIAM/ASA Journal on Uncertainty Quantification}, 6(1):118--150,
  2018.

\bibitem{welper2020transformed}
G~Welper.
\newblock Transformed snapshot interpolation with high resolution transforms.
\newblock {\em SIAM Journal on Scientific Computing}, 42(4):A2037--A2061, 2020.

\bibitem{kirby1992reconstructing}
Michael Kirby and Dieter Armbruster.
\newblock Reconstructing phase space from {PDE} simulations.
\newblock {\em Zeitschrift f{\"u}r angewandte Mathematik und Physik ZAMP},
  43(6):999--1022, 1992.

\bibitem{lee2020model}
Kookjin Lee and Kevin~T Carlberg.
\newblock Model reduction of dynamical systems on nonlinear manifolds using
  deep convolutional autoencoders.
\newblock {\em Journal of Computational Physics}, 404:108973, 2020.

\bibitem{lee2019deep}
Kookjin Lee and Kevin Carlberg.
\newblock Deep conservation: A latent dynamics model for exact satisfaction of
  physical conservation laws.
\newblock {\em arXiv preprint arXiv:1909.09754}, 2019.

\bibitem{kim2020fast}
Youngkyu Kim, Youngsoo Choi, David Widemann, and Tarek Zohdi.
\newblock A fast and accurate physics-informed neural network reduced order
  model with shallow masked autoencoder.
\newblock {\em arXiv preprint arXiv:2009.11990}, 2020.

\bibitem{kim2020efficient}
Youngkyu Kim, Youngsoo Choi, David Widemann, and Tarek Zohdi.
\newblock Efficient nonlinear manifold reduced order model.
\newblock {\em arXiv preprint arXiv:2011.07727}, 2020.

\bibitem{rim2019manifold}
Donsub Rim, Benjamin Peherstorfer, and Kyle~T Mandli.
\newblock Manifold approximations via transported subspaces: Model reduction
  for transport-dominated problems.
\newblock {\em arXiv preprint arXiv:1912.13024}, 2019.

\bibitem{rim2020depth}
Donsub Rim, Luca Venturi, Joan Bruna, and Benjamin Peherstorfer.
\newblock Depth separation for reduced deep networks in nonlinear model
  reduction: Distilling shock waves in nonlinear hyperbolic problems.
\newblock {\em arXiv preprint arXiv:2007.13977}, 2020.

\bibitem{peskin2002immersed}
C.~S. Peskin.
\newblock The immersed boundary method.
\newblock {\em Acta Numerica}, 11:479--517, 2002.

\bibitem{boffi2016discrete}
D.~Boffi and L.~Gastaldi.
\newblock Discrete models for fluid-structure interactions: the finite element
  immersed boundary method.
\newblock {\em Discrete and Continuous Dynamical Systems - Series S},
  9:89--107, 2016.

\bibitem{cheung2018mass}
Siu~Wun Cheung, Eric~T. Chung, and Hyea~Hyun Kim.
\newblock A mass conservative scheme for fluid–structure interaction problems
  by the staggered discontinuous {G}alerkin method.
\newblock {\em Journal of Scientific Computing}, 74:1423--1456, 2018.

\bibitem{yeung2002Lagrangian}
P.~K. Yeung.
\newblock {L}agrangian investigations of turbulence.
\newblock {\em Annual Review of Fluid Mechanics}, 34:115--142, 2002.

\bibitem{luthi2005Lagrangian}
B.~L\"{u}thi, A.~Tsinober, and W.~Kinzelbach.
\newblock {L}agrangian measurement of vorticity dynamics in turbulent flow.
\newblock {\em Journal of Fluid Mechanics}, 528:87--118, 2005.

\bibitem{hirt1974arbitrary}
C.~W. Hirt, A.~A. Amsden, and J.~L. Cook.
\newblock An arbitrary {L}agrangian-{E}ulerian computing method for all flow
  speeds.
\newblock {\em Journal of Computational Physics}, 14:227--253, 1974.

\bibitem{benson1989efficient}
David~J. Benson.
\newblock An efficient, accurate, simple {ALE} method for nonlinear finite
  element programs.
\newblock {\em Computer Methods in Applied Mechanics and Engineering},
  72:305--350, 1989.

\bibitem{guermond2017invariant}
J.~L. Guermond, B.~Popov, L.~Saavedra, and Y.~Yang.
\newblock Invariant domains preserving arbitrary {L}agrangian {E}ulerian
  approximation of hyperbolic systems with continuous finite elements.
\newblock {\em SIAM Journal on Scientific Computing}, 39:A385--A414, 2017.

\bibitem{dobrev2012high}
Veselin~A Dobrev, Tzanio~V Kolev, and Robert~N Rieben.
\newblock High-order curvilinear finite element methods for {L}agrangian
  hydrodynamics.
\newblock {\em SIAM Journal on Scientific Computing}, 34(5):B606--B641, 2012.

\bibitem{harlow1971fluid}
F.H. Harlow and A.A. Amsfen.
\newblock {\em Fluid Dynamics: A LASL Monograph}.
\newblock Tech. rep. LA-4700. Los Alamos Scientific Laboratory, 1971.

\bibitem{chaturantabut2010nonlinear}
Saifon Chaturantabut and Danny~C Sorensen.
\newblock Nonlinear model reduction via discrete empirical interpolation.
\newblock {\em SIAM Journal on Scientific Computing}, 32(5):2737--2764, 2010.

\bibitem{berkooz1993proper}
Gal Berkooz, Philip Holmes, and John~L Lumley.
\newblock The proper orthogonal decomposition in the analysis of turbulent
  flows.
\newblock {\em Annual review of fluid mechanics}, 25(1):539--575, 1993.

\bibitem{hotelling1933analysis}
Harold Hotelling.
\newblock Analysis of a complex of statistical variables into principal
  components.
\newblock {\em Journal of educational psychology}, 24(6):417, 1933.

\bibitem{loeve1955}
Michel Loeve.
\newblock {\em Probability Theory}.
\newblock D. Van Nostrand, New York, 1955.

\bibitem{hinze2005proper}
Michael Hinze and Stefan Volkwein.
\newblock Proper orthogonal decomposition surrogate models for nonlinear
  dynamical systems: Error estimates and suboptimal control.
\newblock In {\em Dimension reduction of large-scale systems}, pages 261--306.
  Springer, 2005.

\bibitem{kunisch2002galerkin}
Karl Kunisch and Stefan Volkwein.
\newblock {G}alerkin proper orthogonal decomposition methods for a general
  equation in fluid dynamics.
\newblock {\em SIAM Journal on Numerical analysis}, 40(2):492--515, 2002.

\bibitem{drmac2016new}
Zlatko Drmac and Serkan Gugercin.
\newblock A new selection operator for the discrete empirical interpolation
  method---improved a priori error bound and extensions.
\newblock {\em SIAM Journal on Scientific Computing}, 38(2):A631--A648, 2016.

\bibitem{drmac2018discrete}
Zlatko Drmac and Arvind~Krishna Saibaba.
\newblock The discrete empirical interpolation method: Canonical structure and
  formulation in weighted inner product spaces.
\newblock {\em SIAM Journal on Matrix Analysis and Applications},
  39(3):1152--1180, 2018.

\bibitem{carlberg2013gnat}
Kevin Carlberg, Charbel Farhat, Julien Cortial, and David Amsallem.
\newblock The {GNAT} method for nonlinear model reduction: Effective
  implementation and application to computational fluid dynamics and turbulent
  flows.
\newblock {\em Journal of Computational Physics}, 242:623--647, 2013.

\bibitem{carlberg2011efficient}
Kevin Carlberg, Charbel Bou‐Mosleh, and Charbel Farhat.
\newblock Efficient non-linear model reduction via a least--squares
  {P}etrov--{G}alerkin projection and compressive tensor approximations.
\newblock {\em International Journal for Numerical Methods in Engineering},
  86:155--181, 2011.

\bibitem{gu1996efficient}
Ming Gu and Stanley~C Eisenstat.
\newblock Efficient algorithms for computing a strong rank-revealing {QR}
  factorization.
\newblock {\em SIAM Journal on Scientific Computing}, 17(4):848--869, 1996.

\bibitem{ahern2013visit}
Sean Ahern, Eric Brugger, Brad Whitlock, Jeremy~S Meredith, Kathleen Biagas,
  Mark~C Miller, and Hank Childs.
\newblock {V}isit: Experiences with sustainable software.
\newblock {\em arXiv preprint arXiv:1309.1796}, 2013.

\bibitem{gresho1990theory}
Philip~M. Gresho and Stevens~T. Chan.
\newblock On the theory of semi‐implicit projection methods for viscous
  incompressible flow and its implementation via a finite element method that
  also introduces a nearly consistent mass matrix. part 2: Implementation.
\newblock {\em International Journal for Numerical Methods in Fluids},
  11:621--659, 1990.

\bibitem{sedov1993similarity}
L.I. Sedov.
\newblock {\em Similarity and Dimensional Methods in Mechanics, Tenth Edition}.
\newblock Taylor \& Francis, 1993.

\bibitem{taylor1937mechanism}
Taylor~Geoffrey Ingram and Green~Albert Edward.
\newblock Mechanism of the production of small eddies from large ones.
\newblock {\em Proceedings of the Royal Society of London A},
  158(895):499--521, 1937.

\bibitem{galera2010two}
S.~Galera, P.-H. Maire, and J.~Breil.
\newblock A two-dimensional unstructured cell-centered multi-material {ALE}
  scheme using {VOF} interface reconstruction.
\newblock {\em Journal of Computational Physics}, 229:621--659, 2010.

\end{thebibliography}
